\documentclass{amsart}
\pdfoutput=1
\usepackage{amssymb}
\usepackage{fullpage}
\usepackage[outline]{contour}
\usepackage{tikz}
\usepackage{pdflscape}
\usepackage{enumitem}
\usepackage{mathdots}
\usepackage[all]{xy}
\usetikzlibrary{decorations.pathreplacing,decorations.markings,decorations.pathreplacing,arrows,shapes}
\usepackage[bookmarks=false]{hyperref}
\hypersetup{colorlinks, linkcolor=blue}

\newtheorem{thm}{Theorem}[section]
\newtheorem{prop}[thm]{Proposition}
\newtheorem{cor}[thm]{Corollary}
\newtheorem{lem}[thm]{Lemma}
\newtheorem{conj}[thm]{Conjecture}

\theoremstyle{definition}
\newtheorem{defn}[thm]{Definition}

\newtheorem{exmp}[thm]{Example}
\newtheorem{notn}[thm]{Notation}

\newtheorem{conv}[thm]{Convention}
\theoremstyle{remark}
\newtheorem{rmk}[thm]{Remark}

\makeatletter
\let\c@equation\c@thm
\makeatother
\numberwithin{equation}{section}

\newcommand{\conf}{\mathrm{Conf}}

\newcommand{\Gr}{\mathrm{Gr}}
\newcommand{\GL}{\mathrm{GL}}
\newcommand{\PGL}{\mathrm{PGL}}
\newcommand{\SL}{\mathrm{SL}}
\newcommand{\mat}{\mathrm{Mat}}
\newcommand{\sgn}{\mathrm{sign}}
\newcommand{\Hom}{\mathrm{Hom}}

\newcommand{\DT}{\mathrm{DT}}
\newcommand{\dt}{\mathrm{dt}}
\newcommand{\trop}{\mathrm{trop}}
\newcommand{\prin}{\mathrm{prin}}

\newcommand{\diag}{\mathrm{diag}}
\newcommand{\Frac}{\mathrm{Frac}}

\newcommand{\tw}{\mathrm{tw}}
\newcommand{\ord}{\mathrm{ord}}
\newcommand{\dGr}{\mathcal{G}r}
\newcommand{\uf}{\mathrm{uf}}
\newcommand{\hc}{\mathrm{hc}}
\newcommand{\cb}{\mathrm{cb}}
\newcommand{\Ad}{\mathrm{Ad}}

\renewcommand{\vec}[1]{\mathbf{#1}}
\newcommand{\inprod}[2]{\left\langle#1,#2\right\rangle}

\tikzset{>=stealth}
\tikzset{->-/.style={decoration={
  markings,
  mark=at position #1 with {\arrow{>}}},postaction={decorate}}, ->-/.default=.5}

\allowdisplaybreaks

\title{Donaldson-Thomas Transformation of Grassmannian}
\author{Daping Weng}
\date{\today}

\begin{document}

\begin{abstract} Let $m$ and $n$ be two integers such that $1<m<n-1$. The configuration space $\conf^\times_n\left(\mathbb{P}^{m-1}\right)$, which is the moduli space of $n$ points in the project space $\mathbb{P}^{m-1}$ satisfying certain general position condition, is closely related to the Grassmannian $\Gr_{m,n}$, and is birationally equivalent to the cluster Poisson variety $\mathcal{X}_{\mathrm{A}_{m-1}\boxtimes \mathrm{A}_{n-m-1}}$.

A Donaldson-Thomas transformation is a special formal automorphism on a cluster Poisson variety which encodes the Donaldson-Thomas invariants of the moduli space of stability conditions on the associated 3d Calabi-Yau category. Existence of a cluster Donaldson-Thomas transformation is part of a sufficient condition given by Gross, Hacking, Keel, and Kontsevich that implies the Fock-Goncharov duality on the cluster ensemble. In this paper we construct the cluster Donaldson-Thomas transformation on the cluster Poisson variety $\mathcal{X}_{\mathrm{A}_{m-1}\boxtimes \mathrm{A}_{n-m-1}}$ and realize it as an isomorphism on the configuration space $\conf^\times_n\left(\mathbb{P}^{m-1}\right)$, proving a conjecture of Goncharov and Shen. By using our result we also prove the Fock-Goncharov duality conjection on the associated cluster ensemble, and give a new proof of the already-proven Zamolodchikov's periodicity conjecture in the case of an $\mathrm{A}_{m-1}\boxtimes \mathrm{A}_{n-m-1}$ quiver.
\end{abstract}

\maketitle

\tableofcontents

\section{Introduction}

Donaldson-Thomas invariants were first introduced by Donaldson and Thomas \cite{DT} as geometric invariants on a Calabi-Yau threefold. Kontsevich and Soibelman \cite{KS} generalized such geometric invariants to the moduli space of stability conditions on 3d Calabi-Yau categories. One family of examples of such 3d Calabi-Yau categories can be constructed from a quiver $Q$, and the associated Donaldson-Thomas invariants can be captured by a formal automorphism on the cluster Poisson variety $\mathcal{X}_Q$ constructed from $Q$, which is known as a Donaldson-Thomas transformation. Keller \cite{KelDT} gave a combinatorial characterization of the Donaldson-Thomas transformation using quiver mutations. Goncharov and Shen \cite{GS} gave an equivalent definition using tropical points of cluster varieties and constructed cluster Donaldson-Thomas transformations on moduli spaces of $G$-local systems on surfaces, which are known to be examples of cluster varieties. 

In this paper we focus on cluster varieties that are closely related to Grassmannians and construct their cluster Donaldson-Thomas transformations, whose existence was conjectured by Goncharov and Shen in \cite{GS}. Our result leads to two corollaries, one is the proof of the Fock-Goncharov duality on these cluster varieties, and the other is a new short proof of the already-proven Zamolodchikov's periodicity conjecture in the case of an $\mathrm{A}\boxtimes \mathrm{A}$ quiver using a reformulation of Keller \cite{Kelperiod}.

\subsection{Main Result}

Let $m$ and $n$ be two integers such that $1<m<n-1$. Let $V$ be an $m$-dimensional vector space. The configuration space $\conf_n^\times\left(\mathbb{P}V\right)$ is defined as
\[
\conf_n^\times\left(\mathbb{P}V\right):=\PGL_m\left\backslash\left\{\left(l_1,\dots, l_n\right)\in \left(\mathbb{P}V\right)^n \ \middle| \ \begin{array}{l}\text{$l_i,\dots, l_{i+m-1}$ are in general position} \\ \text{for each $i$ (indices taken modulo $n$)}\end{array}\right\}\right. ,
\]
where $\PGL_m$ acts diagonally on the product $\left(\mathbb{P}V\right)^n$. Since we have quotient out the $\PGL_m$-action, configuration spaces of $n$ points in projective spaces of the same dimension are actually canonically isomorphic, which allows us to identify them to a single configuration space $\conf_n^\times\left(\mathbb{P}^{m-1}\right)$. This observation is important for our definition of the hyperplane map below.

Let $\left[l_1,\dots, l_n\right]$ be a configuration in $\conf_n^\times\left(\mathbb{P}V\right)$. Due to the general position condition, any consecutive $m-1$ lines $l_i,\dots, l_{i+m-2}$ (indices taken modulo $n$) span a hyperplane in $V$ and hence uniquely define a line $h_{[i,i+m-2]}$ in $V^*$. One can show that the lines $\left(h_{[i,i+m-2]}\right)_{i=1}^n$ also satisfy the general position condition and hence define a configuration in $\conf_n^\times \left(\mathbb{P}V^*\right)$. Since both $\conf_n^\times\left(\mathbb{P}V\right)$ and $\conf_n^\times\left(\mathbb{P}V^*\right)$ are canonically isomorphic to $\conf_n^\times\left(\mathbb{P}^{m-1}\right)$, we can define a map
\begin{align*}
    H:\conf_n^\times\left(\mathbb{P}^{m-1}\right)& \rightarrow \conf_n^\times\left(\mathbb{P}^{m-1}\right)\\
    \left[l_1,\dots, l_n\right]&\mapsto \left[h_{[4-m,2]}, h_{[5-m,3]}, \dots, h_{[3-m,1]}\right].
\end{align*}

It is known that the configuration space $\conf_n^\times\left(\mathbb{P}^{m-1}\right)$ is birationally equivalent to the cluster Poisson variety $\mathcal{X}_{\mathrm{A}_{m-1}\boxtimes \mathrm{A}_{n-m-1}}$ associated the quiver $\mathrm{A}_{m-1}\boxtimes \mathrm{A}_{n-m-1}$ as shown below. 
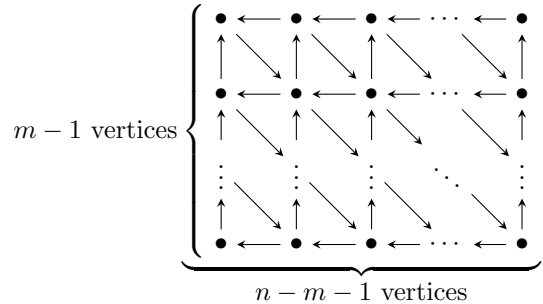
\[
\text{$m-1$ vertices}\underbrace{\left\{\begin{tikzpicture}[baseline=16ex]
\foreach \i in {1,2,3,5}
    {
    \node (\i-2) at (\i,2) [] {$\vdots$};
    \foreach \j in {1,3,4}
        {
        \node (\i-\j) at (\i,\j) [] {$\bullet$};
        }
    }
\foreach \j in {1,3,4}
    {
    \node (4-\j) at (4,\j) [] {$\cdots$};
    }
\node (4-2) at (4,2) [] {$\ddots$};
\foreach \i in {1,2,3,5}
    {
    \draw [->] (\i-1) -- (\i-2);
    \draw [->] (\i-2) -- (\i-3);
    \draw [->] (\i-3) -- (\i-4);
    }
\foreach \j in {1,3,4}
    {
    \draw [->] (5-\j) -- (4-\j);
    \draw [->] (4-\j) -- (3-\j);
    \draw [->] (3-\j) -- (2-\j);
    \draw [->] (2-\j) -- (1-\j);
    }
\foreach \i in {1,...,4}
    {
    \pgfmathtruncatemacro{\k}{\i+1};
    \draw [->] (\i-4) -- (\k-3);
    \draw [->] (\i-3) -- (\k-2);
    \draw [->] (\i-2) -- (\k-1);
    }
\end{tikzpicture}\right.
}_\text{\normalsize{$n-m-1$ vertices}}
\]

Our first main result of this thesis is the following theorem.

\begin{thm} The cluster Donaldson-Thomas transformation $\DT$ exists on the cluster Poisson variety $\mathcal{X}_{\mathrm{A}_{m-1}\boxtimes \mathrm{A}_{n-m-1}}$, and can be realized by the map $H$ defined above in the sense that the following diagram commutes, where the vertical maps are birational equivalence.
\begin{equation}\label{main diagram}
\vcenter{\vbox{\xymatrix{\mathcal{X}_{\mathrm{A}_{m-1}\boxtimes \mathrm{A}_{n-m-1}} \ar[r]^\dt  \ar@{-->}[d]_\cong& \mathcal{X}_{\mathrm{A}_{m-1}\boxtimes \mathrm{A}_{n-m-1}} \ar@{-->}[d]^\cong\\
\conf_n^\times\left(\mathbb{P}^{m-1}\right) \ar[r]_H& \conf_n^\times\left(\mathbb{P}^{m-1}\right) 
}}}
\end{equation}
\end{thm}

We deduce two corollaries from our main result, which we state as follows.

\subsection{Fock-Goncharov Duality}

Let $W$ be an $n$-dimensional vector space and consider the Grassmannian $\Gr_m\left(W^*\right)$. Fix a basis $\left\{e_i\right\}_{i=1}^n$ of $W$; then any $m$-fold exterior product of $e_i$ defines a Pl\"{u}cker coordinate on $\Gr_m\left(W^*\right)$. It is known that the union of the vanishing loci of the Pl\"{u}cker coordinates $\Delta_{e_i\wedge e_{i+1}\wedge \dots \wedge e_{i+m-1}}$ is an anti-canonical divisor on $\Gr_m\left(W^*\right)$, and we denote the complement of such anti-canonical divisor by $\Gr_m^\times\left(W^*\right)$. Note that since a basis $\left\{e_i\right\}$ is already chosen, $\Gr_m^\times\left(W^*\right)$ can be identified with $\Gr_{m,n}^\times:=\Gr_m^\times\left(\mathbb{C}^n\right)$. 

For any point $V^*$ in $\dGr_m^\times\left(W^*\right)$, i.e., an $m$-dimensional subspace $V^*\subset W^*$, we know that the images of the basis vectors $e_i$ under the projection $W\rightarrow V$ are all non-zero and the lines $l_i$ spanned by $e_i$ satisfy the general position condition. Thus we can define a projection map 
\begin{align*}
    \pi:\Gr_{m,n}^\times \cong \Gr_m^\times\left(W^*\right) & \rightarrow \conf_n^\times\left(V\right) \cong \conf_n^\times\left(\mathbb{P}^{m-1}\right)\\
    V& \mapsto [l_1,\dots, l_n].
\end{align*}

It is known that the affine cone over $\Gr_{m,n}^\times$, which we denote by $\dGr_{m,n}^\times$, is birationally equivalent to a cluster $\mathcal{A}$-variety $\mathcal{A}_{m,n}$. Furthermore, Scott \cite{Sco} proved that the birational equivalence induces an isomorphism between their algebras of regular functions. It is also known that the unfrozen cluster $\mathcal{X}$-variety $\mathcal{X}_{m,n}^\uf$ corresponding to $\mathcal{A}_{m,n}$ is precisely the cluster variety $\mathcal{X}_{\mathrm{A}_{m-1}\boxtimes \mathrm{A}_{n-m-1}}$ in our main result. We will prove in Proposition \ref{welldefinedpsi} that the following commutative diagram can be used to construct the birational equivalence $\conf_n^\times\left(\mathbb{P}^{m-1}\right)\dashrightarrow \mathcal{X}_{m,n}^\uf$, where the map $p$ is the canonical map between the two of cluster varieties in an cluster ensemble.
\begin{equation}\label{surjective}
\vcenter{\vbox{\xymatrix{ \dGr_{m,n}^\times \ar@{-->}[r]^\cong \ar[d]_\pi & \mathcal{A}_{m,n} \ar[d]^p \\
\conf_n^\times\left(\mathbb{P}^{m-1}\right) \ar@{-->}[r]_\cong & \mathcal{X}_{m,n}^\uf}}}
\end{equation}

Fock and Goncharov conjectured in \cite{FGensemble} that for any cluster ensemble $\left(\mathcal{A},\mathcal{X}\right)$ associated to a quiver $Q$, the algebra of regular functions $\mathcal{O}\left(\mathcal{A}\right)$ has a canonical basis parametrized by the set of tropical integer points $\mathcal{X}\left(\mathbb{Z}^\trop\right)$, and the algebra of regular functions $\mathcal{O}\left(\mathcal{X}\right)$ has a canonical basis parametrized by the set of tropical integer points $\mathcal{A}\left(\mathbb{Z}^\trop\right)$. This is known as the \emph{Fock-Goncharov duality conjecture}. Gross, Hacking, Keel, and Kontsevich gave a sufficient condition for Fock-Goncharov duality conjecture in \cite{GHKK} (Proposition 8.28), which consists of the following two criteria:
\begin{enumerate}
    \item the cluster complex in $\mathcal{X}_\prin\left(\mathbb{R}^T\right)$ is large enough so that its convex hull contains the whole $\mathcal{X}_\prin\left(\mathbb{R}^T\right)$;
    \item the map $p:\mathcal{A}\rightarrow \mathcal{X}^\uf$ is surjective.
\end{enumerate}

Using this sufficient condition we can prove the following statement.

\begin{cor} The Fock-Goncharov duality conjecture holds for the cluster ensemble $\left(A_{m,n}, \mathcal{X}_{m,n}\right)$. In particular, since $\mathcal{O}\left(\dGr_{m,n}^\times\right)\cong \mathcal{O}\left(\mathcal{A}_{m,n}\right)$ as algebras, it follows that there is a canonical basis of $\mathcal{O}\left(\dGr_{m,n}^\times\right)$ parametrized by the set of tropical integer points $\mathcal{X}_{m,n}\left(\mathbb{Z}^\trop\right)$.
\end{cor}
\begin{proof} Since a cluster Donaldson-Thomas transformation gives rise to a sequence of cluster mutations that goes from one cluster chamber to its opposite chamber in the cluster complex, the existence of cluster Donaldson-Thomas transformation on $\mathcal{X}_\uf$ implies condition (1). Proposition \ref{dominant of psi} shows that the bottom map of the commutative diagram \eqref{surjective} is dominant, which then implies condition (2).
\end{proof}

\subsection{Zamolodchikov's Periodicity Conjecture for \texorpdfstring{$\mathrm{A}_{m-1}\boxtimes \mathrm{A}_{n-m-1}$}{}} The original form of the periodicity conjecture states that for a pair of Dynkin diagrams $\mathrm{D}$ and $\mathrm{D}'$ with vertex set $I$ and $I'$ and incidence matrices $A=\left(a_{ij}\right)$ and $A'=\left(a'_{i'j'}\right)$ respectively (incidence matrix is defined to be $A=2\mathrm{Id}-C$ where $C$ is the Cartan matrix), then the system of variables $\left(Y_{i,j',t}\right)_{t\in \mathbb{Z}}$ defined by the recurrence relation
\[
Y_{i,i',t-1}Y_{i,i',t+1}=\frac{\prod_{j\in I}\left(1+Y_{j,i',t}\right)^{a_{ij}}}{\prod_{j'\in I'} \left(1+Y_{i,j',t}\right)^{a'_{i'j'}}}
\]
is of period (with respect to $t$) dividing $2\left(h+h'\right)$, where $h$ and $h'$ are the Coxeter numbers of $\mathrm{D}$ and $\mathrm{D}'$.

Zamolodchikov's periodicity conjecture has been proved by many people using various methods before, including but not restricted to the following:
\begin{itemize}
    \item $\left(\mathrm{A}_n,\mathrm{A}_1\right)$ by Frenkel and Szenes \cite{FS} and by Giozzi and Tateo \cite{GT};
    \item $\left(\mathrm{D}, \mathrm{A}_1\right)$ by Fomin and Zelevinsky \cite{FZysystem};
    \item $\left(\mathrm{A}_n, \mathrm{A}_m\right)$ by Volkov \cite{Vol}, Szenes \cite{Sze}, and Henriques \cite{Hen};
    \item $\left(\mathrm{D}, \mathrm{D}'\right)$ by Keller \cite{Kelperiod}.
\end{itemize}

When Keller proved Zamolodchikov's periodicity conjecture, he gave a reformulation of the conjecture using Donaldson-Thomas transformations, which reduces the conjecture to show that the order of the Donaldson-Thomas transformation $\DT$ on $\mathcal{X}_{\mathrm{D}\boxtimes \mathrm{D}'}$ divides into $2\left(h+h'\right)$. 

Note that in the case of $\left(\mathrm{A}_{m-1}, \mathrm{A}_{n-m-1}\right)$, the sum of the Coxeter numbers is $h+h'=m+(n-m)=n$; therefore Zamolodchikov's periodicity conjecture in the case of $\left(\mathrm{A}_{m-1}, \mathrm{A}_{n-m-1}\right)$ is equivalent to showing the following statement.

\begin{cor} The cluster Donaldson-Thomas transformation $\DT$ on $\mathcal{X}_{\mathrm{A}_{m-1}\boxtimes \mathrm{A}_{n-m-1}}$ is of order dividing $2n$.
\end{cor}
\begin{proof} We see from the commutative diagram \eqref{main diagram} that it suffices to prove $H^{2n}=\mathrm{Id}$. Note that $H^2:\left[l_1,\dots, l_n\right]\mapsto \left[l_{5-m}, l_{6-m}, \dots, l_{4-m}\right]$, which is just a cyclic shift. Therefore $H^{2n}=\mathrm{Id}$.
\end{proof}

\subsection{Structure on the Paper}. We divide the rest of the paper into five sections.

\begin{itemize}
    \item Section \ref{section2} reviews some basic constructions using bipartite graphs, which are the key combinatorial tool we use to prove our result. The main reference sources for these constructions include \cite{Pos}, \cite{GSVpd}, and \cite{Gon}.
    \item Section \ref{section3} reviews the theory of cluster varieties and their tropicalizations, and state the definition of cluster Donaldson-Thomas transformation. The main reference sourses for these topics include \cite{FZIV}, \cite{FGensemble}, \cite{GS}, and \cite{GHKK}.
    \item Section \ref{section4} describes the connection between the geometric spaces and cluster varieties; the two most important maps are the birational equivalences $\psi:\dGr_{m,n}^\times \dashrightarrow \mathcal{A}_{m,n}$ and $\chi:\mathcal{X}_{m,n}^\uf\dashrightarrow \conf_n^\times\left(\mathbb{P}^{m-1}\right)$. Ideas of these connections are from \cite{Sco}, \cite{Pos}, and \cite{Gon}.
    \item Section \ref{section5} contains the proof of our main result. We first show that the composition $\Xi:=\chi\circ \psi$ maps $\left[l_1,\dots, l_n\right]$ to $\left[h_{[2-m,n]}, h_{[3-m,1]}, \dots, h_{[1-m,n-1]}\right]$, and then we use the map $\Xi$ to show that $\DT:=\psi\circ \chi$ is a cluster transformation and satisfy the degree condition in one version of the cluster Donaldson-Thomas transformation, and finally we induce the $\dt$ version of the cluster Donaldson-Thomas transformation using chiral duality.
    \item Section \ref{section6} discusses how the construction of cluster Donaldson-Thomas transformation presented in this paper is related to the construction of cluster Donaldson-Thomas transformation on double Bruhat cells in semisimple Lie groups.
\end{itemize}

\subsection{Acknowledgements}

The author is deeply grateful to his advisor Alexander Goncharov for suggesting the problem on cluster Donaldson-Thomas transformation and providing insightful comments throughout the problem solving as well as the revision process. The author would like to thank Linhui Shen for his helpful discussion on Donaldson-Thomas transformations. This paper uses many ideas from the work of Fock and Goncharov \cite{FGensemble}, the work of Goncharov and Shen \cite{GS}, the work of Keller \cite{Kelperiod} and \cite{KelDT}, the work of Postnikov \cite{Pos}, and the work of Gekhtman, Shapiro, and Vainshtein \cite{GSVpd}; without the work of all these pioneers in the field this paper would not have been possible, and the author is very thankful to all of them.

\section{Combinatorics of Bipartite Graphs on a Disk}\label{section2}

In this section we introduce bipartite graphs on a disk and some useful combinatorial constructions we can perform on them. We mostly follow the work of Postnikov \cite{Pos}, the work of Gekhtman, Shapiro, and Vainshtein \cite{GSVpd}, and the work of Goncharov \cite{Gon}. 

\subsection{Bipartite Graphs on a Disk}

Fix two positive integers $m$ and $n$ with $1<m<n-1$. Let $\mathbb{D}$ be a disk with marked points on the boundary labeled by $1, \dots, n$ in the clockwise direction (one can think of $\mathbb{D}$ as the closed unit disk with marked points at all the $n$th roots of unity).
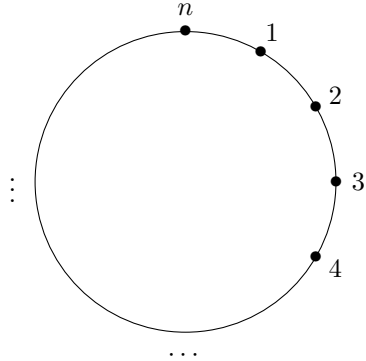
\[
\tikz{
\draw (0,0) circle [radius=2];
\foreach \i in {1,...,4}
    {
    \node at +(90-30*\i:2) [] {$\bullet$};
    \node at +(90-30*\i:2.3) [] {$\i$};
    }
\node at (0,2) [] {$\bullet$};
\node at (0,2.3) [] {$n$};
\node at (0,-2.3) [] {$\cdots$};
\node at (-2.3,0) [] {$\vdots$};
}
\]

\begin{defn} A \emph{bipartite graph} $\Gamma$ on the disk $\mathbb{D}$ is a bipartite graph drawn in the interior of $\mathbb{D}$ with a single \emph{external edge} connecting to each boundary marked point.
\end{defn}

\begin{exmp}\label{bipartite graph example} Below are two examples of bipartite graphs on a disk with parameters $n=4$ and $n=5$ respectively.
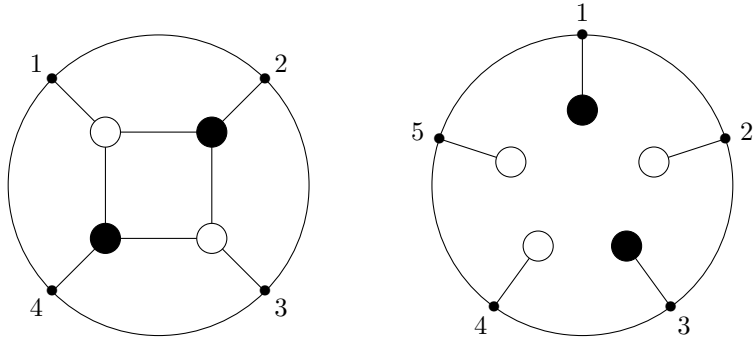
\begin{center}
\begin{tikzpicture}[baseline=0ex]
\draw (0,0) circle [radius=2];
\foreach \i in {1,...,4}
    {
    \coordinate (\i) at +(-135-90*\i:2);
    \node at +(-135-90*\i:2.3) [] {$\i$};
    \node at (\i) [] {$\bullet$};
    }
\foreach \i in {1,2}
    {
    \coordinate (b\i) at +(45-180*\i:1);
    \coordinate (w\i) at +(135-180*\i:1);
    }
\draw (b1) -- (w2) -- (b2) -- (w1) -- cycle;
\draw (1) -- (w2);
\draw (2) -- (b2);
\draw (3) -- (w1);
\draw (4) -- (b1);
\foreach \i in {1,2}
    {
    \draw [fill=white] (w\i) circle [radius=0.2];
    \draw [fill=black] (b\i) circle [radius=0.2];
    }
\end{tikzpicture} \hspace{1cm}
\begin{tikzpicture}[baseline=0ex]
\draw (0,0) circle [radius=2];
\foreach \i in {1,...,5}
    {
    \coordinate (\i) at +(162-72*\i:2);
    \coordinate (i\i) at +(162-72*\i:1);
    \node at +(162-72*\i:2.3) [] {$\i$};
    \node at (\i) [] {$\bullet$};
    \draw (\i) -- (i\i);
    }
\draw [fill] (i1) circle [radius=0.2];
\draw [fill] (i3) circle [radius=0.2];
\draw [fill=white] (i2) circle [radius=0.2];
\draw [fill=white] (i4) circle [radius=0.2];
\draw [fill=white] (i5) circle [radius=0.2];
\end{tikzpicture}
\end{center}
\end{exmp}

Not all bipartite graphs on a disk are useful to us (at least within the scope of this paper); to identify the interesting ones, we need to introduce the notion of \emph{zig-zag strands}. Given a bipartite graph $\Gamma$ on a disk we can draw a set of zig-zag strands, which are directed strands that zig-zags through the edges of $\Gamma$, turning clockwise near white vertices and turning counterclockwise near black vertices.
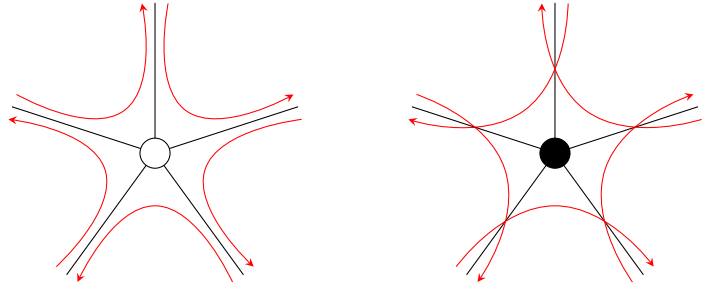
\[
\tikz{
\foreach \i in {0,...,4}
    {
    \draw (0,0) -- +(90-72*\i:2);
    \coordinate (i\i) at +(85-72*\i:2);
    \coordinate (m\i) at +(54-72*\i:0.7);
    \coordinate (o\i) at +(23-72*\i:2);
    \draw [red,->] plot [smooth, tension=1] coordinates {(i\i)(m\i)(o\i)};
    }
\draw [fill=white] (0,0) circle [radius=0.2];
} \quad \quad \quad \quad 
\tikz{
\foreach \i in {0,...,4}
    {
    \draw (0,0) -- +(90-72*\i:2);
    \coordinate (i\i) at +(85-72*\i:2);
    \coordinate (m\i) at +(126-72*\i:0.7);
    \coordinate (o\i) at +(167-72*\i:2);
    \draw [red,->] plot [smooth, tension=1] coordinates {(i\i)(m\i)(o\i)};
    }
\draw [fill] (0,0) circle [radius=0.2];
}
\]
Note that a zig-zag strand either forms an oriented cycle or starts and ends near some boundary marked points. 

\begin{defn}\label{minimal} A bipartite graph $\Gamma$ on a disk is said to be \emph{minimal} if no zig-zag strand form an oriented cycle or intersects itself, and no two zig-zag strands form a parallel bigon. For a minimal bipartite graph $\Gamma$ we label its zig-zag strands as $\zeta_i$ where $i$ is the boundary vertex from which the zig-zag strands starts. A minimal bipartite graph $\Gamma$ on a disk is said to be of \emph{full rank $m$} if the zig-zag strand $\zeta_i$ terminates at the boundary marked point $i+m$ (modulo $n$) for every $i$.
\[
\tikz{
\draw [->, red] plot [smooth, tension=1] coordinates {(-1,1)(1,-1)(1,1)(-1,-1)};
\node at (0,-2) [] {self-intersection};
}\quad \quad \quad \quad \quad \quad \quad \quad 
\tikz{
\draw [->, red] plot [smooth, tension=1] coordinates {(-1.5,-1)(0,0.5)(1.5,-1)};
\draw [->,red] plot [smooth, tension=1] coordinates {(-1.5,1)(0,-0.5)(1.5,1)};
\node at (0,-2) [] {parallel-bigon};
}
\]
\end{defn}

\begin{exmp} As a simple exercise, one can find that the left bipartite graph in Example \ref{bipartite graph example} is minimal and of full rank $4$, whereas the one on the right is not minimal.
\end{exmp}

\begin{conv} For the rest of this paper, by \emph{bipartite graph} we always mean a minimal bipartite graph of full rank $m$ on a disk with $n$ boundary marked points, if not otherwise stated.
\end{conv}

\begin{defn} Let $\Gamma$ be a bipartite graph. A \emph{face} on $\Gamma$ is a connected component of $\mathbb{D}\setminus \Gamma$. A face is said to be a \emph{boundary face} if it contains part of the boundary of the disk $\mathbb{D}$; in particular, we denote the boundary face lying between boundary marked points $i$ and $i+1$ (modulo $n$) by $f_{i\partial}$.
\end{defn}

\begin{defn}\label{dominating set} A face $f$ is said to be \emph{dominated} by a zig-zag strand $\zeta_i$ if a generic point of $f$ lies on the left of $\zeta_i$ with respect to its orientation. The \emph{dominating set} $I(f)$ of a face $f$ is defined to be the set of indices of zig-zag strands that dominate $f$.
\end{defn}

\begin{prop} Let $i$ be a boundary marked point. Then $I\left(f_{i\partial}\right)=\{i-m+1,\dots, i\}$ (modulo $n$).
\end{prop}
\begin{proof} Note that since $f_{i\partial}$ is next to the boundary of the disk, only the zig-zag strands $\zeta_{i-m+1}, \dots, \zeta_i$ can dominate it.
\end{proof}

\begin{notn}\label{bracket} To simplify our notation, for $i,j\in \{1,\dots, n\}$ we define
\[
[i,j]:=\left\{\begin{array}{ll}
    \{i,i+1,\dots, j-1,j\} & \text{if $i\leq j$} \\
    \{i,i+1,\dots, n,1, 2,\dots, j-1,j\} & \text{if $i>j$}.
\end{array}\right.
\]
Therefore the above proposition can be rewritten as $I\left(f_{i\partial}\right)=[i-m+1,i]$ (modulo $n$).
\end{notn}

\begin{prop} Let $\Gamma$ be a biparite graph of full rank $m$. Then $\left|I(f)\right|=m$ for all faces $f$.
\end{prop}
\begin{proof} Note that when we cross an edge to go from one face to an adjacent face, we always cross two oppositely oriented zig-zag strands, and hence the size of the dominating sets are the same. The actual size $m$ can then be easily seen by looking at a boundary face.
\end{proof}

\begin{defn}\label{quiver} From a bipartite graph $\Gamma$ we can also construct a quiver $\vec{i}_\Gamma$ by following the steps below:
\begin{enumerate}
    \item put a vertex for each face of $\Gamma$;
    \item put a counterclockwise cycle with weight $\frac{1}{2}$ around any black vertex and a clockwise cycle with weight $\frac{1}{2}$ around any white vertex (we use dashed arrows to represent half weight);
    \item add up all the weights of the arrows (with orientation considered) between any two faces; the resulting weight is the number of arrows between the two faces.
\end{enumerate}
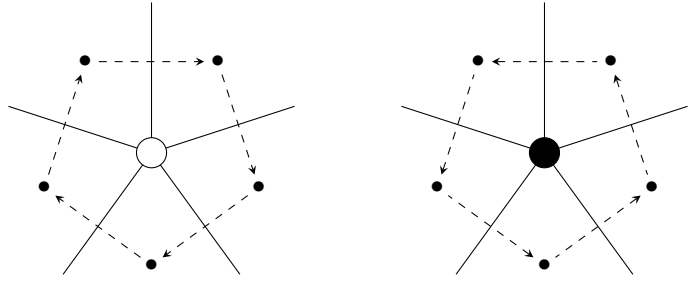
\[
\tikz{
\foreach \i in {0,...,4}
    {
    \draw (0,0) -- +(90-72*\i:2);
    \node at +(54-72*\i:1.5) [] {$\bullet$};
    \draw [dashed, ->] +(48-72*\i:1.4) to +(-12-72*\i:1.4); 
    }
\draw [fill=white] (0,0) circle [radius=0.2];
} \quad \quad \quad \quad 
\tikz{
\foreach \i in {0,...,4}
    {
    \draw (0,0) -- +(90-72*\i:2);
    \node at +(54-72*\i:1.5) [] {$\bullet$};
    \draw [dashed, <-] +(48-72*\i:1.4) to +(-12-72*\i:1.4); 
    }
\draw [fill] (0,0) circle [radius=0.2];
}
\]
\end{defn}

\begin{exmp} Below is an example of a bipartite graph $\Gamma$ with its associated quiver $\vec{i}_\Gamma$.
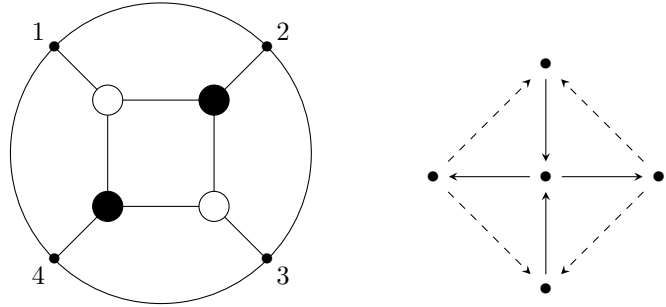
\[ 
\tikz{
\draw (0,0) circle [radius=2];
\foreach \i in {1,...,4}
    {
    \coordinate (\i) at +(-135-90*\i:2);
    \node at +(-135-90*\i:2.3) [] {$\i$};
    \node at (\i) [] {$\bullet$};
    }
\foreach \i in {1,2}
    {
    \coordinate (b\i) at +(45-180*\i:1);
    \coordinate (w\i) at +(135-180*\i:1);
    }
\draw (b1) -- (w2) -- (b2) -- (w1) -- cycle;
\draw (1) -- (w2);
\draw (2) -- (b2);
\draw (3) -- (w1);
\draw (4) -- (b1);
\foreach \i in {1,2}
    {
    \draw [fill=white] (w\i) circle [radius=0.2];
    \draw [fill=black] (b\i) circle [radius=0.2];
    }
}\quad \quad \quad \quad
\tikz{
    \node (0)  at (0,0) [] {$\bullet$};
    \foreach \i in {1,...,4}
        {
        \node (\i) at +(90-90*\i:1.5) [] {$\bullet$};
        }
    \draw [->] (0) -- (1);
    \draw [->] (0) -- (3);
    \draw [->] (2) -- (0);
    \draw [->] (4) -- (0);
    \draw [dashed, ->] (1) -- (4);
    \draw [dashed, ->] (1) -- (2);
    \draw [dashed, ->] (3) -- (4);
    \draw [dashed, ->] (3) -- (2);
}
\]
\end{exmp}

Note that there is no univalent vertex in a bipartite graph as long as we impose the minimality condition in Definition \ref{minimal}. If a black vertex is bivalent, we can delete this black vertex and combine the neighboring two white vertices into one without affecting the dominating sets $I(f)$ or the quiver $\vec{i}_\Gamma$.
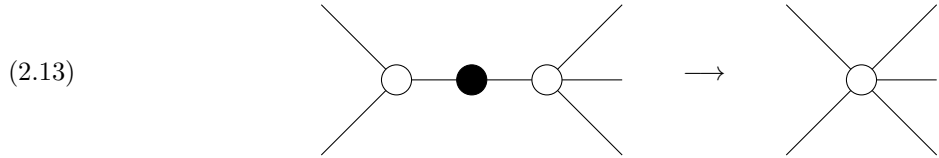
\begin{equation}\label{delete}
\begin{tikzpicture}[baseline=0ex]
\draw (0,0) -- (2,0);
\draw (-1,-1) -- (0,0) -- (-1,1);
\draw (3,-1) -- (2,0) -- (3,1);
\draw (2,0) -- (3,0);
\draw [fill] (1,0) circle [radius=0.2];
\draw [fill=white] (0,0) circle [radius=0.2];
\draw [fill=white] (2,0) circle [radius=0.2];
\end{tikzpicture} \quad \quad \longrightarrow \quad \quad 
\begin{tikzpicture}[baseline=0ex]
\draw (-1,1) -- (1,-1);
\draw (-1,-1) -- (1,1);
\draw (0,0) -- (1,0);
\draw [fill=white] (0,0) circle [radius=0.2];
\end{tikzpicture}
\end{equation}
On the other hand, if a black vertex has valency $v>3$, we can always split it up to $v-2$ black vertices connected by a chain of bivalent white vertices without affecting the dominating sets $I(f)$ or the quiver $\vec{i}_\Gamma$.
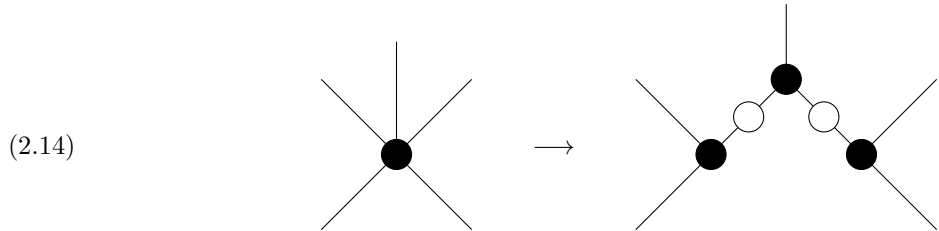
\begin{equation}\label{splitting}
\begin{tikzpicture}[baseline=0ex]
\draw (-1,1) -- (1,-1);
\draw (-1,-1) -- (1,1);
\draw (0,0) -- (0,1.5);
\draw [fill] (0,0) circle [radius=0.2];
\end{tikzpicture}\quad \quad \longrightarrow \quad \quad 
\begin{tikzpicture}[baseline=0ex]
\draw (0,0) -- (1,1) -- (2,0);
\draw (-1,-1) -- (0,0) -- (-1,1);
\draw (3,-1) -- (2,0) -- (3,1);
\draw (1,1) -- (1,2);
\draw [fill=white] (0.5,0.5) circle [radius=0.2];
\draw [fill=white] (1.5,0.5) circle [radius=0.2];
\draw [fill] (1,1) circle [radius=0.2];
\draw [fill] (0,0) circle [radius=0.2];
\draw [fill] (2,0) circle [radius=0.2];
\end{tikzpicture}
\end{equation}
Therefore without loss of generality, we may always assume all black vertices in a bipartite graph are \emph{trivalent}.

By imposing the trivalent condition on black vertices, we get a correspondence between bipartite graphs and a family of Thurston's \emph{minimal triple diagrams} introduced in \cite{Thu}. The process is the following: we first split each boundary marked point $i$ into a pair of neighboring boundary marked points $i'$ and $i$ with $i$ in the clockwise forward direction with respect to $i'$, representing the terminal point of the zig-zag strand $\zeta_{i-m}$ and initial point of $\zeta_i$ respectively (index modulo $n$); we then pull the three zig-zag strands near each black vertex tightly into a triple intersection; lastly we delete the bipartite graph while keeping the zig-zag strands, and the diagram formed by the zig-zag strands is precisely a minimal triple diagram.
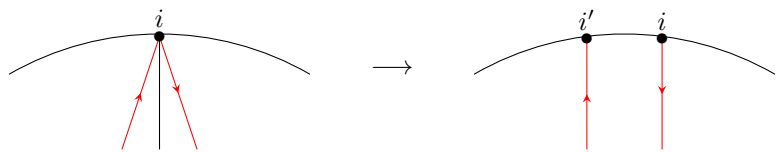
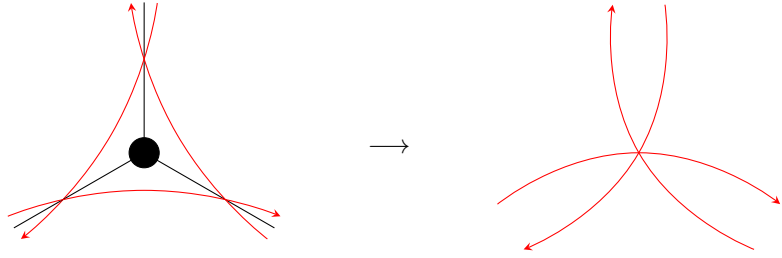
\[
\begin{tikzpicture}[baseline=0ex]
\draw (4,0) arc (60:120:4);
\draw (2,0.5) -- (2,-1);
\draw [red,->-] (2,0.5) -- (2.5,-1);
\draw [red,->-] (1.5,-1) -- (2,0.5);
\node at (2,0.5) [] {$\bullet$};
\node at (2,0.5) [above] {$i$};
\end{tikzpicture} \quad \quad \longrightarrow \quad\quad
\begin{tikzpicture}[baseline=0ex]
\draw (4,0) arc (60:120:4);
\draw [red,->-] (2.5,0.47) -- (2.5,-1);
\draw [red,->-] (1.5,-1) -- (1.5,0.47);
\node at (2.5,0.47) [] {$\bullet$};
\node at (1.5,0.47) [] {$\bullet$};
\node at (2.5,0.47) [above] {$i$};
\node at (1.5,0.47) [above] {$i'$};
\end{tikzpicture}
\]
\[
\begin{tikzpicture}[baseline=0ex]
\foreach \i in {0,1,2}
    {
    \draw (0,0) -- +(90-\i*120:2);
    \coordinate (i\i) at +(85-\i*120:2);
    \coordinate (m\i) at +(150-\i*120:0.5);
    \coordinate (o\i) at +(215-\i*120:2);
    \draw [red, ->] plot [smooth, tension=1] coordinates {(i\i)(m\i)(o\i)};
    }
\draw [fill] (0,0) circle [radius=0.2];
\end{tikzpicture} \quad \quad \quad \longrightarrow \quad \quad \quad
\begin{tikzpicture}[baseline=0ex]
\foreach \i in {0,1,2}
    {
    \coordinate (i\i) at +(80-\i*120:2);
    \coordinate (o\i) at +(220-\i*120:2);
    \draw [red, ->] plot [smooth, tension=1] coordinates {(i\i)(0,0)(o\i)};
    }
\end{tikzpicture} 
\]

There is a reverse operation one can perform to get a bipartite graph out from any minimal triple diagram with assignment $i \mapsto (i+m-1)'$ (modulo $n$): put a black vertex at each triple crossing and put a white vertex at the center of any region with clockwisely oriented strands surrounding; then connect each white vertex to all the intersections of the surrounding strands (which can be either a black vertex or a boundary marked point). Note that the black vertices of a bipartite graph constructed from a minimal triple diagram in this way are necessarily trivalent.

An important result about minimal triple diagram is the following.

\begin{thm}[Thurston, Theorem 6 of \cite{Thu}] Any two minimal triple diagrams obtained from bipartite graphs with the same parameters $m$ and $n$ can be transformed into one another via a sequence of 2-by-2 moves, which come in the following two types.
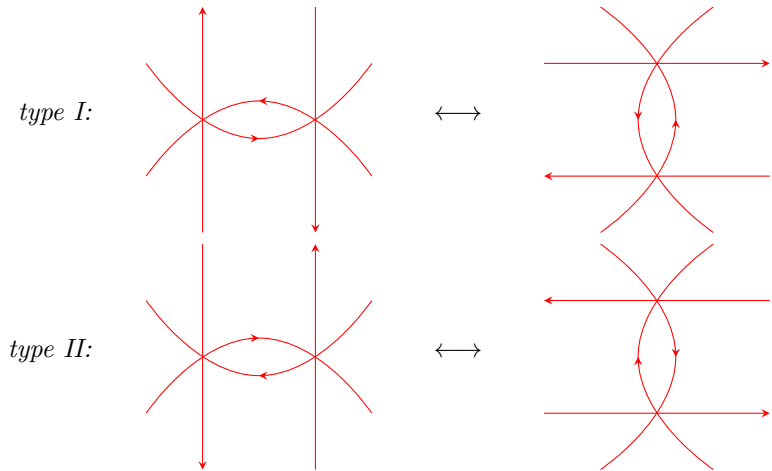
\begin{align*}
    \text{type I:} \quad \quad 
    \begin{tikzpicture}[baseline=0ex]
    \draw [->-,red] plot [smooth, tension=1] coordinates {(-1.5,0.75)(0,-0.25)(1.5,0.75)};
    \draw [->-,red] plot [smooth, tension=1] coordinates {(1.5,-0.75)(0,0.25)(-1.5,-0.75)};
    \draw [->, red] (-0.75,-1.5) -- (-0.75,1.5);
    \draw [->, red] (0.75,1.5) -- (0.75,-1.5);
    \end{tikzpicture} \quad \quad &\longleftrightarrow \quad \quad
    \begin{tikzpicture}[baseline=0ex]
    \draw [->-,red] plot [smooth, tension=1] coordinates {(0.75,1.5)(-0.25,0)(0.75,-1.5)};
    \draw [->-,red] plot [smooth, tension=1] coordinates {(-0.75,-1.5)(0.25,0)(-0.75,1.5)};
    \draw [->, red] (-1.5,0.75) -- (1.5,0.75);
    \draw [->, red] (1.5,-0.75) -- (-1.5,-0.75);
    \end{tikzpicture}\\
    \text{type II:} \quad \quad 
    \begin{tikzpicture}[baseline=0ex]
    \draw [->-,red] plot [smooth, tension=1] coordinates {(1.5,0.75)(0,-0.25)(-1.5,0.75)};
    \draw [->-,red] plot [smooth, tension=1] coordinates {(-1.5,-0.75)(0,0.25)(1.5,-0.75)};
    \draw [<-, red] (-0.75,-1.5) -- (-0.75,1.5);
    \draw [<-, red] (0.75,1.5) -- (0.75,-1.5);
    \end{tikzpicture} \quad \quad &\longleftrightarrow \quad \quad
    \begin{tikzpicture}[baseline=0ex]
    \draw [->-,red] plot [smooth, tension=1] coordinates {(0.75,-1.5)(-0.25,0)(0.75,1.5)};
    \draw [->-,red] plot [smooth, tension=1] coordinates {(-0.75,1.5)(0.25,0)(-0.75,-1.5)};
    \draw [<-, red] (-1.5,0.75) -- (1.5,0.75);
    \draw [<-, red] (1.5,-0.75) -- (-1.5,-0.75);
    \end{tikzpicture}
\end{align*}
\end{thm}

We can translate these two types of 2-by-2 moves back into bipartite graph language, which are the following two moves.
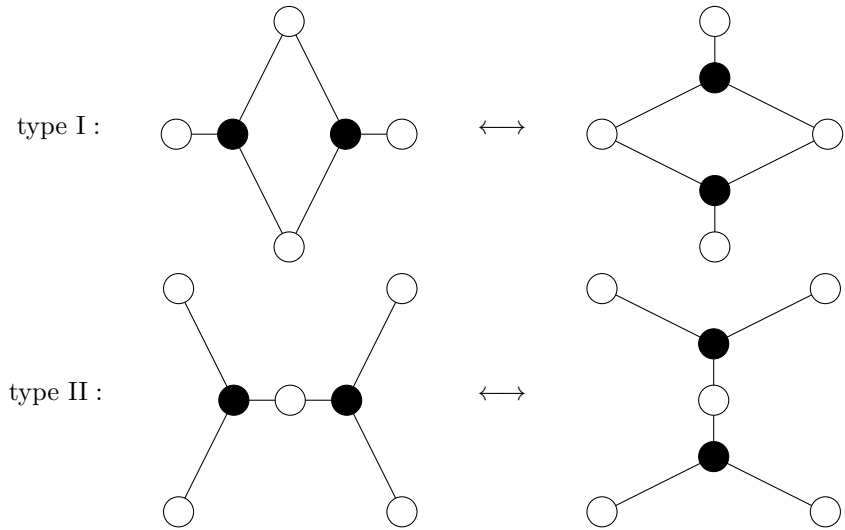
\begin{align*}
    \text{type I}: \quad \quad
    \begin{tikzpicture}[baseline=0ex]
    \draw (1.5,0) -- (0.75,0);
    \draw (-1.5,0) -- (-0.75,0);
    \draw (0.75,0) -- (0,1.5) -- (-0.75,0) -- (0,-1.5) --cycle;
    \draw [fill] (0.75,0) circle [radius=0.2];
    \draw [fill] (-0.75,0) circle [radius=0.2];
    \foreach \i in {1,...,4}
        {
        \draw [fill=white] +(90-\i*90:1.5) circle [radius=0.2];
        }
    \end{tikzpicture} \quad \quad & \longleftrightarrow \quad \quad
    \begin{tikzpicture}[baseline=0ex]
    \draw (0,1.5) -- (0,0.75);
    \draw (0,-1.5) -- (0,-0.75);
    \draw (0,0.75) -- (1.5,0) -- (0,-0.75) -- (-1.5,0) --cycle;
    \draw [fill] (0,0.75) circle [radius=0.2];
    \draw [fill] (0,-0.75) circle [radius=0.2];
    \foreach \i in {1,...,4}
        {
        \draw [fill=white] +(90-\i*90:1.5) circle [radius=0.2];
        }
    \end{tikzpicture}\\
    \text{type II}: \quad \quad 
    \begin{tikzpicture}[baseline=0ex]
    \draw (1.5,-1.5) -- (0.75,0) -- (1.5,1.5);
    \draw (-1.5,-1.5) -- (-0.75,0) -- (-1.5,1.5);
    \draw (0.75,0) -- (-0.75,0);
    \draw [fill] (0.75,0) circle [radius=0.2];
    \draw [fill] (-0.75,0) circle [radius=0.2];
    \foreach \i in {1,...,4}
        {
        \draw [fill=white] +(45-\i*90:2.1) circle [radius=0.2];
        }
    \draw [fill=white] (0,0) circle [radius=0.2];
    \end{tikzpicture}\quad \quad & \longleftrightarrow \quad \quad
    \begin{tikzpicture}[baseline=0ex]
    \draw (1.5,-1.5) -- (0,-0.75) -- (-1.5,-1.5);
    \draw (-1.5,1.5) -- (0,0.75) -- (1.5,1.5);
    \draw (0,0.75) -- (0,-0.75);
    \draw [fill] (0,0.75) circle [radius=0.2];
    \draw [fill] (0,-0.75) circle [radius=0.2];
    \foreach \i in {1,...,4}
        {
        \draw [fill=white] +(45-\i*90:2.1) circle [radius=0.2];
        }
    \draw [fill=white] (0,0) circle [radius=0.2];
    \end{tikzpicture}
\end{align*}
Using these two moves we can rewrite Thurston's theorem as follows.

\begin{thm}[Thurston]\label{thurston} Any two bipartite graphs with the same parameters $m$ and $n$ can be transformed into one another via a sequence of 2-by-2 moves described as above.
\end{thm}

\begin{conv} Note that a type II 2-by-2 move on a bipartite graph is equivalent to an merge-splitting move, which acts trivially on the quiver $\vec{i}_\Gamma$ and the dominating sets. Therefore from now on, by a \emph{2-by-2 move} on a bipartite graph, we always mean a type I 2-by-2 move, unless it is otherwise specified.
\end{conv}

\begin{prop}\label{bipartite graph mutation} A 2-by-2 move changes the quiver $\vec{i}_\Gamma$ and dominating sets $I(f)$ locally as follows, where $J$ is some $(m-2)$-element subset of $\{1,\dots, n\}$ that does not contain the indices $i,j,k,l$.
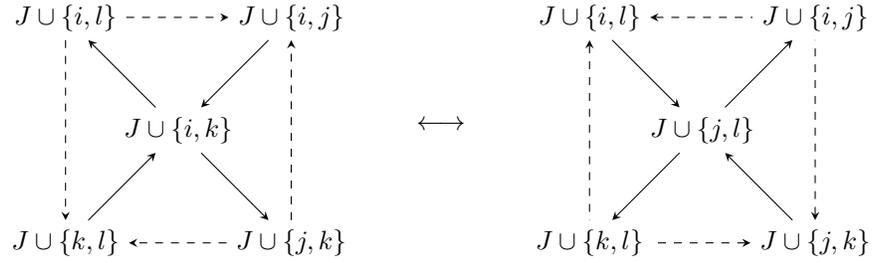
\[
\begin{tikzpicture}[baseline=0ex]
\node (ul) at (-1.5,1.5) [] {$J\cup\{i,l\}$};
\node (ur) at (1.5,1.5) [] {$J\cup \{i,j\}$};
\node (dr) at (1.5,-1.5) [] {$J\cup \{j,k\}$};
\node (dl) at (-1.5,-1.5) [] {$J\cup \{k,l\}$};
\node (c) at (0,0) [] {$J\cup\{i,k\}$};
\draw [->] (c) -- (ul);
\draw [dashed, ->] (ul) -- (dl);
\draw [dashed, ->] (dr) -- (dl);
\draw [->] (dl) -- (c);
\draw [->] (ur) -- (c);
\draw [->] (c) -- (dr);
\draw [dashed, ->] (dr) -- (ur);
\draw [dashed, ->] (ul) -- (ur);
\end{tikzpicture} \quad \quad  \longleftrightarrow \quad \quad
\begin{tikzpicture}[baseline=0ex]
\node (ul) at (-1.5,1.5) [] {$J\cup\{i,l\}$};
\node (ur) at (1.5,1.5) [] {$J\cup \{i,j\}$};
\node (dr) at (1.5,-1.5) [] {$J\cup \{j,k\}$};
\node (dl) at (-1.5,-1.5) [] {$J\cup \{k,l\}$};
\node (c) at (0,0) [] {$J\cup\{j,l\}$};
\draw [<-] (c) -- (ul);
\draw [dashed, <-] (ul) -- (dl);
\draw [dashed, <-] (dr) -- (dl);
\draw [<-] (dl) -- (c);
\draw [<-] (ur) -- (c);
\draw [<-] (c) -- (dr);
\draw [dashed, <-] (dr) -- (ur);
\draw [dashed, <-] (ul) -- (ur);
\end{tikzpicture} 
\]
\end{prop}
\begin{proof} Just consider the following labellings of zig-zag strands.
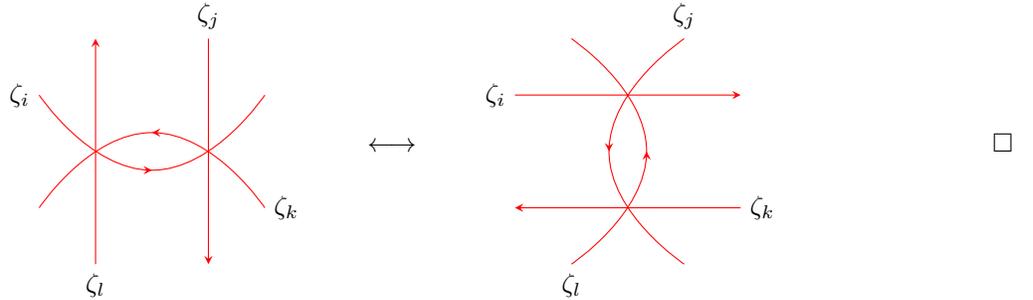
\[
\begin{tikzpicture}[baseline=0ex]
    \draw [->-,red] plot [smooth, tension=1] coordinates {(-1.5,0.75)(0,-0.25)(1.5,0.75)};
    \draw [->-,red] plot [smooth, tension=1] coordinates {(1.5,-0.75)(0,0.25)(-1.5,-0.75)};
    \draw [->, red] (-0.75,-1.5) -- (-0.75,1.5);
    \draw [->, red] (0.75,1.5) -- (0.75,-1.5);
    \node at (-1.5,0.75) [left] {$\zeta_i$};
    \node at (0.75,1.5) [above] {$\zeta_j$};
    \node at (1.5,-0.75) [right] {$\zeta_k$};
    \node at (-0.75,-1.5) [below] {$\zeta_l$};
    \end{tikzpicture} \quad \quad \longleftrightarrow \quad \quad
    \begin{tikzpicture}[baseline=0ex]
    \draw [->-,red] plot [smooth, tension=1] coordinates {(0.75,1.5)(-0.25,0)(0.75,-1.5)};
    \draw [->-,red] plot [smooth, tension=1] coordinates {(-0.75,-1.5)(0.25,0)(-0.75,1.5)};
    \draw [->, red] (-1.5,0.75) -- (1.5,0.75);
    \draw [->, red] (1.5,-0.75) -- (-1.5,-0.75);
    \node at (-1.5,0.75) [left] {$\zeta_i$};
    \node at (0.75,1.5) [above] {$\zeta_j$};
    \node at (1.5,-0.75) [right] {$\zeta_k$};
    \node at (-0.75,-1.5) [below] {$\zeta_l$};
    \end{tikzpicture} \qedhere
\]
\end{proof}

One key lemma about bipartite graph that Scott used to prove his result in \cite{Sco} is somehow missing in the literature. For completeness I present a proof of this lemma below using the following existence theorem of minimal triple diagram.

\begin{thm}[Thurston, Theorem 3 of \cite{Thu}] Given a collection of $2n$ boundary marked points on $\mathbb{D}$ prescribed as the initial points and terminal points of $n$ zig-zag strands, as long as they alternate between initial (in-coming) and terminal (out-going) as one goes around the boundary, there is a minimal triple diagram that satisfies such condition on the boundary data.
\end{thm}

\begin{lem}\label{any I} For any $m$-element subset $I$ of $\{1,\dots, n\}$, there is a bipartite graph $\Gamma$ in which there is a face $f$ such that $I(f)=I$.
\end{lem}
\begin{proof} Let $r$ be the number of elements in $I$ that is less than or equal to $m$. Consider the following rectangle with $2n$ marked points on the top edge and $2r$ marked points on the two vertical edges on the left and right. Our strategy is to find a suitable minimal triple diagram on this rectangle, and then by gluing the left and right edge together and gluing a disk to the bottom edge we hope to recover a minimal triple diagram whose corresponding bipartite graph has $I(f)=I$ for the face $f$ containing the disk that is glued to the bottom edge.
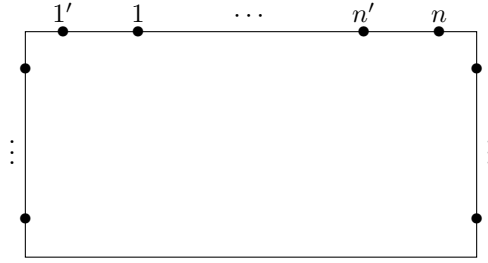
\[
\tikz{
\draw (0,0) rectangle (6,3);
\node at (0,0.5) [] {$\bullet$};
\node at (0,2.5) [] {$\bullet$};
\node at (0.5,3) [] {$\bullet$};
\node at (1.5,3) [] {$\bullet$};
\node at (4.5,3) [] {$\bullet$};
\node at (5.5,3) [] {$\bullet$};
\node at (6,0.5) [] {$\bullet$};
\node at (6,2.5) [] {$\bullet$};
\node at (0.5,3) [above] {$1'$};
\node at (1.5,3) [above] {$1$};
\node at (4.5,3) [above] {$n'$};
\node at (5.5,3) [above] {$n$};
\node at (3,3) [above] {$\cdots$};
\node at (0,1.5) [left] {$\vdots$};
\node at (6,1.5) [right] {$\vdots$};
}
\]

To apply the existence theorem, we need to specify which pair of boundary marked points with opposite orientation corresponds to the same zig-zag strand. For any $i$ in the subset
\[
S:=\left\{m< i\leq n-m\right\}\cup \left\{i\leq m, i\notin I\right\}\cup \left\{i>n-m,i\in I\right\},
\]
we pair up $i$ and $(i+m-1)'$ (modulo $n$) as usual. By simple counting we see that $\left|S^c\right|=2r$, in which exactly $r$ of them are in the set $I$. Denote these $r$ indices in $I$ by $i_1<i_2<\dots<i_r$ and denote the other $r$ indices not in $I$ by $j_1<j_2<\dots<j_r$. Now label the two vertical edges as follows.  
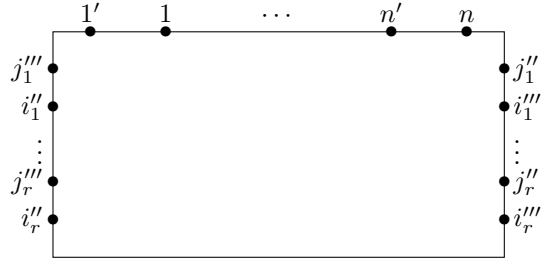
\[
\tikz{
\draw (0,0) rectangle (6,3);
\node at (0,0.5) [] {$\bullet$};
\node at (0,1) [] {$\bullet$};
\node at (0,2) [] {$\bullet$};
\node at (0,2.5) [] {$\bullet$};
\node at (0.5,3) [] {$\bullet$};
\node at (1.5,3) [] {$\bullet$};
\node at (4.5,3) [] {$\bullet$};
\node at (5.5,3) [] {$\bullet$};
\node at (6,0.5) [] {$\bullet$};
\node at (6,1) [] {$\bullet$};
\node at (6,2) [] {$\bullet$};
\node at (6,2.5) [] {$\bullet$};
\node at (0.5,3) [above] {$1'$};
\node at (1.5,3) [above] {$1$};
\node at (4.5,3) [above] {$n'$};
\node at (5.5,3) [above] {$n$};
\node at (6,2.5) [right] {$j''_1$};
\node at (6,2) [right] {$i'''_1$};
\node at (6,1) [right] {$j''_r$};
\node at (6,0.5) [right] {$i'''_r$};
\node at (0,2.5) [left] {$j'''_1$};
\node at (0,2) [left] {$i''_1$};
\node at (0,1) [left] {$j'''_r$};
\node at (0,0.5) [left] {$i''_r$};
\node at (3,3) [above] {$\cdots$};
\node at (0,1.5) [left] {$\vdots$};
\node at (6,1.5) [right] {$\vdots$};
}
\]
Now the rule to pair up boundary mark points for $k\in S^c$ is to pair up $k$ with $k''$ and to pair up $k'''$ with $(k+m-1)'$ (modulo $n$). Then by the existence theorem for minimal triple diagrams, we get a minimal triple diagram for such a pairing assignment, and it is not hard to see that after gluing, the face $f$ containing the glued-in disk has dominating set $I(f)=I$ by construction.

The only remaining thing we need to worry is whether the triple diagram remains minimal after the gluing. For this we need to check two things: no self-intersections and no parallel bigons. 

First let's consider self-intersections. The only kind of self-intersection that may arise from gluing is a tangling between the zig-zag strand going from $k$ to $k''$ with the zig-zag strand going from $k'''$ to $(k+m-1)'$ (modulo $n$). Due to the symmetry of the argument, we may assume without loss of generality that $k\in I$. Then part of the tangle looks like the following.
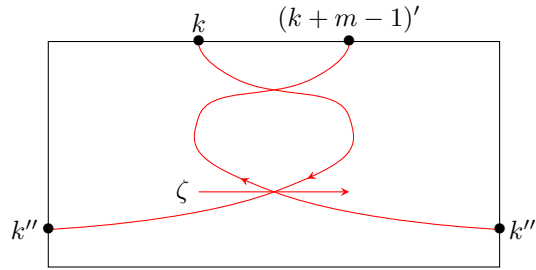
\[
\tikz{
\draw (0,0) rectangle (6,3);
\draw [->-, red] plot [smooth, tension=1] coordinates {(2,3)(2.5,2.5)(4,2)(3,1)(0,0.5)};
\draw [->-,red] plot [smooth,tension=1] coordinates {(6,0.5)(3,1)(2,2)(3.5,2.5)(4,3)};
\draw [->,red] (2,1) node [black,left] {$\zeta$} -- (4,1);
\node at (0,0.5) [] {$\bullet$};
\node at (6,0.5) [] {$\bullet$};
\node at (2,3) [] {$\bullet$};
\node at (4,3) [] {$\bullet$};
\node at (0,0.5) [left] {$k''$};
\node at (6,0.5) [right] {$k'''$};
\node at (2,3) [above] {$k$};
\node at (4,3) [above] {$(k+m-1)'$};
}
\]
Note that the zig-zag strand $\zeta$ must start at some $i<k$ and end at some $j'$ with $j>k+m-1$ to avoid forming parallel bigon with the two tangling zig-zag strands; but this is impossible since we require $i$ to go to $(i+m-1)'$.

Now let's turn to the parallel bigon problem. By using the minimality condition on the triple diagram on the rectangle, we see that for any $i_a, i_b\in S^c$, $\zeta_{i_a}$ and $\zeta_{i_b}$ do not intersect and $\zeta_{i'''_a}$ and $\zeta_{i'''_b}$ always intersect once, and for any $j_a,j_b\in S^c$, $\zeta_{j_a}$ and $\zeta_{j_b}$ always intersect once and $\zeta_{j'''_a}$ and $\zeta_{j'''_b}$ do not intersect. Therefore it is impossible for them to form a parallel bigon. The other possibilities of parallel bigons are already ruled by the minimality condition. This finishes the proof.
\end{proof}

\subsection{Examples of Bipartite Graphs on a Disk}

Before we proceed further, we would like to describe two special families of bipartite graphs for any given pair of parameters $(m,n)$, which will be useful for us in the development of this paper.

\begin{exmp}\label{honeycomb} The first one is constructed via the following process.
\begin{enumerate}
    \item Draw a horizontal pattern as below with $n-m$ black vertices and $n-m$ white vertices.
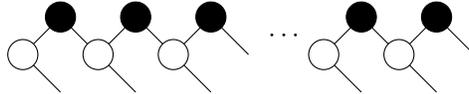
\[
\tikz{
\foreach \i in {0,1,2,4,5}
    {
    \draw (\i+1,0) -- (\i+0.5,0.5) -- (\i,0) -- (\i+0.5,-0.5);
    \draw [fill=white] (\i,0) circle [radius=0.2];
    \draw [fill] (\i+0.5,0.5) circle [radius=0.2];
    }
\node at (3.5,0.25) [] {$\cdots$};
}
\]
\item Stack $m-1$ layers of such horizontal patterns on top of each other.
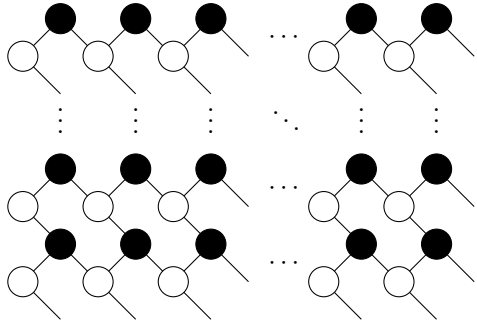
\[
\tikz{
\foreach \i in {0,1,2,4,5}
    {
    \draw (\i+1,3) -- (\i+0.5,3.5) -- (\i,3) -- (\i+0.5,2.5);
    \draw [fill=white] (\i,3) circle [radius=0.2];
    \draw [fill] (\i+0.5,3.5) circle [radius=0.2];
    }
\node at (3.5,3.25) [] {$\cdots$};
\foreach \i in {0,1,2,4,5}
    {
    \node at (\i+0.5,2.25) [] {$\vdots$};
    }
\node at (3.5,2.25) [] {$\ddots$};
\foreach \i in {0,1,2,4,5}
    {
    \draw (\i+1,1) -- (\i+0.5,1.5) -- (\i,1) -- (\i+0.5,0.5);
    \draw [fill=white] (\i,1) circle [radius=0.2];
    \draw [fill] (\i+0.5,1.5) circle [radius=0.2];
    }
\node at (3.5,1.25) [] {$\cdots$};
\foreach \i in {0,1,2,4,5}
    {
    \draw (\i+1,0) -- (\i+0.5,0.5) -- (\i,0) -- (\i+0.5,-0.5);
    \draw [fill=white] (\i,0) circle [radius=0.2];
    \draw [fill] (\i+0.5,0.5) circle [radius=0.2];
    }
\node at (3.5,0.25) [] {$\cdots$};
}
\]
\item Then we add a single white vertex with $n-m$ edges connecting to the $n-m$ black vertices on the top layer together with one extra free edge.
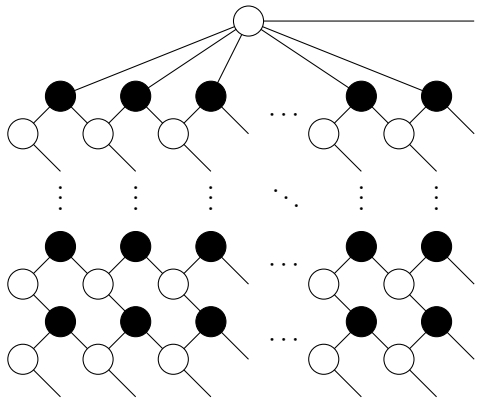
\[
\tikz{
\foreach \i in {0,1,2,4,5}
    {
    \draw (\i+0.5,3.5) -- (3,4.5);
    }
\foreach \i in {0,1,2,4,5}
    {
    \draw (\i+1,3) -- (\i+0.5,3.5) -- (\i,3) -- (\i+0.5,2.5);
    \draw [fill=white] (\i,3) circle [radius=0.2];
    \draw [fill] (\i+0.5,3.5) circle [radius=0.2];
    }
\node at (3.5,3.25) [] {$\cdots$};
\foreach \i in {0,1,2,4,5}
    {
    \node at (\i+0.5,2.25) [] {$\vdots$};
    }
    \node at (3.5,2.25) [] {$\ddots$};
\foreach \i in {0,1,2,4,5}
    {
    \draw (\i+1,1) -- (\i+0.5,1.5) -- (\i,1) -- (\i+0.5,0.5);
    \draw [fill=white] (\i,1) circle [radius=0.2];
    \draw [fill] (\i+0.5,1.5) circle [radius=0.2];
    }
\node at (3.5,1.25) [] {$\cdots$};
\foreach \i in {0,1,2,4,5}
    {
    \draw (\i+1,0) -- (\i+0.5,0.5) -- (\i,0) -- (\i+0.5,-0.5);
    \draw [fill=white] (\i,0) circle [radius=0.2];
    \draw [fill] (\i+0.5,0.5) circle [radius=0.2];
    }
    \draw (3,4.5) -- (6,4.5);
    \draw [fill=white] (3,4.5) circle [radius=0.2];
\node at (3.5,0.25) [] {$\cdots$};
}
\]
\item Draw a frame enclosing such bipartite graph; extend the $n$ free edges as external edges all the way to the frame, and label the intersections as boundary marked points in the following way.
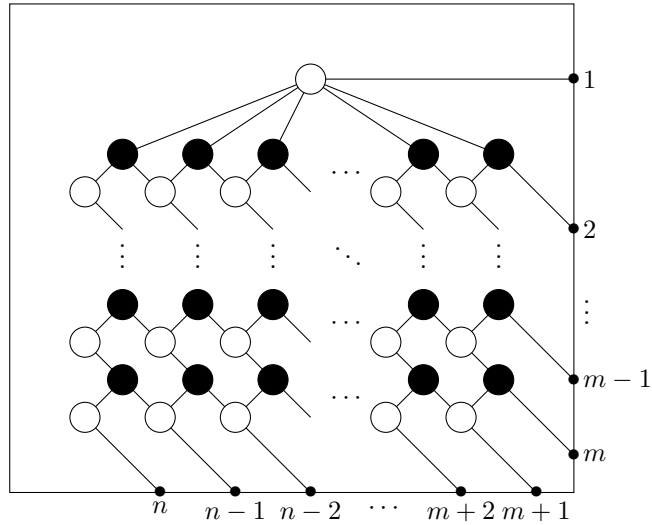
\[
\tikz{
\foreach \i in {0,1,2,4,5}
    {
    \draw (\i+0.5,3.5) -- (3,4.5);
    }
\foreach \i in {0,1,2,4,5}
    {
    \draw (\i+1,3) -- (\i+0.5,3.5) -- (\i,3) -- (\i+0.5,2.5);
    \draw [fill=white] (\i,3) circle [radius=0.2];
    \draw [fill] (\i+0.5,3.5) circle [radius=0.2];
    }
\node at (3.5,3.25) [] {$\cdots$};
\foreach \i in {0,1,2,4,5}
    {
    \node at (\i+0.5,2.25) [] {$\vdots$};
    }
    \node at (3.5,2.25) [] {$\ddots$};
\foreach \i in {0,1,2,4,5}
    {
    \draw (\i+1,1) -- (\i+0.5,1.5) -- (\i,1) -- (\i+0.5,0.5);
    \draw [fill=white] (\i,1) circle [radius=0.2];
    \draw [fill] (\i+0.5,1.5) circle [radius=0.2];
    }
\node at (3.5,1.25) [] {$\cdots$};
\foreach \i in {0,1,2,4,5}
    {
    \draw (\i+1,0) -- (\i+0.5,0.5) -- (\i,0) -- (\i+0.5,-0.5);
    \draw [fill=white] (\i,0) circle [radius=0.2];
    \draw [fill] (\i+0.5,0.5) circle [radius=0.2];
    }
\draw (3,4.5) -- (6.5,4.5);
\draw [fill=white] (3,4.5) circle [radius=0.2];
\draw (6,3) -- (6.5,2.5);
\draw (6,1) -- (6.5,0.5);
\draw (6,0) -- (6.5,-0.5);
\foreach \i in {1,2,3,5,6}
    {
    \draw (\i-0.5,-0.5) -- (\i,-1);
    }
\node at (3.5,0.25) [] {$\cdots$};
\draw (-1,-1) rectangle (6.5,5.5);
\node at (6.5,4.5) [] {$\bullet$};
\node at (6.5,4.5) [right] {$1$};
\node at (6.5,2.5) [] {$\bullet$};
\node at (6.5,2.5) [right] {$2$};
\node at (6.5,1.5) [right] {$\vdots$};
\node at (6.5,0.5) [] {$\bullet$};
\node at (6.5,0.5) [right] {$m-1$};
\node at (6.5,-0.5) [] {$\bullet$};
\node at (6.5,-0.5) [right] {$m$};
\node at (6,-1) [] {$\bullet$};
\node at (6,-1) [below] {$m+1$};
\node at (5,-1) [] {$\bullet$};
\node at (5,-1) [below] {$m+2$};
\node at (4,-1) [below] {$\cdots$};
\node at (3,-1) [] {$\bullet$};
\node at (3,-1) [below] {$n-2$};
\node at (2,-1) [] {$\bullet$};
\node at (2,-1) [below] {$n-1$};
\node at (1,-1) [] {$\bullet$};
\node at (1,-1) [below] {$n$};
}
\]
\end{enumerate}
Since there are lots of hexagons in such a bipartite graph, we call this bipartite graph the \emph{honeycomb bipartite graph} for the pair $(m,n)$ and denote it by $\Gamma_\hc$. Note that the honeycomb bipartite graph has the property that all its black vertices are trivalent. Below is an example of honeycomb bipartite graph for the pair $(m,n)=(3,7)$.
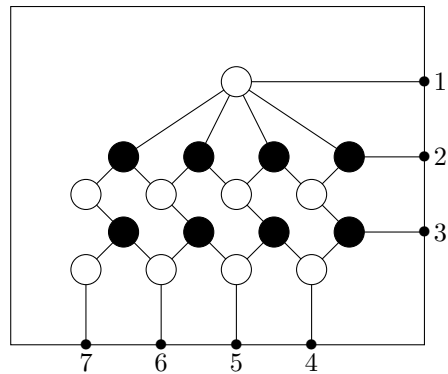
\[
\begin{tikzpicture}
\draw (0,0) rectangle (5.5,4.5);
\draw (3,3.5) -- (1.5,2.5);
\draw (3,3.5) -- (2.5,2.5);
\draw (3,3.5) -- (3.5,2.5);
\draw (3,3.5) -- (4.5,2.5);
\draw (1.5,2.5) -- (1,2);
\draw (2.5,2.5) -- (2,2);
\draw (3.5,2.5) -- (3,2);
\draw (4.5,2.5) -- (4,2);
\draw (1.5,2.5) -- (2,2);
\draw (2.5,2.5) -- (3,2);
\draw (3.5,2.5) -- (4,2);
\draw (1,2) -- (1.5,1.5);
\draw (2,2) -- (2.5,1.5);
\draw (3,2) -- (3.5,1.5);
\draw (4,2) -- (4.5,1.5);
\draw (1.5,1.5) -- (1,1);
\draw (2.5,1.5) -- (2,1);
\draw (3.5,1.5) -- (3,1);
\draw (4.5,1.5) -- (4,1);
\draw (1.5,1.5) -- (2,1);
\draw (2.5,1.5) -- (3,1);
\draw (3.5,1.5) -- (4,1);
\node at (5.5,3.5) [] {$\bullet$};
\node at (5.5,2.5) [] {$\bullet$};
\node at (5.5,1.5) [] {$\bullet$};
\node at (4,0) [] {$\bullet$};
\node at (3,0) [] {$\bullet$};
\node at (2,0) [] {$\bullet$};
\node at (1,0) [] {$\bullet$};
\draw (3,3.5) -- (5.5,3.5) node [right] {$1$};
\draw (4.5,2.5) -- (5.5,2.5) node [right] {$2$};
\draw (4.5,1.5) -- (5.5,1.5) node [right] {$3$};
\draw (4,1) -- (4,0) node [below] {$4$};
\draw (3,1) -- (3,0) node [below] {$5$};
\draw (2,1) -- (2,0) node [below] {$6$};
\draw (1,1) -- (1,0) node [below] {$7$};
\draw [fill=white] (1,1) circle [radius=0.2];
\draw [fill=white] (2,1) circle [radius=0.2];
\draw [fill=white] (3,1) circle [radius=0.2];
\draw [fill=white] (4,1) circle [radius=0.2];
\draw [fill] (1.5,1.5) circle [radius=0.2];
\draw [fill] (2.5,1.5) circle [radius=0.2];
\draw [fill] (3.5,1.5) circle [radius=0.2];
\draw [fill] (4.5,1.5) circle [radius=0.2];
\draw [fill=white] (1,2) circle [radius=0.2];
\draw [fill=white] (2,2) circle [radius=0.2];
\draw [fill=white] (3,2) circle [radius=0.2];
\draw [fill=white] (4,2) circle [radius=0.2];
\draw [fill] (1.5,2.5) circle [radius=0.2];
\draw [fill] (2.5,2.5) circle [radius=0.2];
\draw [fill] (3.5,2.5) circle [radius=0.2];
\draw [fill] (4.5,2.5) circle [radius=0.2];
\draw [fill=white] (3,3.5) circle [radius=0.2];
\end{tikzpicture} 
\]
\end{exmp}

\begin{exmp}\label{chessboard} The other one is constructed via the following process.

\begin{enumerate}
    \item Draw a grid of $(m-1)\times (n-m-1)$ squares; put a vertex at each crossing (including those at the edges and corners) in a bipartite fashion with a white vertex at the upper right hand corner.
    
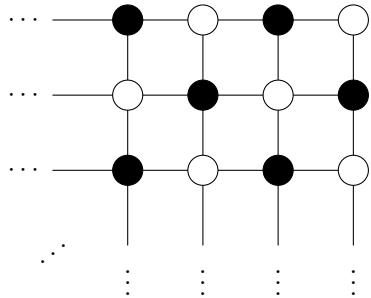
\[
    \tikz{
    \foreach \i in {1,2,3}
        {
        \draw (1,\i) -- (5,\i);
        \node at (1,\i) [left] {$\cdots$};
        }
    \foreach \j in {2,...,5}
        {
        \draw (\j,0) -- (\j,3);
        \node at (\j,0) [below] {$\vdots$};
        }
    \foreach \j in {3,5}
        {
        \draw [fill] (\j,2) circle [radius=0.2];
        \foreach \i in {1,3}
            {
            \draw [fill=white] (\j,\i) circle [radius=0.2];
            }
        }
    \foreach \j in {2,4}
        {
        \draw [fill=white] (\j,2) circle [radius=0.2];
        \foreach \i in {1,3}
            {
            \draw [fill] (\j,\i) circle [radius=0.2];
            }
        }
    \node at (1,0) [] {$\iddots$};
    }
    \]
\item Add $n$ free edges to the boundary of this grid following the following rules: first, each corner has a free edge; the spacing between every two consecutive free edges should be 2, with four and only four exceptions, one for each corner; the exceptional spacing is 1, and the exception at the upper right hand corner is along the vertical edge.
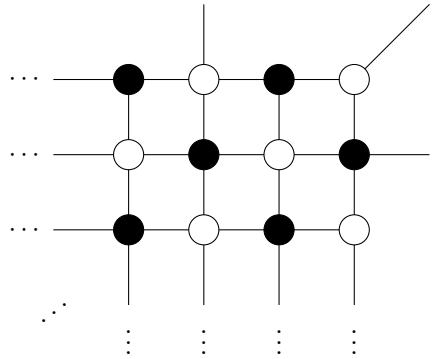
\[
    \tikz{
    \draw (5,3) -- (6,4);
    \draw (5,2) -- (6,2);
    \draw (3,3) -- (3,4);
    \foreach \i in {1,2,3}
        {
        \draw (1,\i) -- (5,\i);
        \node at (1,\i) [left] {$\cdots$};
        }
    \foreach \j in {2,...,5}
        {
        \draw (\j,0) -- (\j,3);
        \node at (\j,0) [below] {$\vdots$};
        }
    \foreach \j in {3,5}
        {
        \draw [fill] (\j,2) circle [radius=0.2];
        \foreach \i in {1,3}
            {
            \draw [fill=white] (\j,\i) circle [radius=0.2];
            }
        }
    \foreach \j in {2,4}
        {
        \draw [fill=white] (\j,2) circle [radius=0.2];
        \foreach \i in {1,3}
            {
            \draw [fill] (\j,\i) circle [radius=0.2];
            }
        }
    \node at (1,0) [] {$\iddots$};
    }
\]
\item Draw a frame enclosing such bipartite graph; extend the $n$ free edges as external edges all the way to the frame, and label the intersections as boundary marked points in the following way.
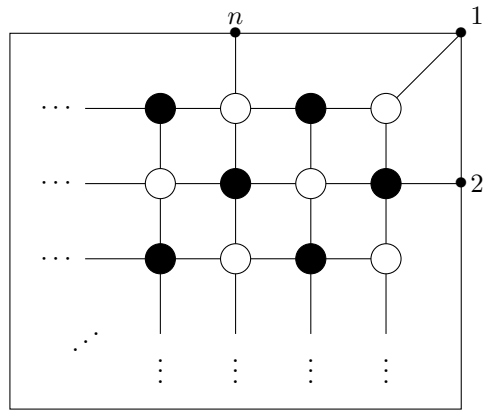
\[
    \tikz{
    \draw (0,-1) rectangle (6,4);
    \node at (6,4) [] {$\bullet$};
    \node at (6,2) [] {$\bullet$};
    \node at (3,4) [] {$\bullet$};
    \node at (6,4) [above right] {$1$};
    \node at (6,2) [right] {$2$};
    \node at (3,4) [above] {$n$};
    \draw (5,3) -- (6,4);
    \draw (5,2) -- (6,2);
    \draw (3,3) -- (3,4);
    \foreach \i in {1,2,3}
        {
        \draw (1,\i) -- (5,\i);
        \node at (1,\i) [left] {$\cdots$};
        }
    \foreach \j in {2,...,5}
        {
        \draw (\j,0) -- (\j,3);
        \node at (\j,0) [below] {$\vdots$};
        }
    \foreach \j in {3,5}
        {
        \draw [fill] (\j,2) circle [radius=0.2];
        \foreach \i in {1,3}
            {
            \draw [fill=white] (\j,\i) circle [radius=0.2];
            }
        }
    \foreach \j in {2,4}
        {
        \draw [fill=white] (\j,2) circle [radius=0.2];
        \foreach \i in {1,3}
            {
            \draw [fill] (\j,\i) circle [radius=0.2];
            }
        }
    \node at (1,0) [] {$\iddots$};
    }
\]
\end{enumerate}
Since there are lots of squares in such a bipartite graph, we call this bipartite graph the \emph{chessboard bipartite graph} for the pair $(m,n)$ and denote it by $\Gamma_\cb$. Note that chessboard bipartite graph has the property that one can perform a type I 2-by-2 move at any of its faces. Below are two examples of chessboard bipartite graphs.
\[
(m,n)=(2,7): \quad \quad \quad \quad 
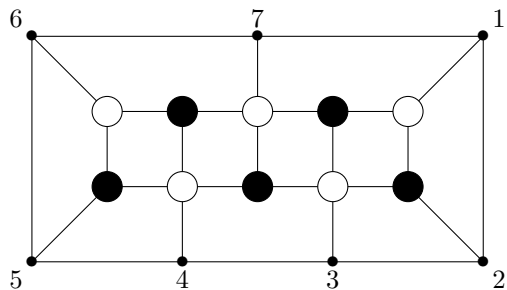
\begin{tikzpicture}[baseline=0ex]
\draw (1,1) -- (5,1);
\draw (1,0) -- (5,0);
\foreach \i in {1,...,5}
    {
    \draw (\i,0)--(\i,1);
    }
\draw (5,1) -- (6,2);
\draw (5,0) -- (6,-1);
\draw (4,0) -- (4,-1);
\draw (2,0) -- (2,-1);
\draw (1,0) -- (0,-1);
\draw (1,1) -- (0,2);
\draw (3,1) -- (3,2);
\draw (0,-1) rectangle (6,2);
\node at (6,2) [] {$\bullet$};
\node at (6,2) [above right] {$1$};
\node at (6,-1) [] {$\bullet$};
\node at (6,-1) [below right] {$2$};
\node at (4,-1) [] {$\bullet$};
\node at (4,-1) [below] {$3$};
\node at (2,-1) [] {$\bullet$};
\node at (2,-1) [below] {$4$};
\node at (0,-1) [] {$\bullet$};
\node at (0,-1) [below left] {$5$};
\node at (0,2) [] {$\bullet$};
\node at (0,2) [above left] {$6$};
\node at (3,2) [] {$\bullet$};
\node at (3,2) [above] {$7$};
\foreach \i in {0,2,4}
    {
    \draw [fill=white] (\i+1,1) circle [radius=0.2];
    \draw [fill] (\i+1,0) circle [radius=0.2];
    }
\foreach \i in {2,4}
    {
    \draw [fill=white] (\i,0) circle [radius=0.2];
    \draw [fill] (\i,1) circle [radius=0.2];
    }
\end{tikzpicture}
\]
\[
(m,n)=(4,8): \quad \quad \quad \quad 
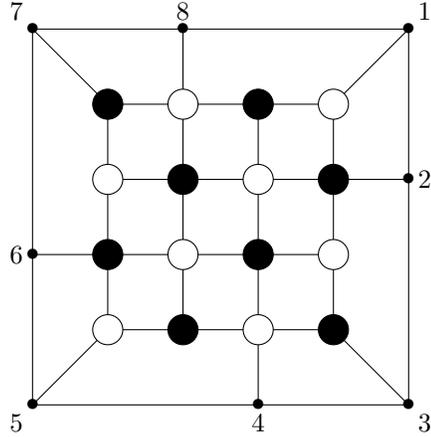
\begin{tikzpicture}[baseline=8ex]
    \draw (1,-1) rectangle (6,4);
    \node at (6,4) [] {$\bullet$};
    \node at (6,2) [] {$\bullet$};
    \node at (3,4) [] {$\bullet$};
    \node at (6,-1) [] {$\bullet$};
    \node at (4,-1) [] {$\bullet$};
    \node at (1,-1) [] {$\bullet$};
    \node at (1,1) [] {$\bullet$};
    \node at (1,4) [] {$\bullet$};
    \node at (6,4) [above right] {$1$};
    \node at (6,2) [right] {$2$};
    \node at (6,-1) [below right] {$3$};
    \node at (4,-1) [below] {$4$};
    \node at (1,-1) [below left] {$5$};
    \node at (1,1) [left] {$6$};
    \node at (1,4) [above left] {$7$};
    \node at (3,4) [above] {$8$};
    \draw (5,3) -- (6,4);
    \draw (5,2) -- (6,2);
    \draw (3,3) -- (3,4);
    \draw (5,0) -- (6,-1);
    \draw (4,0) -- (4,-1);
    \draw (2,0) -- (1,-1);
    \draw (2,1) -- (1,1);
    \draw (2,3) -- (1,4);
    \foreach \i in {0,...,3}
        {
        \draw (2,\i) -- (5,\i);
        }
    \foreach \j in {2,...,5}
        {
        \draw (\j,0) -- (\j,3);
        }
    \foreach \j in {3,5}
        {
        \draw [fill] (\j,2) circle [radius=0.2];
        \foreach \i in {1,3}
            {
            \draw [fill=white] (\j,\i) circle [radius=0.2];
            }
        }
    \foreach \j in {2,4}
        {
        \draw [fill=white] (\j,2) circle [radius=0.2];
        \foreach \i in {1,3}
            {
            \draw [fill] (\j,\i) circle [radius=0.2];
            }
        }
    \foreach \i in {2,4}
        {
        \draw [fill=white] (\i,0) circle [radius=0.2];
        \draw [fill] (\i+1,0) circle [radius=0.2];
        }
\end{tikzpicture}
\]

The chessboard bipartite graph is equivalent to \emph{Postnikov's arrangement} used by Scott in proving the cluster algebra structure on the homogeneous coordinate ring of Grassmannian $\Gr_{m,n}$ \cite{Sco}; in fact, Postnikov's arrangement can be seen as the zig-zag strands of the corresponding chessboard bipartite graph. 
\end{exmp}

\subsection{Perfect Orientation}

Perfect orientation is an assignment of orientations to edges of a bipartite graph. Perfect orientation was introduced by Postnikov to parametrize the strata of positive Grassmannian in \cite{Pos}. Let's first recall the general definition of a perfect orientation.

\begin{defn} A \emph{perfect orientation} on a bipartite graph $\Gamma$ is an assignment of orientation to each edge of $\Gamma$ such that there is only one in-coming edge for each white vertex and only one out-going edge for each black vertex.
\end{defn}

\begin{exmp} Below are two examples of perfect orientations on the same bipartite graph. Note that there can be multiple perfect orientations on the same graph, and some perfect orientations may even have cycles (e.g., the picture on the left).
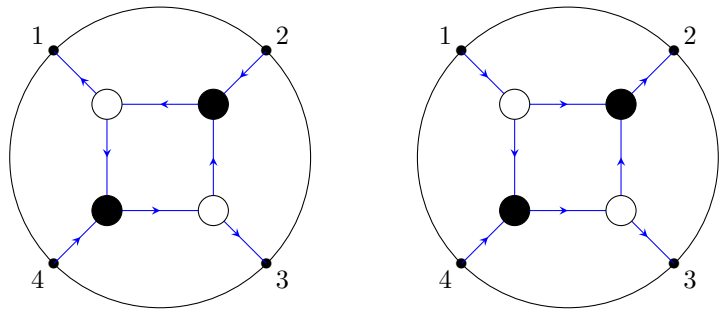
\[
\tikz{
\foreach \i in {1,...,4}
	{
	\coordinate (o\i) at +(225-90*\i:2);
	\coordinate (i\i) at +(225-90*\i:1);
	\node at (o\i) [] {$\bullet$};
	\node at +(225-90*\i:2.3) [] {$\i$};
	}
\draw [->-,blue] (i1) -- (o1);
\draw [->-,blue] (o2) -- (i2);
\draw [->-,blue] (i3) -- (o3);
\draw [->-,blue] (o4) -- (i4);
\draw [->-, blue](i2) -- (i1);
\draw [->-,blue] (i1) -- (i4);
\draw [->-,blue] (i3) -- (i2);
\draw [->-,blue] (i4) -- (i3);
\draw [fill=white] (i1) circle [radius=0.2];
\draw [fill=white] (i3) circle [radius=0.2];
\draw [fill] (i2) circle [radius=0.2];
\draw [fill] (i4) circle [radius=0.2];
\draw (0,0) circle [radius=2]; 
}\quad \quad \quad \quad 
\tikz{
\foreach \i in {1,...,4}
	{
	\coordinate (o\i) at +(225-90*\i:2);
	\coordinate (i\i) at +(225-90*\i:1);
	\node at (o\i) [] {$\bullet$};
	\node at +(225-90*\i:2.3) [] {$\i$};
	}
\draw [->-,blue] (o1) -- (i1);
\draw [->-,blue] (i2) -- (o2);
\draw [->-,blue] (i3) -- (o3);
\draw [->-,blue] (o4) -- (i4);
\draw [->-, blue](i1) -- (i2);
\draw [->-,blue] (i1) -- (i4);
\draw [->-,blue] (i3) -- (i2);
\draw [->-,blue] (i4) -- (i3);
\draw [fill=white] (i1) circle [radius=0.2];
\draw [fill=white] (i3) circle [radius=0.2];
\draw [fill] (i2) circle [radius=0.2];
\draw [fill] (i4) circle [radius=0.2];
\draw (0,0) circle [radius=2]; 
}
\]
\end{exmp}

In this paper we will only focus on a special choice of perfect orientation which we call the \emph{standard perfect orientation}. The construction of the standard perfect orientation on a bipartite graph $\Gamma$ goes as follows. First draw the zig-zag strands on $\Gamma$. Then a white vertex of valence $v$ there are $v$ zig-zag strands surrounding it. Due to the minimality condition, the indices of these zig-zag strands must satisfy the same cyclic ordering as the labellings of the boundary marked points. This implies that there must be a unique pair of consecutive zig-zag strands $\zeta_i$ and $\zeta_j$ with $\zeta_j$ succeeding $\zeta_i$ in the clockwise direction but $i>j$ with respect to the linear ordering on $\{1,\dots, n\}$ (as integers), and we demand to orient the edge between these two zig-zag strands to be pointing towards this white vertex $i$ while all other edges are pointing away.
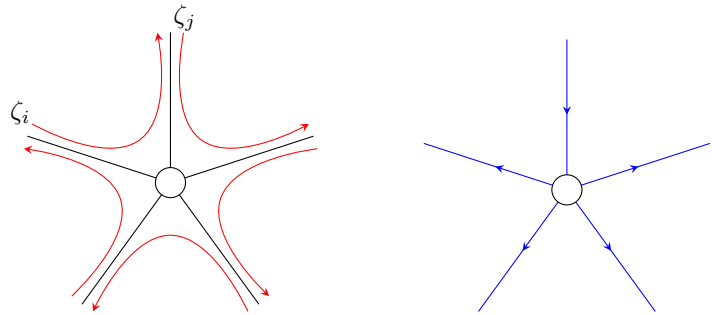
\[
\tikz{
\node at +(155:2.2) [] {$\zeta_i$};
\node at +(85:2.2)[] {$\zeta_j$};
\foreach \i in {0,...,4}
    {
    \draw (0,0) -- +(90-72*\i:2);
    \coordinate (i\i) at +(85-72*\i:2);
    \coordinate (m\i) at +(54-72*\i:0.7);
    \coordinate (o\i) at +(23-72*\i:2);
    \draw [red,->] plot [smooth, tension=1] coordinates {(i\i)(m\i)(o\i)};
    }
\draw [fill=white] (0,0) circle [radius=0.2];
}\quad \quad \quad \quad
\tikz{
\foreach \i in {1,...,4}
    {
    \draw [blue, ->-] (0,0) -- +(90-72*\i:2);
    }
\draw [blue, ->-] (0,2) -- (0,0);
\draw [fill=white] (0,0) circle [radius=0.2];
}
\]
If there is an external edge directly connecting to a black vertex, one can determine its orientation by adding a bivalent white vertex between the boundary marked point and the black vertex it is connecting to.

\begin{prop} Each black vertex has exactly one out-going edge.
\end{prop}
\begin{proof} Similar to white vertices, the indices of the zig-zag strands near a black vertex also satisfy the same cyclic ordering as the boundary marked points because of the minimality condition, and the out-going edge is again exactly the edge where the linear ordering breaks down.
\end{proof}

\begin{cor} The deleting move \eqref{delete} and the splitting move \eqref{splitting} are compatible with standard perfect orientation, i.e., there is a natural and unique way of deciding the orientations on the edges that are changed under either move. 
\end{cor}
\begin{proof} This follows from the fact that there is exactly one out-going edge for each black vertex.
\end{proof}

Let's investigate how standard perfect orientation changes under a 2-by-2 move. Without loss of generality let's assume that $i<j<k<l$ with respect to the linear ordering on $\{1,\dots, n\}$.
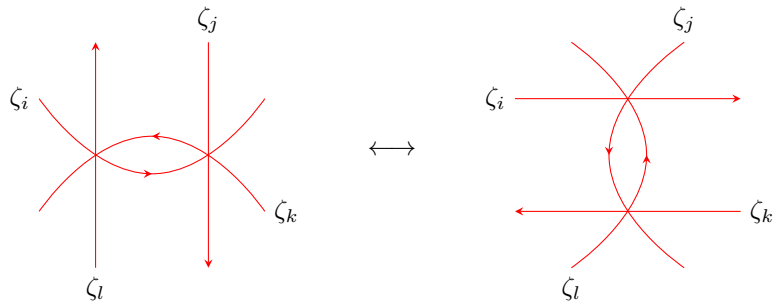
\[
\begin{tikzpicture}[baseline=0ex]
    \draw [->-,red] plot [smooth, tension=1] coordinates {(-1.5,0.75)(0,-0.25)(1.5,0.75)};
    \draw [->-,red] plot [smooth, tension=1] coordinates {(1.5,-0.75)(0,0.25)(-1.5,-0.75)};
    \draw [->, red] (-0.75,-1.5) -- (-0.75,1.5);
    \draw [->, red] (0.75,1.5) -- (0.75,-1.5);
    \node at (-1.5,0.75) [left] {$\zeta_i$};
    \node at (0.75,1.5) [above] {$\zeta_j$};
    \node at (1.5,-0.75) [right] {$\zeta_k$};
    \node at (-0.75,-1.5) [below] {$\zeta_l$};
    \end{tikzpicture} \quad \quad \longleftrightarrow \quad \quad
    \begin{tikzpicture}[baseline=0ex]
    \draw [->-,red] plot [smooth, tension=1] coordinates {(0.75,1.5)(-0.25,0)(0.75,-1.5)};
    \draw [->-,red] plot [smooth, tension=1] coordinates {(-0.75,-1.5)(0.25,0)(-0.75,1.5)};
    \draw [->, red] (-1.5,0.75) -- (1.5,0.75);
    \draw [->, red] (1.5,-0.75) -- (-1.5,-0.75);
    \node at (-1.5,0.75) [left] {$\zeta_i$};
    \node at (0.75,1.5) [above] {$\zeta_j$};
    \node at (1.5,-0.75) [right] {$\zeta_k$};
    \node at (-0.75,-1.5) [below] {$\zeta_l$};
    \end{tikzpicture} 
\]
It is not hard to see that the induced standard perfect orientations before and after the 2-by-2 move are the following.
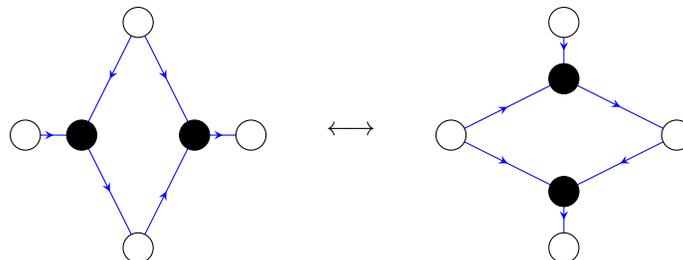
\[
\begin{tikzpicture}[baseline=0ex]
    \draw [blue, ->-] (-1.5,0) -- (-0.75,0);
    \draw [blue, ->-] (0.75,0) -- (1.5,0);
    \draw [blue,->-] (-0.75,0) -- (0,-1.5);
    \draw [blue,->-] (0,-1.5) -- (0.75,0);
    \draw [blue,->-] (0,1.5) --(-0.75,0);
    \draw [blue,->-] (0,1.5) -- (0.75,0);
    \draw [fill] (0.75,0) circle [radius=0.2];
    \draw [fill] (-0.75,0) circle [radius=0.2];
    \foreach \i in {1,...,4}
        {
        \draw [fill=white] +(90-\i*90:1.5) circle [radius=0.2];
        }
    \end{tikzpicture} \quad \quad  \longleftrightarrow \quad \quad
    \begin{tikzpicture}[baseline=0ex]
    \draw [blue, ->-] (-1.5,0) -- (0,0.75);
    \draw [blue, ->-] (0,1.5) -- (0,0.75);
    \draw [blue,->-] (0,0.75) -- (1.5,0);
    \draw [blue,->-] (-1.5,0) -- (0,-0.75);
    \draw [blue, ->-] (0,-0.75) -- (0,-1.5);
    \draw [blue, ->-] (1.5,0) -- (0,-0.75);
    \draw [fill] (0,0.75) circle [radius=0.2];
    \draw [fill] (0,-0.75) circle [radius=0.2];
    \foreach \i in {1,...,4}
        {
        \draw [fill=white] +(90-\i*90:1.5) circle [radius=0.2];
        }
    \end{tikzpicture}
\]
In particular, we observe that within such a local picture, a white vertex has an in-coming edge before the 2-by-2 move if and only if it has an in-coming edge after the 2-by-2 move. Therefore we see that 2-by-2 moves are really local: they only change the quiver $\vec{i}_\Gamma$, dominating sets $I(f)$, and the standard perfect orientation on $\Gamma$ locally.

Another important feature about the standard perfect orientation is its acyclicity, which directly implies the polynomial nature of the boundary measurement functions which we are going to define later.

\begin{prop} The standard perfect orientation on a bipartite graph $\Gamma$ is acyclic.
\end{prop}
\begin{proof} First we want to show that no boundary of a non-boundary face is cyclically oriented under the standard perfect orientation. Let's assume the contrary: let $f$ be a face with cyclically oriented edges surrounding it as depicted below.
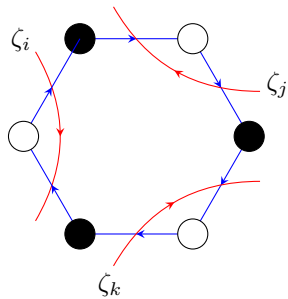
\[
\tikz{
\foreach \i in {0,...,2}
    {
    \draw [blue,->-] +(120-120*\i:1.5) -- +(60-120*\i:1.5);
    \draw [blue,->-] +(60-120*\i:1.5) -- +(-120*\i:1.5);
    \draw [fill] +(120-120*\i:1.5) circle [radius=0.2];
    \draw [fill=white] +(60-120*\i:1.5) circle [radius=0.2];
    \draw [red, ->-] +(140-120*\i:1.75) to [bend left] +(220-120*\i:1.75);
    }
\node at +(140:2) [] {$\zeta_i$};
\node at +(260:2) [] {$\zeta_k$};
\node at +(20:2) [] {$\zeta_j$};
}
\]
Since at each black vertex, all edges other than the two next to $f$ are in-coming, it follows that we  have $i<j<k<i$ with respect to the linear ordering, which is absurd.

Now suppose we have some oriented cycle that bounds a region with more than one face. Then some vertex of this oriented cycle must have edges going into the bounded region. Fix such an edge, which may be either inward- or outward-pointing under the standard perfect orientation; the proof of the two cases are symmetric, so without loss of generality we may assume that it is oriented inward. But then since each vertex, regardless of the color, has at least one out-going edge, we can continue traveling in this inward pointing direction and continue into some path inside the bounded region.
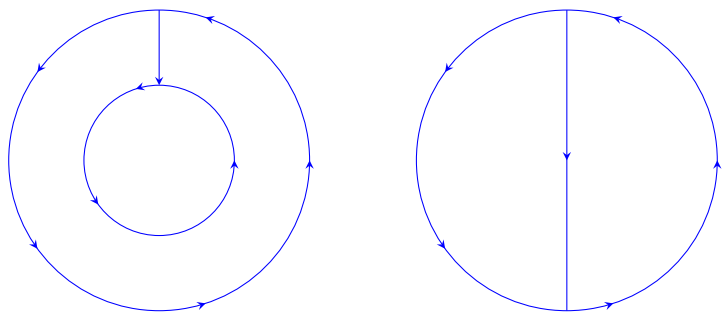
\[
\tikz{
\draw [decoration={markings, mark=at position 0.2 with {\arrow{>}},mark=at position 0.4 with {\arrow{>}},mark=at position 0.6 with {\arrow{>}},mark=at position 0.8 with {\arrow{>}},mark=at position 1 with {\arrow{>}}}, postaction={decorate},  blue] (0,0) circle [radius=2];
\draw [decoration={markings, mark=at position 0.3 with {\arrow{>}},mark=at position 0.6 with {\arrow{>}},mark=at position 1 with {\arrow{>}}}, postaction={decorate},  blue] (0,0) circle [radius=1];
\draw [ ->, blue] +(90:2) -- +(90:1);
} \quad \quad \quad \quad
\tikz{
\draw [decoration={markings, mark=at position 0.2 with {\arrow{>}},mark=at position 0.4 with {\arrow{>}},mark=at position 0.6 with {\arrow{>}},mark=at position 0.8 with {\arrow{>}},mark=at position 1 with {\arrow{>}}}, postaction={decorate},  blue] (0,0) circle [radius=2];
\draw [decoration={markings, mark=at position 0.5 with {\arrow{>}}}, postaction={decorate},  blue] +(90:2) -- +(-90:2);
}
\]
But since there are only finitely many edges inside this bounded region, we must end with one of the two cases above: either the path closes up and gives us a smaller oriented cycle with fewer faces inside (left), or it exits the bounded region by joining some other point on the oriented cycle (right). Either case gives us an oriented cycle with fewer faces inside, and the proof is then complete by induction on the number of faces the oriented cycle bounds. 
\end{proof}

\begin{defn} A boundary marked point is called a \emph{source} with respect to the perfect orientation if the external edge connecting to it points into the bipartite graph; otherwise it is called a \emph{sink}. 
\end{defn}

\begin{prop} The source set of the standard perfect orientation on a bipartite graph is always $\{1,\dots, m\}$.
\end{prop}
\begin{proof} By Thurston's theorem we know that any two bipartite graphs are related by a sequence of 2-by-2 moves (both types). But then since 2-by-2 moves (both types) are local operations, they do not change the source set. Therefore it suffices to find a bipartite graph with source set $\{1,\dots, m\}$. Here we choose the honeycomb bipartite graph, whose standard perfect orientation looks like the following.
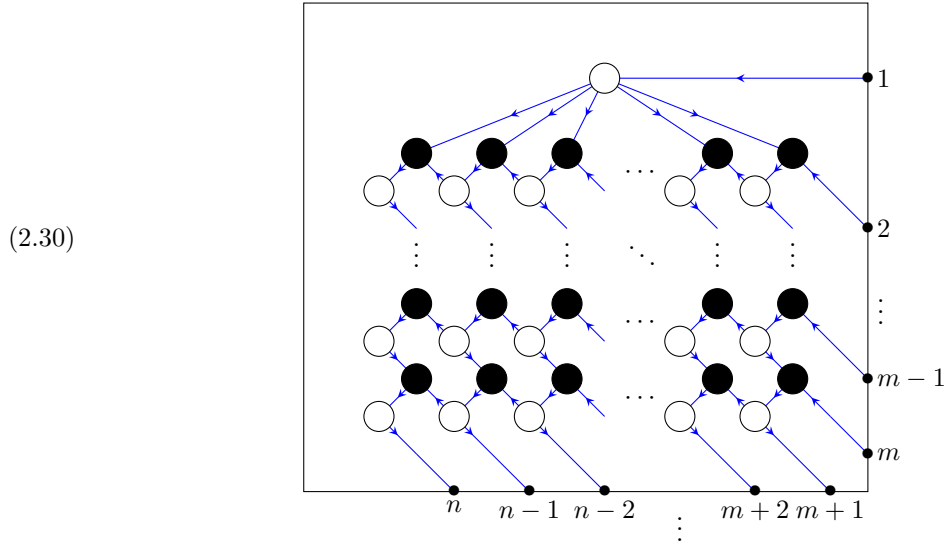
\begin{equation}\label{hc perfect orientation}
\begin{tikzpicture}[baseline=15ex]
\foreach \i in {0,1,2,4,5}
    {
    \draw [blue, ->-] (3,4.5)--(\i+0.5,3.5);
    }
\foreach \i in {0,1,2,4,5}
    {
    \draw [blue, ->-] (\i+1,3) -- (\i+0.5,3.5);
    \draw [blue, ->-] (\i+0.5,3.5) -- (\i,3);
    \draw [blue, ->-] (\i,3) -- (\i+0.5,2.5);
    \draw [fill=white] (\i,3) circle [radius=0.2];
    \draw [fill] (\i+0.5,3.5) circle [radius=0.2];
    }
\node at (3.5,3.25) [] {$\cdots$};
\foreach \i in {0,1,2,4,5}
    {
    \node at (\i+0.5,2.25) [] {$\vdots$};
    }
    \node at (3.5,2.25) [] {$\ddots$};
\foreach \i in {0,1,2,4,5}
    {
    \draw [blue, ->-] (\i+1,1) -- (\i+0.5,1.5);
    \draw [blue, ->-] (\i+0.5,1.5) -- (\i,1);
    \draw [blue, ->-] (\i,1) -- (\i+0.5,0.5);
    \draw [fill=white] (\i,1) circle [radius=0.2];
    \draw [fill] (\i+0.5,1.5) circle [radius=0.2];
    }
\node at (3.5,1.25) [] {$\cdots$};
\foreach \i in {0,1,2,4,5}
    {
    \draw [blue, ->-] (\i+1,0) -- (\i+0.5,0.5);
    \draw [blue, ->-] (\i+0.5,0.5) -- (\i,0);
    \draw [blue, ->-] (\i,0) -- (\i+0.5,-0.5);
    \draw [fill=white] (\i,0) circle [radius=0.2];
    \draw [fill] (\i+0.5,0.5) circle [radius=0.2];
    }
\draw [blue, ->-] (6.5,4.5)-- (3,4.5);
\draw [fill=white] (3,4.5) circle [radius=0.2];
\node at (3.5,0.25) [] {$\cdots$};
\foreach \i in {0,1,3}
    {
    \draw [blue] (6.5,\i-0.5) -- (6,\i);
    }
\foreach \i in {1,2,3,5,6}
    {
    \draw [blue] (\i, -1) -- (\i-0.5,-0.5);
    }
\draw (-1,-1) rectangle (6.5,5.5);
\node at (6.5,4.5) [] {$\bullet$};
\node at (6.5,4.5) [right] {$1$};
\node at (6.5,2.5) [] {$\bullet$};
\node at (6.5,2.5) [right] {$2$};
\node at (6.5,1.5) [right] {$\vdots$};
\node at (6.5,0.5) [] {$\bullet$};
\node at (6.5,0.5) [right] {$m-1$};
\node at (6.5,-0.5) [] {$\bullet$};
\node at (6.5,-0.5) [right] {$m$};
\node at (6,-1) [] {$\bullet$};
\node at (6,-1) [below] {$m+1$};
\node at (5,-1) [] {$\bullet$};
\node at (5,-1) [below] {$m+2$};
\node at (4,-1) [below] {$\vdots$};
\node at (3,-1) [] {$\bullet$};
\node at (3,-1) [below] {$n-2$};
\node at (2,-1) [] {$\bullet$};
\node at (2,-1) [below] {$n-1$};
\node at (1,-1) [] {$\bullet$};
\node at (1,-1) [below] {$n$};
\end{tikzpicture}
\end{equation}
Note that the source set is exactly $\{1,\dots, m\}$.
\end{proof}

\begin{rmk} Through out this paper we choose to work with the linear ordering on $\{1,\dots, n\}$ induced from the ordering on integers, but feel free to replace it by a cyclically shifted one, say 
\[
i\leq i+1\leq \dots \leq n\leq 1\leq 2\leq \dots \leq i-1.
\]
One can construct standard perfect orientation analogously with respect to such linear ordering, and the source set will always be $\{i,i+1,\dots, i+m-1\}$.
\end{rmk}

\section{Review of Cluster Theory}\label{section3}

Fock and Goncharov introduced cluster ensembles and their tropicalizations in \cite{FGensemble}, and formulated a duality conjecture describing a canonical basis of cluster algebras. In \cite{GHKK}, Gross, Hacking, Keel, and Kontsevich proved a weak version of the duality conjecture and gave a sufficient condition under which the full duality conjecture holds. In this section we will briefly review concepts and statements that are related to the development of this paper.

\subsection{Cluster Varieties and Cluster Transformations}

Let's start with the definition of a seed. 

\begin{defn} A \emph{seed} $\vec{i}$ is a triple $\left(I,I^\uf,\epsilon\right)$ satisfying the following properties:
\begin{enumerate}
    \item $I$ is a finite set;
    \item $I^\uf\subset I$;
    \item $\epsilon$ is an $\left(\frac{1}{2}\mathbb{Z}\right)$-valued $I\times I$ matrix such that $\epsilon_{ij}\in \mathbb{Z}$ unless both $i$ and $j$ are not in $I^\uf$.
\end{enumerate}
We call $I$ the set of \emph{vertices}, $I^\uf$ the set of \emph{unfrozen vertices}, and $\epsilon$ the \emph{exchange matrix}; it follows naturally that the set $I\setminus I^\uf$ is called the set of \emph{frozen vertices}. 
\end{defn}

It is obvious that the data of a seed is equivalent to the data of a quiver with possibly half edges and vertex set $I$ and exchange matrix $\epsilon$ but no loops or 2-cycles. Because of this fact, we will use the same symbol to denote seeds and quivers, and sometimes even use these two words interchangeably.

\begin{defn} Let $\vec{i}=\left(I,I^\uf,\epsilon\right)$ be a seed and let $k\in I^\uf$ be an unfrozen vertex. A \emph{mutation} in the direction of $k$ produces a new seed $\vec{i}'=\left(I',I'^\uf, \epsilon'\right)$ with $I':=I$, $I'^\uf:=I^\uf$, and
\[
\epsilon'_{ij}:=\left\{\begin{array}{ll}
    -\epsilon_{ij} & \text{if $k\in \{i,j\}$},  \\
    \epsilon_{ij}+\left[\epsilon_{ik}\right]_+\left[\epsilon_{kj}\right]_+-\left[-\epsilon_{ik}\right]_+\left[-\epsilon_{kj}\right]_+ & \text{if $k\notin\{i,j\}$}, 
\end{array}\right.
\]
where $\left[a\right]_+:=\max\{0,a\}$. We sometimes denote the mutated seed $\vec{i}'$ as $\mu_k\vec{i}$. One can verify easily that $\mu^2_k\vec{i}=\vec{i}$.
\end{defn}

If one translate the above formula into quiver language, one can see that seed mutation defined above coincide with quiver mutation defined by Derksen, Weyman, and Zelevinsky in \cite{DWZ}.

Given an initial seed $\vec{i}=\left(I,I^\uf,\epsilon\right)$ we define the \emph{mutation tree} $\mathfrak{T}$ to be an $\left|I^\uf\right|$-regular tree with edges labeled by elements of $I^\uf$ such that edges connected to the same vertex all have distinct labellings. Fix a vertex $v_0$ of $\mathfrak{T}$ (it doesn't matter which one because all of them are the same) and associate the initial seed $\vec{i}$ to this vertex. Note that for any other vertex $v$ there is a unique shortest path going from $v_0$ to $v$; by following such path and mutating $v_0$ according to the labellings of the edges on this path we get a seed new seed $\vec{i}_v$ which we associate to the vertex $v$. At the end each vertex of $\mathfrak{T}$ has a seed associated to it, and the labellings of edges along any paths between two vertices (not necessarily the shortest one) gives a sequence of mutations that turn one associated seed into the other.

Now we move to the geometry level to define cluster varieties. 

For each vertex $v$ on the mutation tree $\mathfrak{T}$ we construct two \emph{seed tori}: $\mathcal{A}_v:=\mathbb{G}_m^I$ and $\mathcal{X}_v:=\mathbb{G}_m^I$. Elements of $I$ specifies coordinate systems on $\mathcal{A}_v$ and $\mathcal{X}_v$, which we denote as $\left(A_{i;v}\right)_{i\in I}$ and $\left(X_{i;v}\right)_{i\in I}$ respectively. Between $\mathcal{A}_v$ and $\mathcal{X}_v$ there is a map
\begin{align*}
p_v:\mathcal{A}_v&\rightarrow \mathcal{X}_v\\
p_v^*\left(X_{i;v}\right)&=\prod_jA_{j;v}^{\epsilon_{ij;v}}
\end{align*}
where $\epsilon_{;v}$ denotes the exchange matrix of the seed $\vec{i}_v$. Note that the map $p_v$ may not be algebraic in general because $\epsilon_{ij}$ may not be integers when $i,j\notin I^\uf$.

For an edge $\begin{tikzpicture}[baseline=-0.5ex]\draw (0,0) node [left] {$v$} -- node [above] {$k$} (1,0) node [right] {$w$};\end{tikzpicture}$ in the mutation tree $\mathfrak{T}_{\vec{i}_0}$ we also define two birational maps
\[
\mu_k:\mathcal{A}_v\dashrightarrow \mathcal{A}_w \quad \quad \text{and} \quad \quad \mu_k:\mathcal{X}_v\dashrightarrow \mathcal{X}_w,
\]
which are called \emph{cluster mutations}, by setting
\begin{align}
\mu_k^*\left(A_{i;w}\right)=&\left\{\begin{array}{ll}
A_{k;v}^{-1}\prod_jA_{j;v}^{\left[-\epsilon_{kj;v}\right]_+}\left(1+\prod_jA_{j;v}^{\epsilon_{kj;v}}\right) & \text{if $i=k$},\\
    A_{i;v} & \text{if $i\neq k$},
\end{array}\right. \label{a mutation} \\
\mu_k^*\left(X_{i;w}\right)=&\left\{\begin{array}{ll}
    X_{k;v}^{-1} & \text{if $i=k$}, \\
    X_{i;v}X_{k;v}^{\left[\epsilon_{ik;v}\right]_+}\left(1+X_{k;v}\right)^{-\epsilon_{ik;v}} & \text{if $i\neq k$}. 
\end{array}\right. \label{x mutation}
\end{align}
One can verify that $\mu_k^2$ is birationally equivalent to the identity map in both cases. 

For any path $\gamma$ in the mutation tree $\mathfrak{T}$ going from a vertex $v$ to another vertex $w$ (not necessarily the shortest) we can compose the corresponding cluster mutations to get birational maps $\mu_\gamma$ between the respective seed tori. Since $\mu_k^2$ is birationally equivalent to the identity map, the birational map $\mu_\gamma$ actually only depends on the initial vertex and the terminal vertex of $\gamma$; therefore sometimes we also denote $\mu_\gamma$ as $\mu_{v\rightsquigarrow w}$.

For any two vertices $v$ and $w$ in the mutation tree $\mathfrak{T}$, Fomin and Zelevinsky \cite{FZIV} defined two $I\times I$ matrices $c_{;v\rightsquigarrow w}$ and $g_{;v\rightsquigarrow w}$ with integer entries and a collection of polynomials $F_{i;v\rightsquigarrow w}$ with positive integer coefficients in terms of cluster variables $X_{j;v}$ such that
\begin{equation}\begin{split}
\mu_{v\rightsquigarrow w}^*\left(A_{i;w}\right)=&\left.F_{i;v\rightsquigarrow w}\right|_{X_{j;v}=\prod_k A_{k;v}^{\epsilon_{jk;v}}}\prod_jA_{j;v}^{g_{ij;v\rightsquigarrow w}}\\
\mu_{v\rightsquigarrow w}^*\left(X_{i;w}\right)=&\prod_j F_{j;v\rightsquigarrow w}^{\epsilon_{ij;w}}\prod_j X_{j;v}^{c_{ij;v\rightsquigarrow w}}.\label{cf}
\end{split}
\end{equation}
The two matrices $c_{;v\rightsquigarrow w}:=\left(c_{ij;v\rightsquigarrow w}\right)$ and $g_{;v\rightsquigarrow w}:=\left(g_{ij;v\rightsquigarrow w}\right)$ are called the $c$-matrix and the $g$-matrix (associated to $w$ with respect to the initial seed $v$) respectively, and the polynomials $F_{i;v\rightsquigarrow w}$ are called the $F$-polynomials (associated to $w$ with respect to the initial seed $v$). The following statement describes some of their most important properties.

\begin{thm}\label{cgf} The following results about $c$-matrices, $g$-matrices, and $F$-polynomials are due to Gross, Hacking, Keel, and Kontsevich \cite{GHKK} and Nakanishi and Zelevinsky \cite{NZ}.
\begin{enumerate}
    \item Row vectors of $c$-matrices are sign-coherent \cite{GHKK}.
    \item Column vectors of $g$-matrices are sign-coherent \cite{GHKK}.
    \item Row vectors of a $g$-matrix generate the corresponding cluster chamber \cite{GHKK}.
    \item All $F$-polynomials have positive integer coefficients with constant terms equal to 1 \cite{GHKK}.
    \item Topical duality $g_{;v\rightsquigarrow w}^t=c_{;v\rightsquigarrow w}^{-1}$ \cite{NZ}.
\end{enumerate}
\end{thm}

\begin{defn} \emph{Cluster varieties} are varieties obtained by gluing the seed tori in each family via the cluster mutations, i.e.,
\[
\mathcal{A}:=\left.\bigsqcup_{v\in \mathfrak{T}}\mathcal{A}_v\right/\left\{\mu_\gamma\right\} \quad \quad \text{and}\quad \quad \mathcal{X}:=\left.\bigsqcup_{v\in \mathfrak{T}}\mathcal{X}_v\right/\left\{\mu_\gamma\right\}.
\]
It follows that each cluster torus gives a coordinate chart on the corresponding cluster variety, and cluster mutations become gluing maps. In particular, the coordinate functions $\left(A_{i;v}\right)_{i\in I}$ and $\left(X_{i;v}\right)_{i\in I}$ are called \emph{cluster coordinates} on $\mathcal{A}$ and $\mathcal{X}$ respectively.
\end{defn}

\begin{rmk}\label{exchange graph} In fact many of the gluings above are repetitive. One can reduce the amount of gluing by first defining an equivalence relation on $\mathfrak{T}$: for two vertices $v$ and $w$ in $\mathfrak{T}$ we say $v\sim w$ if there is a bijection $\sigma:I\rightarrow I$ with $\sigma|_{I\setminus I^\uf}=\mathrm{Id}_{I\setminus I^\uf}$ such that 
\[
\epsilon_{ij;w}=\epsilon_{\sigma(i)\sigma(j);v}, \quad \quad \mu_{v\rightsquigarrow w}^*\left(A_{i;w}\right)=A_{\sigma(i);v}, \quad \quad \text{and} \quad \quad \mu^*_{v\rightsquigarrow w} \left(X_{i;w}\right)=X_{\sigma(i);v}.
\]
In other words, up to a relabelling of the vertices, $\mu_{v\rightsquigarrow w}$ are birationally equivalent to the identity maps. The quotient $\mathfrak{E}:=\mathfrak{T}/\sim$ is called the \emph{exchange graph} associated to the seed $\vec{i}$. Fock and Goncharov gave four examples of such equivalences in \cite{FGensemble} corresponding to rank 2 Dynkin types $\mathrm{A}_1\times \mathrm{A}_1$, $\mathrm{A}_2$, $\mathrm{B}_2$, and $\mathrm{G}_2$.
\end{rmk}

For an edge $\begin{tikzpicture}[baseline=-0.5ex]\draw (0,0) node [left] {$v$} -- node [above] {$k$} (1,0) node [right] {$w$};\end{tikzpicture}$ in the mutation tree $\mathfrak{T}_{\vec{i}_0}$, one can verify that the following square commutes.
\[
\xymatrix{\mathcal{A}_v \ar@{-->}[r]^{\mu_k} \ar[d]_{p_v} & \mathcal{A}_w \ar[d]^{p_w} \\
\mathcal{X}_v \ar@{-->}[r]_{\mu_k} & \mathcal{X}_w
}
\]
Therefore after gluing the seed tori, the above commutative diagram gives rise to a map (may not be algebraic) between the two cluster varieties:
\[
p:\mathcal{A}\rightarrow \mathcal{X}.
\]
The pair $\left(\mathcal{A},\mathcal{X}\right)$ is called a \emph{cluster ensemble}.

On the cluster variety $\mathcal{X}$ there is a canonical Poisson structure given by 
\[
\left\{X_{i;v},X_{j;v}\right\}=\epsilon_{ij;v}X_{i;v}X_{j;v}
\]
on each seed torus. Thus the cluster variety $\mathcal{X}$ is also known as a \emph{cluster Poisson variety}. In particular, the image $p(\mathcal{A})$ is a symplectic leaf with respect to such Poisson structure (see \cite{FGensemble}). 

\begin{defn}\label{cluster transformation} For two vertices $v$ and $w$ in the mutation tree $\mathfrak{T}$, a \emph{seed isomorphism} is a bijection $\sigma:I\rightarrow I$ such that $\epsilon_{ij;w}=\epsilon_{\sigma(i)\sigma(j);v}$. Given a seed isomorphism $\sigma$ between $\vec{i}_v$ and $\vec{i}_w$ we define two automorphisms 
\[
\sigma:\mathcal{A}\rightarrow \mathcal{A} \quad \quad \text{and} \quad \quad\sigma:\mathcal{X}\rightarrow \mathcal{X},
\]
which are called \emph{cluster transformations}, by setting
\[
\sigma^*\left(A_{i;w}\right)=A_{\sigma(i);v} \quad \quad \text{and} \quad \quad \sigma^*\left(X_{i;w}\right)=X_{\sigma(i);v}.
\]
In particular, it follows that if we try to express the pull-backs of cluster coordinates in terms of cluster coordinates on the same seed torus, we get
\[
\sigma^*\left(A_{i;w}\right)=\mu_{w\rightsquigarrow v}^*\left(A_{\sigma(i);v}\right) \quad \quad \text{and} \quad \quad \sigma^*\left(X_{i;w}\right)=\mu_{w\rightsquigarrow v}^*\left(X_{\sigma(i);v}\right).
\]
\end{defn}

\begin{rmk} Please be aware that cluster mutations and cluster transformations are vastly different things: cluster mutations are change of coordinates between different coordinate charts whereas cluster transformations are actual automorphisms acting on the cluster variety.
\end{rmk}

\begin{rmk} Please also note the similarities and differences between the definition of seed isomorphisms and cluster transformation above and the equivalence defined in Remark \ref{exchange graph}: the equivalence relation defined in Remark \ref{exchange graph} precisely corresponds to the trivial cluster transformations (the ones act by identity on cluster varieties), and therefore it is safe to remove those extra gluings. The group of cluster transformations modulo the trivial ones is called the \emph{cluster modular group} and it acts on both cluster varieties. 
\end{rmk}

\begin{defn} On a cluster variety (either $\mathcal{A}$ or $\mathcal{X}$), a \emph{global monomial} is a regular function on the cluster variety which can be written as a monomial in terms of cluster coordinates on some seed torus. A global monomial on $\mathcal{A}$ is also called a \emph{cluster monomial}. It is known a monomial of the form $\prod_i A_{i;v}^{m_i}$ on some seed torus $\mathcal{A}_v$ is a cluster monomial if and only if $m_i$ are all non-negative integers.
\end{defn}

The algebra of regular functions $\mathcal{O}(\mathcal{A})$ is called the \emph{upper cluster algebra} associated to the mutation equivalent family of seeds $\left|\vec{i}\right|$. The subalgebra inside $\mathcal{O}(\mathcal{A})$ generated by the cluster monomials is called the \emph{cluster algebra} associated to the mutation equivalent family of seeds $\left|\vec{i}\right|$.

By forgetting the frozen part of the initial seed $\vec{i}=\left(I,I^\uf,\epsilon\right)$ we get a seed without frozen vertices $\vec{i}^\uf=\left(I^\uf,\emptyset, \epsilon|_{I^\uf\times I^\uf}\right)$. By going through the same construction we get another cluster ensemble $\left(\mathcal{A}^\uf,\mathcal{X}^\uf\right)$. We can relate the two cluster ensembles (defined by $\vec{i}$ and $\vec{i}^\uf$ respectively) via the following maps
\begin{equation}\label{uf}
\vcenter{\vbox{\xymatrix{ \mathcal{A}^\uf \ar[r]^e \ar[drr]_(0.3){p^\uf} & \mathcal{A} \ar[d]^(0.3){p} & \\
& \mathcal{X} \ar[r]_{q} &\mathcal{X}^\uf}}}
\end{equation}
where the restrictions of $e$ and $q$ to each seed torus are given by
\[
e^*\left(A_{i;v}\right)=\left\{\begin{array}{ll}
    A_{i;v} & \text{if $i\in I^\uf$},  \\
    1 & \text{if $i\notin I^\uf$}, 
\end{array}\right. \quad \quad \text{and} \quad \quad q^*\left(X_{i;v}\right)=X_{i;v}.
\]
Note that although the map $p$ may not be algebraic, the composition $q\circ p$ and $p\circ e$ are both algebraic. In this paper we will mainly focus on some cluster variety $\mathcal{A}_{m,n}$ and its unfrozen counterpart $\mathcal{X}_{m,n}^\uf$, so to save energy we will abuse notation and denote the composition $q\circ p$ simply by $p$.

\subsection{Chiral Duality}

\begin{defn} \label{chiral dual} Given an initial seed $\vec{i}=\left(I,I^\uf,\epsilon\right)$ we define its \emph{chiral dual} to be $\vec{i}^\circ:=\left(I,I^\uf,-\epsilon\right)$. Let $\mathcal{A}^\circ$ and $\mathcal{X}^\circ$ be the cluster varieties corresponding to the chiral dual seed $\vec{i}^\circ$. There is an obvious one-to-one correspondence between seed tori on cluster varieties associated to $\vec{i}$ and those associated to $\vec{i}^\circ$. By using corresponding seed tori on the two cluster ensembles, we can define two involutive morphisms
\begin{align*}
    i_\mathcal{A}:\mathcal{A}\rightarrow & \mathcal{A}^\circ &  i_\mathcal{X}:\mathcal{X}\rightarrow & \mathcal{X}^\circ\\
    i_\mathcal{A}^*\left(A_{i;v}\right)=&A_{i;v}, & i_\mathcal{X}^*\left(X_{i;v}\right)=&\mathcal{X}_{i;v}^{-1}.
\end{align*}
\end{defn}

It is not hard to seed that a seed isomorphism $\vec{i}_v\overset{\sigma}{\cong}\vec{i}_w$ also gives rise to a seed isomorphism between their chiral duals $\vec{i}_v^\circ\overset{\sigma^\circ}{\cong} \vec{i}_w^\circ$. Correspondingly we get a pair of cluster transformations, $\sigma^\circ:\mathcal{A}^\circ\rightarrow \mathcal{A}^\circ$ and $\sigma^\circ:\mathcal{X}^\circ\rightarrow \mathcal{X}^\circ$, and they are intertwine to their counterparts via the morphisms $i_\mathcal{A}$ and $i_\mathcal{X}$; in other words, the following diagrams commute.
\[
\xymatrix{\mathcal{A} \ar[r]^\sigma \ar[d]_{i_\mathcal{A}} & \mathcal{A} \ar[d]^{i_\mathcal{A}}\\ \mathcal{A}^\circ \ar[r]_{\sigma^\circ} & \mathcal{A}^\circ} \quad \quad \quad \quad \xymatrix{\mathcal{X} \ar[r]^\sigma \ar[d]_{i_\mathcal{X}} & \mathcal{X} \ar[d]^{i_\mathcal{X}}\\ \mathcal{X}^\circ \ar[r]_{\sigma^\circ} & \mathcal{X}^\circ}
\]

\subsection{Canonical Quantization of a Cluster Poisson Variety}

In this subsection we recall the canonical quantization of a cluster Poisson variety as introduced by Fock and Goncharov in \cite{FGensemble}.

Let $\vec{i}=\left(I,I^\uf,\epsilon\right)$ be a seed. Let $\Lambda$ be a rank $|I|$ lattice together with a $\mathbb{Q}$-valued skew-symmetric form $\{\cdot, \cdot\}$ and a basis $\left\{e_i\right\}_{i\in I}$ such that $\left\{e_i,e_j\right\}=\epsilon_{ij}$. We then define the \emph{quantum torus} associated to $\vec{i}$ to be the following non-commutative division ring
\[
\mathbb{X}_\vec{i}=\Frac\frac{\mathbb{C}\left[q^{\pm \frac{1}{2}}\right]\left\langle X^v \ \middle| \ v\in \Lambda\right\rangle}{X^vX^w=q^{\{v,w\}}X^{v+w}}
\]

Note that in the classical limit $q\rightarrow 1$, $\mathbb{X}_\vec{i}$ recovers the field of rational functions on the seed torus $\mathcal{X}_\vec{i}$.

Suppose $\vec{i}'$ is the seed obtained from $\vec{i}$ via a mutation in the direction $k$. We keep the lattice $\Lambda$ and the skew-symmetric form $\{\cdot, \cdot\}$ unchanged and change the basis from $\left\{e_i\right\}_{i\in I}$ to $\left\{e'_i\right\}_{i\in I}$ where
\begin{equation}\label{e mutation}
e'_i= \left\{\begin{array}{ll}
    -e_k & \text{if $i=k$},\\
    e_i+\left[\epsilon_{ik}\right]_+e_k & \text{if $i\neq k$}.
\end{array}\right.
\end{equation}
It is not hard to verify that 
\[
\left\{e'_i,e'_j\right\}=\epsilon'_{ij}.
\] 
By going through the same construction we obtain another quantum torus $\mathbb{X}_{\vec{i}'}$, and there is a natural identification map
\[
\mu_k^\#:\mathbb{X}'_\vec{i}\rightarrow \mathbb{X}_\vec{i}.
\]

Next consider the \emph{quantum dilogarithm series}
\[
\Psi_q(X):=\prod_{a=1}^\infty \frac{1}{1+q^{2a-1}X}.
\]
One amazing fact is that, although $\Psi_q(X)$ itself is an infinite series, the adjoint action of $\Psi_q\left(X^{e_k}\right)$ on $\mathbb{X}_\vec{i}$ is actually rational; by a simple computation one can verify that
\begin{equation}\label{conjugation}
\Ad_{\Psi_q\left(X^{e_k}\right)}\left(X^v\right)=\prod_{a=1}^{\left|\left\{v,e_k\right\}\right|}\left(1+q^{(2a-1)\sgn\left\{v,e_k\right\}}X^{e_k}\right)^{-\sgn\left\{v,e_k\right\}}X^v.
\end{equation}

Using this adjoint action we define the \emph{quantum mutation map} $\mu_k^*:=\Ad_{\Psi_q\left(X^{e_k}\right)}\circ \mu_k^\#$.

For any path $\gamma=v\rightsquigarrow w$ on the mutation tree $\mathfrak{T}$ we can compose the quantum mutation maps backward according to the labeling of the edges and get a map $\mu_\gamma^*:\mathbb{X}_{w}\rightarrow \mathbb{X}_{v}$. We call the system of quantum tori $\left\{\mathbb{X}_v\right\}$ with maps $\left\{\mu_\gamma^*\right\}$ the \emph{canonical quantization} of the cluster Poisson variety $\mathcal{X}$, which we denote by $\mathbb{X}$. 

One should really think of $\mathbb{X}_\vec{i}$ as the quantum analogue of the field of rational functions on the seed torus $\mathcal{X}_\vec{i}$. In particular, by taking the classical limit $q\rightarrow 1$ and setting $X_i=X^{e_i}$, we see that we actually recover the mutation formula \eqref{x mutation} from the quantum mutation map $\mu_k^*$. Furthermore, the Poisson bracket on $\mathcal{X}$ can be recovered from the commutator on the canonical quantization $\mathbb{X}$ as
\[
\left\{X_i,X_j\right\}=\lim_{q\rightarrow 1}\frac{\left[X^{e_i},X^{e_j}\right]}{q-q^{-1}}
\]

\subsection{Tropicalization and Fock-Goncharov Duality Conjecture}\label{tropicalization}

In this subsection we recall tropicalizatoin of cluster variety introduced by Fock and Goncharov in \cite{FGensemble}, and state the Fock and Goncharov duality conjecture.

Given an algebraic torus $T$, we define \emph{positive rational functions} to be functions that are ratios of linear combinations of characters with positive integer coefficients. Positive rational functions on $T$ form a semifield under ordinary addition and ordinary multiplication, which we denote by $P(T)$. A rational function $f:T\dashrightarrow T'$ is said to be \emph{positive} if the restriction of the pull-back map $f^*$ to $P\left(T'\right)$ is a semifield homomorphism from $P\left(T'\right)$ to $P(T)$.

A \emph{positive space} is an algebraic variety $X$ together with an atlas of algebraic tori $\left\{T_\alpha\right\}$ and positive birational gluing maps $\left\{\phi_{\alpha\beta}:T_\alpha\dashrightarrow T_\beta\right\}$. Note that by definition $\phi_{\alpha\beta}^*$ are necessarily semifield isomorphisms, and hence we can identify all $P\left(T_\alpha\right)$ as $P(X)$ and call it the semifield of positive functions on $X$. A rational map $f:X\dashrightarrow Y$ between positive spaces is said to be \emph{positive} if the restriction of the pull-back map $f^*$ to $P(Y)$ is a semifield homomorphism from $P(Y)$ to $P(X)$.

Given a positive space $X$ and a semifield $S$, we define its space of $S$-points, denoted by $X(S)$, to be the set of semifield homomorphisms from $P(X)$ to $S$, i.e.,
\[
X(S):=\Hom_\text{semifield}\left(P(X),S\right).
\]
Given a positive rational map $f:X\dashrightarrow Y$ between positive spaces we can induce a map 
\[
f(S):X(S)\rightarrow Y(S)
\]
between the space of $S$-points by precomposing with the pull-back map $f^*$. Note that this is completely analogous to functor of points in algebraic geometry.

\begin{exmp} The simplest example of positive space is an algebraic torus $T$. Let $\Lambda$ be the character lattice of $T$ and let $\Lambda^*$ be the cocharacter lattice. Let's try to compute $T(S)$, the space of $S$-points for some semifield $S$. Since the semifield of rational functions $P(T)$ is freely generated by its characters as a $\mathbb{Z}_{\geq 0}$-module, it follows that
\[
T\left(S\right)=\Hom_\text{semifield} \left(P(X),S\right)=\Hom_\text{abelian group}\left(\Lambda,S\right)=\Lambda^*\otimes_\mathbb{Z}S.
\]

The story for general positive space $X$ only adds a small twist. From above we know that for each chart $T_\alpha$ we know that $T_\alpha(S)=\Lambda_\alpha^*\otimes_\mathbb{Z} S$; since gluing maps $\phi_{\alpha\beta}$ are positive birational maps, $\phi_{\alpha\beta}(S)$ are going to be bijections as maps between sets. Now we can identify all $T_\alpha(S)$ together via these bijections, which gives precisely $X(S)$. However, please be aware that $X(S)$ no longer possesses linear structure as each $\Lambda_\alpha\otimes_\mathbb{Z} S$.
\end{exmp}

From the definition of cluster varieties $\mathcal{A}$ and $\mathcal{X}$, it is obvious that both are positive spaces. Therefore for any semifield $S$, the spaces of $S$-points $\mathcal{A}(S)$ and $\mathcal{X}(S)$ make sense.

In this paper, we will mainly use two semifields, which are the two versions of \emph{tropical integers} $\mathbb{Z}^T:=\left(\mathbb{Z}, +, \max\right)$ and $\mathbb{Z}^t:=\left(\mathbb{Z},+, \min\right)$. These two semifields are isomorphic via the isomorphism $a\mapsto -a$. Therefore for any positive space $X$ we can identify $X\left(\mathbb{Z}^T\right)$ and $X\left(\mathbb{Z}^t\right)$ canonically; we call the resulting space the space of \emph{tropical points} and just denote it as $X\left(\mathbb{Z}^\trop\right)$.

The spaces of tropical points $\mathcal{A}\left(\mathbb{Z}^\trop\right)$ and $\mathcal{X}\left(\mathbb{Z}^\trop\right)$ are important in cluster theory because of the following conjecture.

\begin{conj}[Fock-Goncharov Duality Conjecture] The algebra of regular functions $\mathcal{O}(\mathcal{A})$ admits a basis parametrized by $\mathcal{X}\left(\mathbb{Z}^\trop\right)$ and the algebra of regular functions $\mathcal{O}(\mathcal{A})$ admits a basis parametrized by $\mathcal{A}\left(\mathbb{Z}^\trop\right)$.
\end{conj}

\begin{rmk} The above formulation of Fock-Goncharov duality conjecture is only for skew-symmetric cases. In general there is the notion of Langlands dual cluster varieties, and one needs to replace the spaces of tropical points $\mathcal{X}\left(\mathbb{Z}^\trop\right)$ and $\mathcal{A}\left(\mathbb{Z}^\trop\right)$ by their Langlands dual counterparts.
\end{rmk}

Fock and Goncharov first observe such a duality on decorated spaces of $G$-local systems in \cite{FGteich}. Gross, Hacking, Keel, and Kontsevich proved a weak version of the conjecture in \cite{GHKK}, which states that there are subalgebras of $\mathcal{O}(\mathcal{A})$ and $\mathcal{O}(\mathcal{X})$ admitting bases parametrized by subsets of $\mathcal{X}\left(\mathbb{Z}^\trop\right)$ and $\mathcal{A}\left(\mathbb{Z}^\trop\right)$. In the same paper Gross, Hacking, Keel, and Kontsevich also gave sufficient conditions on the cluster ensemble for the full Fock-Goncharov duality to hold, and one of them can be reformulated as follows.

\begin{thm}[Gross-Hack-Keel-Kontsevich, Proposition 8.28 in \cite{GHKK}]\label{sufficient condition} The Fock-Goncharov duality holds for the cluster ensemble $(\mathcal{A},\mathcal{X})$ when the following two conditions are satisfied:
\begin{enumerate}
    \item the cluster complex in $\mathcal{X}_\prin\left(\mathbb{R}^\trop\right)$ is large enough so that its convex hull contains the whole $\mathcal{X}_\prin\left(\mathbb{R}^\trop\right)$ (under an identification with a linear space using an initial seed);
    \item the map $p:\mathcal{A}\rightarrow \mathcal{X}^\uf$ is surjective (see \eqref{uf} for the definition of this maps).
\end{enumerate}
\end{thm}

We will see in Proposition \ref{dt=sufficient condition} that condition (1) is implied by the existence of cluster Donaldson-Thomas transfromation, so constructing cluster Donaldson-Thomas transformation is a key to proving the Fock-Goncharov duality conjecture.

\subsection{Cluster Donaldson-Thomas Transformation}

We still owe the readers the definition of cluster Donaldson-Thomas transformation. Recall from Definition \ref{cluster transformation} that a cluster transformation $\sigma$ on $\mathcal{X}$ is given by a seed isomorphism between two seeds $\vec{i}_v$ and $\vec{i}_w$ and is given in terms of cluster variables as $\sigma^*\left(X_{i;w}\right)=X_{i;v}$; if one wants to write the cluster transformation $\sigma$ via its pull-back on the same seed torus, we get the following expression
\[
\sigma^*\left(X_{i;w}\right)=\mu^*_{w\rightsquigarrow v} \left(X_{\sigma(i);v}\right).
\]
Since the rational expression $\sigma^*\left(X_{i;w}\right)$ in terms of the cluster variables $\left(X_{j;w}\right)$ is the same as the pull-back via some cluster mutation sequence, it makes sense to talk about the $c$-matrix of the cluster transformation of $\sigma$ with respect to any seed $\vec{i}$.

\begin{defn} A cluster transformation $\dt$ on $\mathcal{X}^\uf$ is called a \emph{cluster Donaldson-Thomas transformation} if its $c$-matrix with respect to some seed $\vec{i}$ is $-\mathrm{Id}$.
\end{defn}

\begin{prop}\label{dt=sufficient condition} The existence of cluster Donaldson-Thomas transformation implies condition (1) in Theorem \ref{sufficient condition}.
\end{prop}
\begin{proof} From part (5) of Theorem \ref{cgf} we know that when the $c$-matrix equals to $-\mathrm{Id}$, so does the $g$-matrix; but then by part (3) of Theorem \ref{cgf} we know that the cluster chamber corresponding to the $g$-matrix that is $-\mathrm{Id}$ must be the opposite chamber. But then this implies that the opposite chamber is part of the cluster complex and hence the convex hull of the cluster complex must equal to all of $\mathcal{X}_\prin\left(\mathbb{R}^\trop\right)$.
\end{proof}

Recall from Equation \eqref{cf} that we can write the pull-backs against cluster Donaldson-Thomas transformation as a rational expression in terms of the $c$-matrix and $F$-polynomials, and since $F$-polynomials have constant term 1, the above definition of cluster Donaldson-Thomas transformation is equivalent to the following one.

\begin{defn}\label{dt} A cluster transformation $\dt$ on $\mathcal{X}^\uf$ is called a \emph{cluster Donaldson-Thomas transformation} if on some (equivalently any) seed torus $\mathcal{X}_v$,
\[
\ord_{X_{i;v}}\dt^*\left(X_{j;v}\right)=-\delta_{ij},
\]
where $\ord_zf$ is defined to be the lowest power of the variable $z$ if $f$ is a polynomial, and $\ord_z\frac{f}{g}:=\ord_zf-\ord_zg$ for a rational function $f/g$.
\end{defn}

Recall from Definition \ref{chiral dual} that any cluster transformation on $\mathcal{X}^\uf$ also gives rise to a cluster transformation on $\left(\mathcal{X}^\circ\right)^\uf$. We will denote the cluster transformation on $\left(\mathcal{X}^\circ\right)^\uf$ that is the chiral dual of the cluster Donaldson-Thomas transformation $\dt$ by $\DT$, and call it the cluster Donaldson-Thomas transformation as well. One should understand $\DT$ is the analogue of $\dt$ with respect to the tropical integer $\mathbb{Z}^T$ instead of $\mathbb{Z}^t$. In particular, there is an analogous degree condition similar to that of Definition \ref{dt}.

\begin{prop}\label{degree condition for DT} A cluster transformation $\DT$ on $\left(\mathcal{X}^\circ\right)^\uf$ is a cluster Donaldson-Thomas transformation if and only if for some (and equivalently any) seed torus $\mathcal{X}_v^\circ$,
\[
\deg_{X_{i;v}^\circ}\DT^*\left(X_{j;v}^\circ\right)=-\delta_{ij}.
\]
\end{prop}
\begin{proof} Note that for any polynomial $f(z)$, $\deg_z f\left(z^{-1}\right)=-\ord_zf(z)$. Therefore if $\dt^*\left(X_{j;v}\right)=\frac{f\left(X_{1;v},\dots, X_{n;v}\right)}{g\left(X_{1;v}, \dots, X_{n;v}\right)}$ for polynomials $f$ and $g$, then we have
\begin{align*}
\deg_{X_{i;v}^\circ}\DT^*\left(X_{j;v}^\circ\right)=&\deg_{X_{i;v}^\circ}g\left(\frac{1}{X_{1;v}^\circ}, \dots, \frac{1}{X_{n;v}^\circ}\right)-\deg_{X_{i;v}^\circ}f\left(\frac{1}{X_{1;v}^\circ}, \dots, \frac{1}{X_{n;v}^\circ}\right)\\
=&-\ord_{X_{i;v}}g\left(X_{1;v},\dots, X_{n;v}\right)+\ord_{X_{i;v}}f\left(X_{1;v}, \dots, X_{n;v}\right)\\
=&\ord_{X_{i;v}}\dt^*\left(X_{j;v}\right)\\
=&-\delta_{ij}. \qedhere
\end{align*}
\end{proof}

By definition the two versions of cluster Donaldson-Thomas transformation are intertwined by the chiral duality map $i_{\mathcal{X}}$, so knowing the construction of one of them automatically gives a construction for the other.

\section{Decorated Grassmannian and Configuration Space}

\label{section4}

In this section we introduce the main geometric spaces we consider in this paper, namely the decorated Grassmannian and configuration spaces. These two spaces are birationally equivalent to a pair of cluster varieties $\mathcal{A}_{m,n}$ and $\mathcal{X}_{m,n}^\uf$ respectively. As we will see in the next section that these birational equivalences will provide us a factorization of one version of the cluster Donaldson-Thomas transformation and serve as key ingredient for our main result.

\subsection{Decorated Grassmannian and \texorpdfstring{$\mathcal{A}_{m,n}$}{}}

Let's begin by describing the cluster variety structure on an open subset of the affine cone over the Grassmannian, and the main reference we will follow is \cite{Sco} and \cite{Gon}.

Fix a pair of integers $1<m<n-1$. Let $W$ be an $n$-dimensional vector space and denote its dual space as $W^*$. The \emph{Grassmannian} $\Gr_m\left(W^*\right)$ is defined to be the moduli space of $m$-dimensional subspace in $W^*$. We define the \emph{decorated Grassmannian} to be
\[
\dGr_m\left(W^*\right):=\left\{\left(V^*,\alpha\right)\ \middle| \ \begin{array}{l} \text{$V^*\subset W^*$ is an $m$-dimensional subspace}\\ \text{and $\alpha$ is a non-zero element in $\bigwedge^m V^*$}\end{array}\right\}.
\]
There is an obvious forgetful map $\dGr_m\left(W^*\right)\rightarrow \Gr_m\left(W^*\right)$.

For any element $g\in \bigwedge^m W$ we can define a regular function $\Delta_g$ on $\dGr_m\left(W^*\right)$ by setting
\[
\Delta_g\left(V^*,\alpha\right):=\inprod{\alpha}{g}.
\]

Fix a basis $\left\{e_i\right\}_{i=1}^n$ of $W$. Then for any $m$-element subset $I=\left\{i_1<i_2<\dots<i_m\right\}$ of $\{1,\dots, m\}$ we get a regular function 
\[
\Delta_I:=\Delta_{e_{i_1}\wedge e_{i_2}\wedge \dots \wedge e_{i_m}}
\]
on $\dGr_m\left(W^*\right)$; such regular functions are known as the \emph{Pl\"{u}cker coordinates}. It is known that the ring of regular functions $\mathcal{O}\left(\dGr_m\left(W^*\right)\right)$ is generated by Pl\"{u}cker coordinates, with relations generated by the 3-term \emph{Pl\"{u}cker relations}:
\begin{equation}\label{plucker relation}
\Delta_{J\cup \{i,l\}}\Delta_{J\cup\{j,k\}}+\Delta_{J\cup \{i,j\}}\Delta_{J\cup \{k,l\}}=\Delta_{J\cup \{i,k\}}\Delta_{J\cup \{j,l\}}
\end{equation}
for any $i<j<k<l$ and any $(m-2)$-element subset $J\subset \{1,\dots, n\}\setminus \{i,j,k,l\}$.

Let $D_i$ be the vanishing locus of the regular function $\Delta_{[i-m+1,i]}$ on $\Gr_m\left(W^*\right)$ (see Notation \ref{bracket} for the notation $[i-m+1,i]$). It is known that $D:=\sum_{i=1}^n D_i$ is an anticanonical divisor $\Gr_m\left(W^*\right)$. Lastly, we define the open subset
\[
\Gr_m^\times\left(W^*\right):=\Gr_m\left(W^*\right)\setminus D.
\]
and define $\dGr_m^\times \left(W^*\right)$ to be the preimage of $\Gr_m^\times \left(W^*\right)$ against the projection map $\dGr_m\left(W^*\right)\rightarrow \Gr_m\left(W^*\right)$. In particular, it follows that $\mathcal{O}\left(\dGr_m^\times\left(W^*\right)\right)$ can be obtained from $\mathcal{O}\left(\dGr_m\left(W^*\right)\right)$ by localizing at Pl\"{u}cker coordinates $\Delta_{[i,i+m-1]}$ for all $1\leq i\leq n$.

For the rest of this section we will describe a birational equivalence between $\dGr_m^\times\left(W^*\right)$ and an cluster variety $\mathcal{A}_{m,n}$ which is an isomorphism on a subset up to codimension 2; in particular, this implies that 
\[
\mathcal{O}\left(\dGr_m^\times \left(W^*\right)\right)\cong \mathcal{O}\left(\mathcal{A}_{m,n}\right).
\]

Let $\Gamma$ be a bipartite graph on a disk with $n$ boundary marked points and of full rank $m$. Recall that from $\Gamma$ we obtain a quiver $\vec{i}_\Gamma$ and a dominating set $I(f)$, which is of size $m$, for each face $f$. Since quiver carries the same data as a seed, we get a seed torus $\mathcal{A}_\Gamma:=\mathcal{A}_{\vec{i}_\Gamma}$. We then define a rational map \begin{align*}
    \psi_\Gamma:\dGr_m^\times \left(W^*\right)&\dashrightarrow \mathcal{A}_{\Gamma}\\
    \psi_\Gamma\left(A_f\right)&=\Delta_{I(f)}.
\end{align*}

\begin{prop} As $\Gamma$ varies through all bipartite graph of full rank $m$ on the disk with $n$ boundary marked points, the rational maps $\psi_\Gamma$ can be glued into a rational map $\psi:\dGr_m^\times\left(W^*\right)\dashrightarrow \mathcal{A}_{m,n}$,
where $\mathcal{A}_{m,n}$ is the cluster variety defined by the seed $\vec{i}_\Gamma$ for any bipartite graph $\Gamma$ of full rank $m$ on the disk with $n$ boundary marked points.
\end{prop}
\begin{proof} We know from Thurston's theorem \ref{thurston} that any two bipartite graphs with the same parameters $m$ and $n$ are related by a sequence of 2-by-2 moves, and we know from Proposition \ref{bipartite graph mutation} that each 2-by-2 move gives a mutation between the corresponding quivers. So now all we need to show is that if $\Gamma$ and $\Gamma'$ are related by a single 2-by-2 move at a non-boundary face $f$, then $\psi_\Gamma$ and $\psi_{\Gamma'}$ fit in the following commutative diagram:
\[
\vcenter{\vbox{\xymatrix{ & \mathcal{A}_\Gamma \ar@{-->}[dd]^{\mu_f} \\ \dGr_m^\times\left(W^*\right) \ar@{-->}[ur]^{\psi_\Gamma} \ar@{-->}[dr]_{\psi_{\Gamma'}} & \\
& \mathcal{A}_{\Gamma'}}}}.
\]
But this follows from a comparison between \eqref{plucker relation}, \eqref{a mutation}, and Proposition \ref{bipartite graph mutation}.
\end{proof}

Scott proved that $\mathcal{O}\left(\dGr_m\left(W^*\right)\right)\cong \ord(\mathcal{A}_{m,n})$ in \cite{Sco}, by modifying his proof we can prove the following statement.

\begin{prop}\label{psi isomorphism} The pull-back of any regular function on $\mathcal{O}(\mathcal{A}_{m,n})$ to $\dGr_m^\times\left(W^*\right)$ via $\psi^*$ is regular, and 
\[
\psi^*:\mathcal{O}(\mathcal{A}_{m,n})\rightarrow \mathcal{O}\left(\dGr_m^\times\left(W^*\right)\right)
\]
is an algebra isomorphism. In particular, this implies that the birational equivalence $\psi:\dGr_m^\times\left(W^*\right)\dashrightarrow \mathcal{A}_{m,n}$ is an isomorphism between the complements of a subset of codimension at least 2. 
\end{prop}
\begin{proof} Consider the chessboard bipartite graph $\Gamma_\cb$ (Example \ref{chessboard}). For each non-boundary face $f$ of $\Gamma_\cb$ we define $\Gamma_\cb^f$ to be the bipartite graph obtained by applying a 2-by-2 move to the face $f$. Let $U_{\Gamma_\cb}$ be the preimage $\psi^{-1}\left(\mathcal{A}_{\Gamma_\cb}\right)$, and let $U_{\Gamma_\cb^f}$ be the preimage $\psi^{-1}\left(\mathcal{A}_{\Gamma_\cb^f}\right)$. Define
\[
U:=U_{\Gamma_\cb}\cup\bigcup_f U_{\Gamma_\cb^f}\subset \dGr_m^\times\left(W^*\right).
\]
We claim that the complement of $\dGr_m^\times\left(W^*\right)$ is of codimension at least 2. Note that in order to be in the complement, it has to lie in the intersection of the vanishing loci of $\Delta_{I(f);\Gamma_\cb}$ and $\Delta_{I(f);\Gamma_\cb^f}$ for some non-boundary face $f$; but since $\Delta_{I(f);\Gamma_\cb}$ and $\Delta_{I(f);\Gamma_\cb^f}$ are coprime in $\mathcal{O}\left(\dGr_m^\times\left(W^*\right)\right)$, it follows that the codimension of the complement is at least 2.

Now if $F$ is a regular function on $\mathcal{A}_{m,n}$, then it is regular on $\mathcal{A}_{\Gamma_\cb}\cup \bigcup_f\mathcal{A}_{\Gamma_\cb^f}$, and hence its pull-back $\psi^*(F)$ is regular on $U$. Recall a fact from algebraic geometry that if a function is regular if and only if it is regular outside of a codimension no less than 2 subset; therefore we can conclude that $\psi^*(F)$ is regular on $\dGr_m^\times\left(W^*\right)$.

To see that $\psi^*:\mathcal{O}(\mathcal{A}_{m,n})\rightarrow \mathcal{O}\left(\dGr_m^\times\left(W^*\right)\right)$ is an algebra isomorphism, we first observe that $\psi^*$ is surjective because all generators of $\mathcal{O}\left(\dGr_m^\times\left(W^*\right)\right)$ are in the image of $\psi^*$ (Proposition \ref{any I}). But then this implies that $\psi:U\rightarrow \mathcal{A}_{m,n}$ is an embedding. On the other hand we know that $\mathcal{A}_{m,n}$ is irreducible, and since 
\[
\dim U=\dim \dGr^\times_m\left(W^*\right)=m(n-m)+1=(m-1)(n-m-1)+n=\dim \mathcal{A}_{m,n},
\]
we know that the embedding has to be dominant. From this we know that $\psi^*$ is injective; combining with the surjectivity proved above we get that $\psi^*$ is an algebra isomorphism.
\end{proof}

There is another way to think of the open set $\dGr_m^\times\left(W^*\right)$. The choice of basis $\left\{e_i\right\}$ for $W$ defines a dual basis $\left\{e_i^*\right\}$ for $W^*$; this gives an isomorphism $W^*\cong \mathbb{C}^n$, and we can represent an element $\left(V^*,\alpha\right)\in \dGr_m\left(W^*\right)$ by an $m\times n$ matrix $M_{\left(V^*,\alpha\right)}$ over $\mathbb{C}$ with row vectors $\alpha_i$ such that $\alpha=\alpha_1\wedge \dots \wedge \alpha_m$. It is not hard to see that the choices of such $\alpha_i$ are unique up to a left $\SL_m$-action, and hence we can deduce that
\[
\dGr_m^\times\left(W^*\right)\cong \SL_m\left\backslash \mat_{m,n}^\times\right.,
\]
where $\mat_{m,n}^\times$ denotes the set of $m\times n$ matrices whose column vectors $v_i,\dots, v_{i+m-1}$ (indices modulo $n$) are in general position. Now for any $m$-element subset $I\subset\{1,\dots, n\}$, we see that $\Delta_I\left(V^*,\alpha\right)$ is precisely the minor of $M_{\left(V^*,\alpha\right)}$ associated to the column vectors labeled by elements of $I$.

\subsection{Configuration Space}

\begin{defn} For an $m$-dimensional vector space $V$, we define the \emph{configuration space} of $n$ non-zero vectors in $V$ to be
\[
\conf_n\left(V\right):=\GL(V)\left\backslash \left(V\setminus \vec{0}\right)^n\right.,
\]
where $\GL(V)$ acts diagonally on all the $n$ vectors. There is a projective counterpart of such space, which can be defined as
\[
\conf_n\left(\mathbb{P}V\right):=\PGL(V)\left\backslash\left(\mathbb{P}V\right)^n\right..
\]
There is an obvious projection map $\conf_n(V)\rightarrow \conf_n(\mathbb{P}V)$
\end{defn}

One important fact about configuration spaces is that they are canonically isomorphic, as stated in the proposition below.

\begin{prop}\label{canonical isomorphism} If $V$ and $V'$ are two vector spaces over the same base field and $\dim V=\dim V'<\infty$, then $\conf_n(V)\cong \conf_n\left(V'\right)$ and $\conf_n\left(\mathbb{P}V\right)\cong \conf_n\left(\mathbb{P}V'\right)$ canonically. As a corollary, we can canonically identify $\conf_n(V)$ with $\conf_n\left(\mathbb{C}^m\right)$, and canonically identify $\conf_n(\mathbb{P}V)$ with $\conf_n\left(\mathbb{P}^{m-1}\right)$, for any $m$-dimensional vector space $V$.
\end{prop}
\begin{proof} Recall that any two vector spaces of the same finite dimensions are isomorphic, but the isomorphism is not canonical because it involves a choice of bases for the two vector spaces. However, since we quotient out the automorphism groups when we pass to configuration spaces, such isomorphisms become canonical.
\end{proof}

\begin{notn} Analogous to the case of decorated Grassmannian, we use a superscript $\times$ to denote the open subset where any $n$ consecutive vectors (resp. lines) are in general position. It is not hard to see that the projection map $\conf_n(V)\rightarrow \conf_n(\mathbb{P}V)$ restricts to a projection $\conf_n^\times(V)\rightarrow \conf_n^\times(\mathbb{P}V)$. 
\end{notn}

For $\left(V^*,\alpha\right)\in \dGr_m^\times\left(W^*\right)$, we let $V$ be the dual space of $V^*$. Then we have a projection map $W\rightarrow V$ dual to the inclusion map $V^*\rightarrow W^*$. We can project our basis vectors $\left\{e_i\right\}$ of $W$ to $V$ and denote their images as $\left\{\overline{e}_i\right\}$. Due to the condition that $\Delta_{[i,i+m-1]}\left(V^*,\alpha\right)\neq 0$, we know the configuration $\left[\overline{e}_1,\dots, \overline{e}_n\right]$ is in $\conf_n^\times (V)\cong \conf_n^\times\left(\mathbb{C}^m\right)$. Furthermore, since the map $\left(V^*,\alpha\right)\mapsto \left[\overline{e}_1,\dots, \overline{e}_n\right]$ does not depend on the decoration $\alpha$, it factors through the ordinary Grassmannian as $\dGr_m^\times\left(W^*\right)\rightarrow \Gr_m^\times\left(W^*\right)\rightarrow \conf_n^\times\left(\mathbb{C}^m\right)$.

\begin{prop} The map $V^*\mapsto \left[\overline{e}_1,\dots, \overline{e}_n\right]$ is an isomorphism between $\Gr_m^\times\left(W^*\right)$ and $\conf_n^\times\left(\mathbb{C}^m\right)$.
\end{prop}
\begin{proof} Similar to $\dGr_m^\times\left(W^*\right)\cong \SL_m\left\backslash \mat_{m,n}^\times\right.$, we have $\Gr_m^\times\left(W^*\right)\cong \GL_m\left\backslash \mat_{m,n}^\times\right.$, which is by definition $\conf_n^\times\left(\mathbb{C}^m\right)$. To see that the map $V^*\mapsto \left[\overline{e}_1,\dots, \overline{e}_n\right]$ actually provides such an isomorphism, fix a matrix representative $M_{\left(V^*,\alpha\right)}=\left(m_{ij}\right)$. Then we know that $\alpha_i:=\sum_jm_{ij}e_j^*$ is a basis of $V^*$ and hence the images $\overline{e}_i$ in $V$ can be expressed as 
\[
\overline{e}_i=\sum_j\inprod{\alpha_j}{e_i}\alpha_j^*=\sum_jm_{ji}\alpha_j^*.
\]
In other words, the column vectors of $M_{\left(V^*,\alpha\right)}$ are precisely $\overline{e}_i$ expressed with respect to the basis $\left\{\alpha_j^*\right\}$ of $V$; but this is exactly how we think of the configuration space $\conf_n^\times\left(\mathbb{C}^m\right)$ in terms of column vectors of matrices, and hence the claim is true.
\end{proof}

Last but not least, by composing the projection map $\dGr_m^\times\left(W^*\right)\rightarrow \Gr_m^\times\left(W^*\right)\cong \conf_n^\times\left(\mathbb{C}^m\right)$ with another projection map $\conf_n^\times\left(\mathbb{C}^m\right)\rightarrow \conf_n^\times\left(\mathbb{P}^{m-1}\right)$ we get a surjective map
\[
\pi:\dGr_m^\times\left(W^*\right)\rightarrow \conf_n^\times\left(\mathbb{P}^{m-1}\right).
\]

\subsection{Quantum Matrices and the Canonical Quantization of \texorpdfstring{$\mathcal{X}_{m,n}$}{}}\label{quantum}

Following the notation we use in this paper, we denote the cluster Poisson variety defined by any bipartite graph of full rank $m$ on a disk with $n$ boundary marked points as $\mathcal{X}_{m,n}$, and its corresponding canonical quantization as $\mathbb{X}_{m,n}$. In this section we relate the canonical quantization $\mathbb{X}_{m,n}$ to quantum matrices $\mathcal{O}_q\left(\mat_{m,n}\right)$, which is a $q$-deformation of the algebra of regular functions on the space of $m\times n$ matrices.

Let's first recall the definition of quantum $m\times (n-m)$ matrices $\mathcal{O}_q\left(\mat_{m,n-m}\right)$. Consider a formal $m\times (n-m)$ matrix with entries $m_{ij}$. The quantum $m\times (n-m)$ matrices is the associative algebra generated by these entries $m_{ij}$ over the field $\mathbb{C}(q)$ subject to the following relations whenever $i<k$ and $j<l$:
\begin{equation}\label{quantum matrix relations} \begin{split}
    m_{ij}m_{il}&=qm_{il}m_{ij}, \\
    m_{ij}m_{kj}&=qm_{kj}m_{ij}, \\
    m_{il}m_{kj}&=m_{kj}m_{il}, \\
    m_{ij}m_{kl}&=m_{kl}m_{ij}+\left(q-q^{-1}\right)m_{il}m_{kj}.
\end{split}
\end{equation}

By a ``quantum rational map'' from $\mathbb{X}_{m,n}$ to the quantum $m\times (n-m)$ matrices, we mean a collection of maps $\chi_\Gamma^*:\mathcal{O}_q\left(\mat_{m,n-m}\right)\rightarrow \mathbb{X}_\Gamma$ such that the following diagram commutes for any 2-by-2 move that takes place at a non-boundary face $f$, which turns a bipartite graph $\Gamma$ into $\Gamma'$:
\begin{equation}\label{quantum diagram}
\vcenter{\vbox{\xymatrix{& \mathbb{X}_\Gamma   \\ \mathcal{O}_q\left(\mat_{m,n-m}\right) \ar[ur]^{\chi_\Gamma^*}  \ar[dr]_{\chi^*_{\Gamma'}} & & \\
& \mathbb{X}_{\Gamma'} \ar[uu]_{\mu_f^*}}}}.
\end{equation}

Let's now define such $\chi_\Gamma^*$. Let $\Gamma$ be a bipartite graph of full rank $m$ on a disk with $n$ boundary marked points. Draw the standard perfect orientation on $\Gamma$. Change the labellings of the boundary marked points according to
\[
i\mapsto \left\{\begin{array}{ll}
    m+i+1 & \text{if $i\leq m$}, \\
    (i-m)' & \text{if $i>m$}.
\end{array}\right.
\]
By doing such a change, all sources of the standard perfect orientation are denoted with unprimed indices and all sinks of the standard perfect orientation are denoted with primed indices. Now for a path $\gamma$ on $\Gamma$ going from a source $i$ to a sink $j'$, which we denote as $\gamma:i\rightarrow j'$, we define $\hat{\gamma}$ to be the collection of faces lying on the right hand side of $\gamma$ with respect to its orientation.
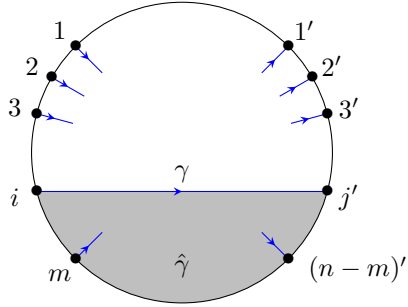
\[
\tikz{
\fill [fill=lightgray] +(195:2) arc (195:345:2) -- cycle;
\node at (0,-1.5) [] {$\hat{\gamma}$};
\draw (0,0) circle [radius=2];
\foreach \i in {1,...,3}
    {
    \draw [blue, ->-] +(120+15*\i:2) to +(120+15*\i:1.5);
    \draw [blue, ->-] +(60-15*\i:1.5) to +(60-15*\i:2);
    \node at +(120+15*\i:2) [] {$\bullet$};
    \node at +(60-15*\i:2) [] {$\bullet$};
    \node at +(120+15*\i:2.3) [] {$\i$};
    \node at +(60-15*\i:2.3) [] {$\i'$};
    }
\draw [blue, ->-] +(225:2) to +(225:1.5);
\draw [blue, ->-] +(-45:1.5) to + (-45:2);
\draw [blue, ->-] +(195:2) to node [black, above] {$\gamma$} + (-15:2);
\node at +(195:2) [] {$\bullet$};
\node at +(-15:2) [] {$\bullet$};
\node at +(225:2) [] {$\bullet$};
\node at +(-45:2) [] {$\bullet$};
\node at +(195:2.3) [] {$i$};
\node at +(-15:2.3) [] {$j'$};
\node at +(225:2.3) [] {$m$};
\node at +(-45:2.2) [right] {$(n-m)'$};
}
\]
Then we define
\[
\chi^*_\Gamma\left(m_{ij}\right):=\sum_{\gamma:i\rightarrow j'}X^{\sum_{f\in \hat{\gamma}}e_f}.
\]

\begin{prop} $\chi^*_\Gamma$ is an algebra homomorphisms.
\end{prop}
\begin{proof} To prove this proposition, we need to verify that the relations \eqref{quantum matrix relations} still hold with $m_{ij}$ replaced by their images $\chi^*_\Gamma\left(m_{ij}\right)$. The key is to show the following three identites, where each letter represents the sum of the basis vectors corresponding to the faces in the respective connected region in the schematic picture below (note that each blue arrow is part of a path on $\Gamma$ compatible with the standard perfect orientation, and there may be other vertices and edges which we have not drawn explicitly).
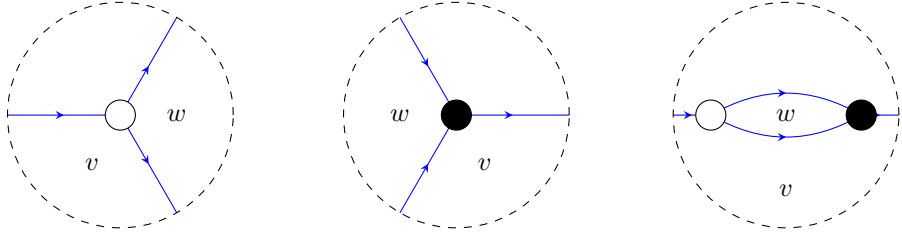
\[
\begin{tikzpicture}[baseline=-0.5ex]
\draw [blue, ->-] (-1.5,0) -- (0,0);
\draw [blue, ->-] (0,0) to +(60:1.5);
\draw [blue, ->-] (0,0) to +(-60:1.5);
\draw [dashed] (0,0) circle [radius=1.5];
\draw [fill=white] (0,0) circle [radius=0.2];
\node at +(-120:0.75) [] {$v$};
\node at +(0.75,0) [] {$w$};
\node at (0,-2.5) [] {$X^{v+w}X^v=qX^vX^{v+w},$};
\end{tikzpicture}
\quad\quad
\begin{tikzpicture}[baseline=-0.5ex]
\draw [blue, ->-] +(120:1.5) to (0,0);
\draw [blue,->-] +(-120:1.5) to (0,0);
\draw [blue,->-] (0,0) -- (1.5,0);
\draw [dashed] (0,0) circle [radius=1.5];
\draw [fill] (0,0) circle [radius=0.2];
\node at (-0.75,0) [] {$w$};
\node at +(-60:0.75) [] {$v$};
\node at (0,-2.5) [] {$X^{v+w}X^v=qX^vX^{v+w},$};
\end{tikzpicture}
\quad \quad
\begin{tikzpicture}[baseline=-0.5ex]
\draw [blue,->-] (-1.5,0) -- (-1,0);
\draw [blue, ->-] (1,0) -- (1.5,0);
\draw [blue, ->-] (-1,0) to [bend left] (1,0);
\draw [blue, ->-] (-1,0) to [bend right] (1,0);
\draw [dashed] (0,0) circle [radius=1.5];
\draw [fill] (1,0) circle [radius=0.2];
\draw [fill=white] (-1,0) circle [radius=0.2];
\node at (0,0) [] {$w$};
\node at (0,-1) [] {$v$};
\node at (0,-2.5) {$X^{v+w}X^v=X^vX^{v+w},$};
\end{tikzpicture}
\]
From the $q$-commutation relation defined on $\mathbb{X}_\Gamma$ we know that the left hand side of these three identities are equal to $q^{2\{w,v\}}X^vX^{v+w}$; therefore these three identities come down to computing the skew-symmetric product $\{w,v\}$. In the first case we see that the center white vertex contributes a weight $\frac{1}{2}$ from $w$ to $v$, and every additional vertex along the lower right blue arrow contribute two arrows with weight $\frac{1}{2}$ in the opposite directions; therefore we can conclude that $\{w,v\}=\frac{1}{2}$. By similar computation we get that $\{w,v\}=\frac{1}{2}$ in the second case and $\{w,v\}=0$ in the last case.

Since that the skew-symmetric form $\{\cdot, \cdot\}$ are defined locally, for any pair of paths $\gamma$ and $\eta$ that we want to commute, we just need to paste these identities together to find the right identity. To demonstrate we will do one example here. Consider the following two paths in $\Gamma$ that are compatible with the standard perfect orientation (we use blue for $\gamma$ and red for $\eta$).
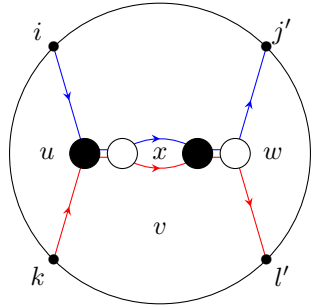
\[
\tikz{
\draw [->-, blue] +(135:2) -- (-1,0);
\draw [->-, red] +(-135:2) -- (-1,0);
\draw [ blue] (-1,0.05) -- (-0.5,0.05);
\draw [red] (-1,-0.05) -- (-0.5,-0.05);
\draw [->-, red] (-0.5,-0.05) to [bend right] (0.5,-0.05);
\draw [->-, blue] (-0.5,0.05) to [bend left] (0.5,0.05);
\draw [ blue] (1,0.05) -- (0.5,0.05);
\draw [red] (1,-0.05) -- (0.5,-0.05);
\draw [blue,->-] +(0:1) to +(45:2);
\draw [red, ->-] +(0:1) to +(-45:2);
\draw (0,0) circle [radius=2];
\draw [fill] (-1,0) circle [radius=0.2];
\draw [fill=white] (-0.5,0) circle [radius=0.2];
\draw [fill] (0.5,0) circle [radius=0.2];
\draw [fill=white] (1,0) circle [radius=0.2];
\node at (-1.5,0) [] {$u$};
\node at (0,-1) [] {$v$};
\node at (1.5,0) [] {$w$};
\node at (0,0) [] {$x$};
\foreach \i in {0,...,3}
	{
	\node at (45-90*\i:2) [] {$\bullet$};
	}
\node at (45:2.3) [] {$j'$};
\node at (-45:2.3) [] {$l'$};
\node at (135:2.3) [] {$i$};
\node at (-135:2.3) [] {$k$};
}
\]
By using the three identities above (in fact, the first two are sufficient), we see that
\begin{align*}X^{u+v+w+x}X^v=&qX^{u+v+x}X^{v+w} \\
=&qX^{u+v+x}X^{v+w}+X^vX^{u+v+w+x}-X^vX^{u+v+w+x}\\
=&qX^{u+v+x}X^{v+w}+X^vX^{u+v+w+x}-q^{-1}X^{u+v}X^{v+w+x};
\end{align*} 
note that at the first step we have applied the first identity above to the white vertex on the farthest right and at the last step we have applied the second identity to the black vertex on the farthest left. This may not look like what we expect yet; but we can use the last identity above to show that $X^{u+v}X^{v+w+x}=X^{u+v+x}X^{v+w}$. Therefore we may conclude that
\[
X^{u+v+w+x}X^v=X^vX^{u+v+w+x}+\left(q-q^{-1}\right)X^{u+v+x}X^{v+w},
\]
which gives one term in the identity $\chi_\Gamma^*\left(m_{ij}\right)\chi_\Gamma^*\left(m_{kl}\right)=\chi_\Gamma^*\left(m_{kl}\right)\chi_\Gamma^*\left(m_{ij}\right)+\left(q-q^{-1}\right)\chi_\Gamma^*\left(m_{il}\right)\chi_\Gamma^*\left(m_{jk}\right)$ that we want to show. All other identities can be proved in a similar way, and we will leave the details as an exercise to the readers.
\end{proof}

\begin{prop} The diagram \eqref{quantum diagram} commutes.
\end{prop}
\begin{proof} Since a 2-by-2 move acts only locally, it suffices to consider the parts of $\Gamma$ and $\Gamma'$ near the mutating face $f$.
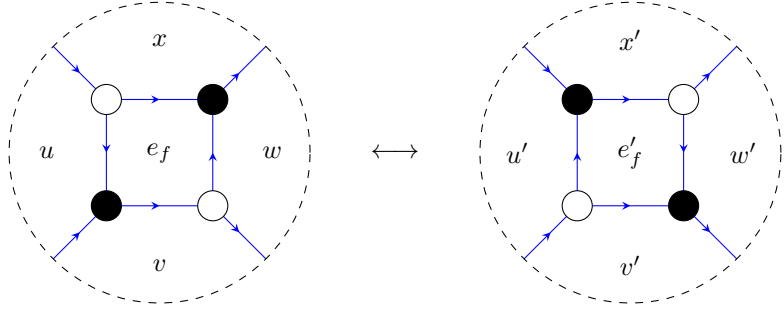
\[
\begin{tikzpicture}[baseline=-0.5ex]
\foreach \i in {1,...,4}
	{
	\coordinate (o\i) at +(225-90*\i:2);
	\coordinate (i\i) at +(225-90*\i:1);
	}
\draw [->-,blue] (o1) -- (i1);
\draw [->-,blue] (i2) -- (o2);
\draw [->-,blue] (i3) -- (o3);
\draw [->-,blue] (o4) -- (i4);
\draw [->-, blue](i1) -- (i2);
\draw [->-,blue] (i1) -- (i4);
\draw [->-,blue] (i3) -- (i2);
\draw [->-,blue] (i4) -- (i3);
\draw [fill=white] (i1) circle [radius=0.2];
\draw [fill=white] (i3) circle [radius=0.2];
\draw [fill] (i2) circle [radius=0.2];
\draw [fill] (i4) circle [radius=0.2];
\draw [dashed] (0,0) circle [radius=2]; 
\node at (-1.5,0) [] {$u$};
\node at (0,-1.5) [] {$v$};
\node at (1.5,0) [] {$w$};
\node at (0,1.5) [] {$x$};
\node at (0,0) [] {$e_f$};
\end{tikzpicture}\quad \quad \longleftrightarrow \quad \quad
\begin{tikzpicture}[baseline=-0.5ex]
\foreach \i in {1,...,4}
	{
	\coordinate (o\i) at +(225-90*\i:2);
	\coordinate (i\i) at +(225-90*\i:1);
	}
\draw [->-,blue] (o1) -- (i1);
\draw [->-,blue] (i2) -- (o2);
\draw [->-,blue] (i3) -- (o3);
\draw [->-,blue] (o4) -- (i4);
\draw [->-, blue](i1) -- (i2);
\draw [->-,blue] (i4) -- (i1);
\draw [->-,blue] (i2) -- (i3);
\draw [->-,blue] (i4) -- (i3);
\draw [fill] (i1) circle [radius=0.2];
\draw [fill] (i3) circle [radius=0.2];
\draw [fill=white] (i2) circle [radius=0.2];
\draw [fill=white] (i4) circle [radius=0.2];
\draw [dashed] (0,0) circle [radius=2]; 
\node at (-1.5,0) [] {$u'$};
\node at (0,-1.5) [] {$v'$};
\node at (1.5,0) [] {$w'$};
\node at (0,1.5) [] {$x'$};
\node at (0,0) [] {$e'_f$};
\end{tikzpicture}
\]
There are essentially six cases to verify, depending on the position of the path $\gamma$ relative to the parts of the bipartite graphs we drawn above:
\begin{enumerate}
    \item all faces shown above are on the right of $\gamma$;
    \item all faces shown above are on the left of $\gamma$;
    \item $\gamma$ enters via the upper left edge and exits via the upper right edge;
    \item $\gamma$ enters via the upper left edge and exits via the lower right edge;
    \item $\gamma$ enters via the lower left edge and exits via the upper right edge;
    \item $\gamma$ enters via the lower left edge and exits via the lower right edge.
\end{enumerate}
Out of the six cases, (2) is trivially true. To save space, we will verify (1) and (3) below and leave other cases as exercises for the readers.

For case (1), we want to show that $\mu_f^*\left(X^{u'+v'+w'+x'+e'_f}\right)=X^{u+v+w+x+e_f}$. First recall from \eqref{e mutation} that as elements of the lattice $\Lambda$,
\[
u'=u, \quad v'=v+e_f, \quad w'=w, \quad x'=x+e_f, \quad e'_f=-e_f;
\]
therefore we have $u'+v'+w'+x'+e'_f=u+v+w+x+e_f$ as elements of $\Lambda$. Now by \eqref{conjugation} we know that
\begin{align*}
&\mu_f^*\left(X^{u'+v'+w'+x'+e'_f}\right)\\
=&\Ad_{\Psi_q\left(X^{e_f}\right)}\left(X^{u+v+w+x+e_f}\right)\\
=&\prod_{a=1}^{\left|\left\{u+v+w+x+e_f,e_f\right\}\right|}\left(1+q^{(2a-1)\sgn\left\{u+v+w+x+e_f,e_f\right\}}X^{e_f}\right)^{-\sgn\left\{u+v+w+x+e_f,e_f\right\}}X^{u+v+w+x+e_f}\\
=&X^{u+v+w+x+e_f},
\end{align*}
where the last equality is because $\left\{u+v+w+x+e_f,e_f\right\}=0$.

For case (3), we want to show that $\mu_f^*\left(X^{u'+v'+w'+e'_f}\right)=X^{u+v+w}+X^{u+v+w+e_f}$. The reason we have two terms on the right is because in the picture on the left, there are two paths entering via the upper left edge and exiting via the upper right edge; in some sense, the path in $\Gamma'$ ``splits'' into two in $\Gamma$. Again from \eqref{e mutation} we know that $u'+v'+w'+e'_f=u+v+w$, and from \eqref{conjugation} we get
\begin{align*}
    \mu_f^*\left(X^{u'+v'+w'+e'_f}\right)=&\Ad_{\Psi_q\left(X^{e_f}\right)}\left(X^{u+v+w}\right)\\
    =&\prod_{a=1}^{\left|\left\{u+v+w,e_f\right\}\right|}\left(1+q^{(2a-1)\sgn\left\{u+v+w,e_f\right\}}X^{e_f}\right)^{-\sgn\left\{u+v+w,e_f\right\}}X^{u+v+w}\\
    =&\left(1+q^{-1}X^{e_f}\right)X^{u+v+w}\\
    =&X^{u+v+w}+X^{u+v+w+e_f},
\end{align*}
which is exactly what we want to show.
\end{proof}

\subsection{Classical Limit}

By taking the classical limit $q\rightarrow 1$,  on the one hand we recover the cluster Poisson variety $\mathcal{X}_{m,n}$ from its canonical quantization $\mathbb{X}_{m,n}$, and on the other hand we recover the ordinary space of matrices $\mat_{m,n-m}$ on the other hand. Therefore in the classical limit $q\rightarrow 1$, the collection of maps $\chi^*_\Gamma$ gives rise to a map
\[
\chi:\mathcal{X}_{m,n}\rightarrow \mat_{m,n-m}.
\]
In particular, for a given bipartite graph $\Gamma$, we have
\begin{equation}\label{chi}
\chi^*_\Gamma\left(m_{ij}\right)=\sum_{\gamma:i\rightarrow j'} \prod_{f\in \hat{\gamma}}X_f
\end{equation}
in the classical limit $q\rightarrow 1$.

\begin{rmk} The map $\chi$ can be seen as a special case of Postnikov's boundary measurement map defined in \cite{Pos}. In Postnikov's original version there is an extra factor $(-1)^{w(\gamma)}$ in each term; such a factor does not exist in our version because our choice of standard perfect orientation is always acyclic, and hence $\chi^*_\Gamma\left(m_{ij}\right)$ as defined above are always polynomials.
\end{rmk}

\begin{rmk} The map $\chi$ can also be seen as a generalization of Fock and Gocharov's amalgamation map in the case of $\GL_n$ (see \cite{FGamalgamation} for the general definition).
\end{rmk}

Let's now relate the cluster variety $\mathcal{X}_{m,n}$ with configuration space $\conf_n\left(\mathbb{C}^{m-1}\right)$. Let $\overline{w}_0$ be the $m\times m$ matrix
\[
\overline{w}_0:=\begin{pmatrix} 0 & \cdots & 0 & 0 & 1 \\
0 & \cdots & 0 & -1 & 0 \\
0 &  \cdots & 1 & 0 & 0 \\
\vdots & \iddots & \vdots & \vdots & \vdots\\
(-1)^{m-1} & \cdots & 0 & 0 & 0 \end{pmatrix}.
\]
(We denote it by $\overline{w}_0$ because it is a lift of the longest Weyl group element $w_0$ associated to the Lie group $\GL_m$.) There is an embedding $\mat_{m,n-m}\rightarrow \mat_{m,n}$ by juxtaposition 
\[
M\mapsto \left(\overline{w}_0, M\right).
\]

By composing $\chi$ with this embedding we get a map from $\mathcal{X}_{m,n}\rightarrow \mat_{m,n}$. We adopt the convention that we label the rows of the resulting matrix with $m,\dots, 1$ from top to bottom. 
\begin{equation}\label{matrix chi}
\begin{tikzpicture}[baseline=5ex]
\draw (0,0) rectangle (5,2);
\draw (2,0) -- (2,2);
\node at (0,2) [below left] {$m$};
\node at (0,0) [above left] {$1$};
\node at (0,2) [above right] {$1$};
\node at (2,2) [above left] {$m$};
\node at (2,2) [above right] {$m+1$};
\node at (5,2) [above left] {$n$};
\node at (0,1) [left] {$\vdots$};
\node at (1,2) [above] {$\cdots$};
\node at (4,2) [above] {$\cdots$};
\node at (1,1) [] {$\overline{w}_0$};
\node at (3.5,1) [] {$m_{ij}$};
\end{tikzpicture}
\end{equation}

One may wonder why we use the matrix $\overline{w}_0$ instead of something simpler, like an identity matrix. The reason is because we want a certain positivity property, which was proved by Postnikov in \cite{Pos}; for completeness we include a proof here as well.

\begin{prop}\label{positive} Let $I$ be an $m$-element set. Let $I^s:= \{1,\dots, m\}\setminus I$ and let $I^t:=I\setminus \{1,\dots, m\}$. Then the minor 
\[
\Delta_I\left(\overline{w}_0,\left(m_{ij}\right)_{m\geq i\geq 1}^{m<j\leq n}\right)=\sum_{\{\gamma\}:I^s\rightarrow I^t}\prod_{\gamma\in \{\gamma\}}\prod_{f\in \left\{\hat{\gamma}\right\}}X_f,
\]
where $\{\gamma\}:I^s\rightarrow I^t$ means a collection of pairwise non-intersecting paths going from the boundary marked points in $I^s$ to the boundary marked points in $I^t$. In particular, it follows that all minors of the matrix $\left(\overline{w}_0,\left(m_{ij}\right)_{m\geq i\geq 1}^{m<j\leq n}\right)$ are polynomials in $X_f$ with positive integer coefficients.
\end{prop}
\begin{proof} First we know from simple linear algebra that $\Delta_I\left(\overline{w}_0,\left(m_{ij}\right)_{m\geq i\geq 1}^{m<j\leq n}\right)=\Delta_{I^s,I^t}\left(\left(m_{ij}\right)_{m\geq i\geq 1}^{m<j\leq n}\right)$ where $\Delta_{I^s,I^t}$ denotes the minor whose rows are picked by elements of $I^s$ and whose columns are picked by elements of $I^t$. But then it follows from \eqref{chi} that $\Delta_{I^s,I^t}\left(\left(m_{ij}\right)_{m\geq i\geq 1}^{m<j\leq n}\right)$ should be a sum of $\prod_\gamma\prod_{f\in \hat{\gamma}}X_f$ over all collection of paths from $I^s$ to $I^t$. 

So the remaining thing to show is that if there are two paths among the collection of paths going from $I^s$ to $I^t$ intersect, then there is another term in $\Delta_{I^s,I^t}\left(\left(m_{ij}\right)_{m\geq i\geq 1}^{m<j\leq n}\right)$ that cancels it out. Let $\gamma$ and $\eta$ be two intersecting paths in one such collection; then they must share some edges in $\Gamma$ due to the fact that each white vertex only has one in-coming edge and each black vertex only has one out-going edge. 
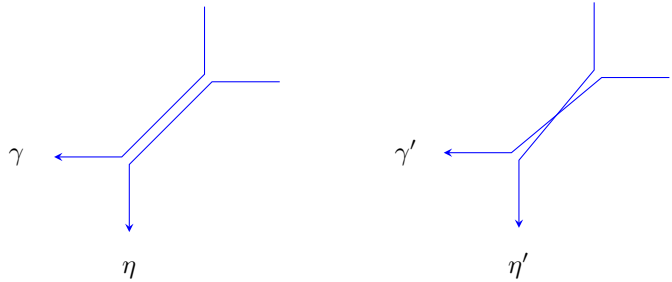
\[
\begin{tikzpicture}
\draw [ ->, blue] (3,2) -- (2.1,2) -- (1,0.9) -- (1,0);
\draw [ ->, blue] (2,3) -- (2,2.1) -- (0.9,1) --(0,1);
\node at (-0.5,1) [] {$\gamma$};
\node at (1,-0.5) [] {$\eta$};
\end{tikzpicture}
\quad \quad \quad \quad
\begin{tikzpicture}
\draw [ ->, blue] (3,2) -- (2.1,2) -- (0.9,1) -- (0,1);
\draw [->, blue] (2,3) -- (2,2.1) -- (1,0.9) --(1,0);
\node at (-0.5,1) [] {$\gamma'$};
\node at (1,-0.5) [] {$\eta'$};
\end{tikzpicture}
\]
But then we can swtich the tails of $\gamma$ and $\eta$ to get another set of paths going from $I^s$ to $I^t$; note that the sets of faces dominated by these two sets of paths are the same with the same multiplicities, but the tail switching produces an extra minus sign, so the two terms cancel out. This leaves us with the only non-vanishing terms in $\Delta_{I^s,I^t}\left(\left(m_{ij}\right)_{m\geq i\geq 1}^{m<j\leq n}\right)$ being from collection of paths that are pairwise non-intersecting.
\end{proof}

Assuming that the image of the map $\mathcal{X}_{m,n}\rightarrow \mat_{m,n}$ intersects $\mat_{m,n}^\times$ non-trivially (recall that $\mat_{m,n}^\times$ consists of matrices whose consecutive $m$ column vectors are in general position), which we will show in Corollary \ref{nontrivial intersection}, we can further compose $\mathcal{X}_{m,n}\dashrightarrow \mat_{m,n}^\times$ with the quotient map $\mat_{m,n}^\times\rightarrow \conf_n^\times\left(\mathbb{C}^m\right)$ to get a rational map, which we abuse notation and denote by $\chi$ as well:
\[
\chi:\mathcal{X}_{m,n}\dashrightarrow \conf_n^\times\left(\mathbb{C}^m\right).
\]

\section{The Hyperplane Map and Cluster Donaldson-Thomas Transformation of Grassmannian}\label{section5}

Let's begin by defining the map $\Xi$ that will play the center role in this section.

\begin{defn} The \emph{hyperplane map} is the map 
\begin{align*} \Xi:\conf^\times_n\left(\mathbb{P}^{m-1}\right)&\rightarrow \conf_n^\times\left(\mathbb{P}^{m-1}\right)\\
\left[l_1,\dots, l_n\right]&\mapsto \left[h_{[2-m,n]}, h_{[3-m,1]},  \dots, h_{[1-m,n-1]}\right],
\end{align*}
where $h_{[i,i+m-2]}$ is defined to be the line dual to the hyperplane spanned by the $m-1$ general positioned lines $l_i,l_{i+1},\dots, l_{i+m-2}$. 
\end{defn}

\begin{rmk} Note that the image of the hyperplane map a priori lives in $\conf_n^\times\left(\left(\mathbb{P}^{m-1}\right)^*\right)$; however, by Proposition \ref{canonical isomorphism} we know that $\conf_n^\times\left(\left(\mathbb{P}^{m-1}\right)^*\right)$ can be canonically identified with $\conf_n^\times\left(\mathbb{P}^{m-1}\right)$ because $\mathbb{P}^{m-1}$ and its dual projective space $\left(\mathbb{P}^{m-1}\right)^*$ have the same dimension.
\end{rmk}

Based on the definition of $\Xi$ it is not hard to prove the following statement.

\begin{prop} $\Xi$ is an isomorphism on $\conf_n^\times\left(\mathbb{P}^{m-1}\right)$.
\end{prop}
\begin{proof} Consider another isomorphism called \emph{reflection}
\begin{align*} \iota:\conf_n^\times\left(\mathbb{P}^{m-1}\right)&\rightarrow \conf_n^\times\left(\mathbb{P}^{m-1}\right)\\
\left[l_1,\dots, l_n\right]&\mapsto \left[l_m,\dots, l_1, l_n,\dots, l_{m+1}\right].
\end{align*}
It is not hard to see that $\iota^2=\mathrm{Id}$. On the other hand, the composition $\tw:=\iota\circ \Xi$, which we call the \emph{twist map}, maps as follows
\[
\tw:\left[l_1,\dots, l_n\right] \mapsto \left[h_{[1,m-1]}, h_{[n,m-2]}, \dots, h_{[2,m]}\right].
\]
Moreover, note that if we take the twist map twice, the first line is a line that is dual to the span of $h_{[1,m-1]}, h_{[n,m-2]}, \dots, h_{[3-m,1]}$, which is obviously $l_1$. By applying the same computation to all other lines it is not hard to see that $\tw^2=\mathrm{Id}$. Therefore $\Xi=\iota\circ \tw$ must also be an isomorphism.
\end{proof}

\begin{rmk} The twist map $\tw$ constructed in the above proof is closely related to the twist map defined by Marsh and Scott in \cite{MStw}.
\end{rmk}

\subsection{A Factorization of the Hyperplane Map \texorpdfstring{$\Xi$}{}}

In this subsection we would like to draw a connection between the hyperplane map $\Xi$ and maps relating to cluster varieties such that $\psi$ and $\chi$. 

Let $\left[l_1,\dots, l_n\right]$ be a configuration in which any $m$ lines are in general position (note that this is strictly stronger that $\left[l_1,\dots, l_n\right]\in \conf_n^\times\left(\mathbb{C}^m\right)$). Fix a representative $\left(l_1,\dots, l_n\right)$ consisting of $n$ lines in $V:= \mathbb{C}^m$, and let $v_i$ be a non-zero vector in the line $l_i$. Then $\left[v_1,\dots, v_n\right]$ is a preimage of $\left[l_1,\dots, l_n\right]$ under the projection map $\conf_n^\times\left(\mathbb{C}^m\right)\rightarrow \conf_n^\times\left(\mathbb{P}^{m-1}\right)$. We denote the matrix with column vectors $v_1,\dots, v_n$ as $M$. Note that by taking the span of the row vectors, $M$ also defines an element in $\dGr_{m,n}^\times:=\dGr_m^\times\left(\mathbb{C}^n\right)$ with $\Delta_I(M)\neq 0$ for all $m$-element subset $I\subset\{1,\dots, n\}$.

Fix a bipartite graph $\Gamma$ of full rank $m$ on a disk with $n$ boundary marked points. We assume without loss of generality that all black vertices of $\Gamma$ are trivalent (see \eqref{delete} and \eqref{splitting}), and assume that every external edge of $\Gamma$ is connected to a unique white vertex (this is obvious because if both $i$ and $i+1$ are connected to the same white vertex, then there is a zig-zag strand going from $i$ to $i+1$ which then implies $m=1$, contradicting the assumption that $m>1$); label the white vertex connecting to the boundary marked point $i$ by $w_i$.

Similar to faces, we define a white vertex $w$ is \emph{dominated} by a zig-zag strand $\zeta_i$ if $w$ lies on the left of $\zeta_i$ with respect to the orientation of $\zeta_i$. Then we define the \emph{dominating set} $I(w)$ to be the indices of all zig-zag strands $\zeta_i$ that dominate $w$. It is not hard to see that $|I(w)|=m-1$ for all white vertex $w$. Now we define a vector $\xi_w\in V^*$ by
\[
\inprod{\xi_w}{v}:=\det\left(v\wedge \overset{\rightarrow}{\bigwedge}_{i\in I(w)}v_i\right),
\]
where $\overset{\rightarrow}{\bigwedge}$ means we are taking the wedge product in the order of ascending indices; in other words, if $I=\left\{i_1<i_2<\dots<i_m\right\}$, then $\overset{\rightarrow}{\bigwedge}_{i\in I}v_i:=v_{i_1}\wedge v_{i_2}\wedge \dots \wedge v_{i_m}$.

Let's investigate what happens locally near a black vertex. Since we have assumed that all black vertices are trivalent, its neighbor must look like the following, assuming that $i<j<k$.
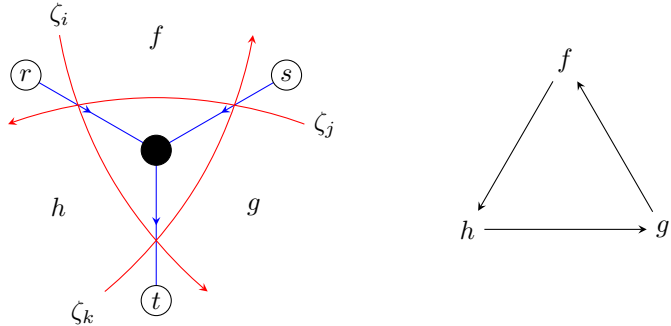
\[
\begin{tikzpicture}[baseline=0ex]
\draw [->-,blue] +(150:2) -- (0,0);
\draw [->-,blue] +(30:2) -- (0,0);
\draw [->-, blue] (0,0) -- (0,-2);
\foreach \i in {0,1,2}
    {
    \draw [fill=white] +(150-120*\i:2) circle [radius=0.2];
    \coordinate (i\i) at +(130-120*\i:2);
    \coordinate (m\i) at +(210-120*\i:0.7);
    \coordinate (o\i) at +(290-120*\i:2);
    \draw [red, ->] plot [smooth, tension=1] coordinates {(i\i)(m\i)(o\i)};
    }
\draw [fill] (0,0) circle [radius=0.2];
\node at +(150:2) [] {$r$};
\node at +(30:2) [] {$s$};
\node at (0,-2) [] {$t$};
\node at +(130:2) [above] {$\zeta_i$};
\node at +(10:2) [right] {$\zeta_j$};
\node at +(250:2) [below left] {$\zeta_k$};
\node at (0,1.5) [] {$f$};
\node at +(-30:1.5) [] {$g$};
\node at +(-150:1.5) [] {$h$};
\end{tikzpicture}
\quad \quad \quad \quad 
\begin{tikzpicture}[baseline=2ex]
\node (f) at (0,1.5) [] {$f$};
\node (g) at +(-30:1.5) [] {$g$};
\node (h) at +(-150:1.5) [] {$h$};
\draw [->] (f) -- (h);
\draw [->] (h) -- (g);
\draw [->] (g) -- (f);
\end{tikzpicture}
\]

\begin{prop} $A_f(M)\xi_t=A_g(M)\xi_r+A_h(M)\xi_s$, where $A_f$, $A_g$, and $A_h$ are understood as their pull-backs under $\psi:\dGr_{m,n}^\times\dashrightarrow \mathcal{A}_{m,n}$.
\end{prop}
\begin{proof} Let $J$ be the $(m-2)$-element set shared by $I(f)$, $I(g)$, and $I(h)$. It is obvious that for any $l\in J$, 
\[
\inprod{-A_f(M)\xi_t+A_g(M)\xi_r+A_h(M)\xi_s}{v_l}=0.
\]
Since $\det \overset{\rightarrow}{\bigwedge}_{l\in I(f)}v_l\neq 0$, it suffices to show
\[
\inprod{-A_f(M)\xi_t+A_g(M)\xi_r+A_h(M)\xi_s}{v_i}=\inprod{-A_f(M)\xi_t+A_g(M)\xi_r+A_h(M)\xi_s}{v_k}=0.
\]
Let $\nu(i)$ and $\nu(k)$ be the number of elements in $J$ that is smaller than $i$ and $k$ respectively. Then it follows that
\[
\inprod{\xi_t}{v_i}=(-1)^{\nu(i)}A_g(M),  \quad \quad \inprod{\xi_r}{v_i}=(-1)^{\nu(i)}A_f(M),\quad \quad \inprod{\xi_s}{v_i}=0,
\]
\[
\inprod{\xi_t}{v_k}=(-1)^{\nu(k)+1}A_h(M), \quad \quad \inprod{\xi_s}{v_k}=(-1)^{\nu(k)+1}A_f(M), \quad \quad \inprod{\xi_r}{v_k}=0.
\]
By plugging these into our desired equalities above, we see that both of them vanish.
\end{proof}

By dividing the identity in the statement of last proposition by $A_f(M)$, we see that
\begin{equation}\label{Y}
\xi_t=\frac{A_g(M)}{A_f(M)}\xi_r+\frac{A_h(M)}{A_f(M)}\xi_s=\frac{A_g(M)}{A_f(M)}\xi_r+\frac{A_g(M)}{A_f(M)}\frac{A_h(M)}{A_g(M)}\xi_s.
\end{equation}
Note that the coefficients of $\xi_r$ and $\xi_s$ resemble a product of cluster $\mathcal{X}$-variables lying on the right hand side of the paths $r\rightarrow t$ and $s\rightarrow t$ respectively, which reminds us of the map $\chi$ we defined at the end of last chapter.

Let's try to make this observation more precise. Suppose $j$ is a sink (which is equivalent to saying that $j>m$). Then we can use the above observation to write $\xi_j:=\xi_{w_j}$ in terms of $\xi_w$ for other white vertices $w$ that are more ``upstream'' in the standard perfect orientation. By going against the standard perfect orientation and repeatedly applying \eqref{Y} we eventually get a linear expansion of $\xi_j$ in terms of $\xi_i:=\xi_{w_i}$ for $1\leq i\leq m$; by carefully comparing \eqref{Y} with \eqref{chi} we can see that
\begin{equation}\label{xi expansion}
\left((-1)^{m-1}\xi_1, (-1)^{m-2}\xi_2,\dots, -\xi_{m-1},\xi_m, \xi_{m+1},\xi_{m+2},\dots, \xi_n\right)=\left(\xi_m,\xi_{m-1},\dots, \xi_1\right)\chi_\Gamma\left(\left(X_f\right)_f\right).
\end{equation}
with suitable choices of $X_f$ that are ratios of cluster $\mathcal{A}$-variables described as in \eqref{Y}. Note that the left hand side of the above equation is a $1\times n$ matrix whose entries are vectors and the right hand side is a vector-valued $1\times m$ matrix $\left(\xi_m,\dots, \xi_1\right)$ multiplied by a scalar valued $m\times n$ matrix that is $\chi_\Gamma\left(\left(X_f\right)_f\right)$. Is $\left(X_f\right)_f=p\left(\left(A_f\right)_f\right)$? Well, the answer is almost correct; the only problem occurs at the boundary faces where, instead of having a full (weight 1) arrow, we only have a weight $\frac{1}{2}$ arrow (for example, $\epsilon_{fg}=-\frac{1}{2}$ in the picture below).
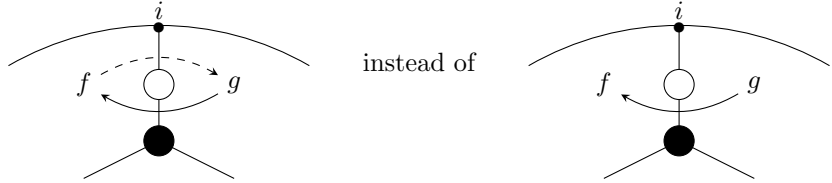
\[
\begin{tikzpicture}[baseline=-0.5ex]
\draw (4,0) arc (60:120:4);
\draw (2,0.5) -- (2,-1);
\node at (2,0.5) [] {$\bullet$};
\node at (2,0.5) [above] {$i$};
\draw (1,-1.5) -- (2,-1) -- (3,-1.5);
\draw [fill=white] (2,-0.25) circle [radius=0.2];
\draw [fill] (2,-1) circle [radius=0.2];
\node (f) at (1,-0.25) [] {$f$};
\node (g) at (3,-0.25) [] {$g$};
\draw [->, dashed] (f) to [bend left] (g);
\draw [->] (g) to [bend left] (f);
\end{tikzpicture}\quad \quad \text{instead of} \quad \quad 
\begin{tikzpicture}[baseline=-0.5ex]
\draw (4,0) arc (60:120:4);
\draw (2,0.5) -- (2,-1);
\node at (2,0.5) [] {$\bullet$};
\node at (2,0.5) [above] {$i$};
\draw (1,-1.5) -- (2,-1) -- (3,-1.5);
\draw [fill=white] (2,-0.25) circle [radius=0.2];
\draw [fill] (2,-1) circle [radius=0.2];
\node (f) at (1,-0.25) [] {$f$};
\node (g) at (3,-0.25) [] {$g$};
\draw [->] (g) to [bend left] (f);
\end{tikzpicture}
\]
So if we use $\left(X_f\right)_f=p\left(\left(A_f\right)_f\right)$, then the values of the cluster $\mathcal{X}$-variables $X_{i\partial}:=X_{f_{i\partial}}$ for boundary faces $f_{i\partial}$ are going to be off by a bit. However, if we only care about the configuration of the spans of $\xi_i$, then we are fine, because of the following proposition.

\begin{prop} Suppose Equation \eqref{xi expansion} holds. Let $\alpha\in \mathbb{C}^\times$.
\begin{enumerate}
\item For $1\leq i<m$, if we let $X'_f=\left\{\begin{array}{ll}
    X_f & \text{if $f\neq f_{i\partial}$}, \\
    \alpha X_{i\partial} & \text{if $f=f_{i\partial}$}
\end{array}\right.$, then the following equality holds:
\begin{align*}
&\left((-1)^{m-1}\xi_1,\dots, (-1)^{m-i}\xi_i, (-1)^{m-i-1}\alpha^{-1}\xi_{i+1}, \dots, \alpha^{-1}\xi_m,\xi_{m+1} \dots, \xi_n\right)\\
=&\left(\alpha^{-1}\xi_m,\dots,  \alpha^{-1}\xi_{i+1},\xi_i, \dots, \alpha\xi_1\right)\chi_\Gamma\left(\left(X'_f\right)_f\right).
\end{align*}
\item For $m<j\leq n$, if we let $X'_f=\left\{\begin{array}{ll}
    X_f & \text{if $f\neq f_{j\partial}$}, \\
    \alpha X_{j\partial} & \text{if $f=f_{j\partial}$}
\end{array}\right.$, then the following equality holds:
\[
\left((-1)^{m-1}\xi_1,\dots, \xi_m, \alpha\xi_{m+1},\dots, \alpha\xi_j, \xi_{j+1},\dots, \xi_n\right)=\left(\xi_m,\dots, \xi_1\right)\chi_\Gamma\left(\left(X'_f\right)_f\right).
\]
\item If we let $X'_f=\left\{\begin{array}{ll}
    X_f & \text{if $f\neq f_{m\partial}$}, \\
    \alpha X_{m\partial} & \text{if $f=f_{m\partial}$}
\end{array}\right.$ only, then the following equality holds:
\[
\left((-1)^{m-1}\xi_1, (-1)^{m-2}\xi_2,\dots, -\xi_{m-1},\xi_m, \xi_{m+1},\xi_{m+2},\dots, \xi_n\right)=\left(\xi_m,\dots, \xi_1\right)\chi_\Gamma\left(\left(X'_f\right)_f\right).
\]
\end{enumerate}
\end{prop}
\begin{proof} Consider the picture below.
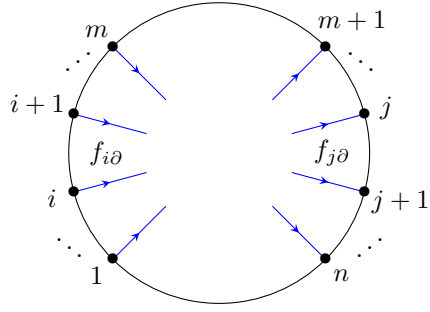
\[
\tikz{
\draw (0,0) circle [radius=2];
\foreach \i in {0,...,3}
    {
    \draw [->-, blue] +(135+\i*30:2) to +(135+\i*30:1);
    \draw [->-, blue] +(45-\i*30:1) to +(45-\i*30:2);
    \node at +(135+\i*30:2) [] {$\bullet$};
    \node at +(45-\i*30:2) [] {$\bullet$};
    }
\node at +(135:2.3) [] {$m$};
\node at +(165:2.5) [] {$i+1$};
\node at +(-165:2.3) [] {$i$};
\node at +(-135:2.3) [] {$1$};
\node at +(145:2.3) [] {$\iddots$};
\node at +(-150:2.3) [] {$\ddots$};
\node at +(45:2.5) [] {$m+1$};
\node at +(15:2.3) [] {$j$};
\node at +(-15:2.5) [] {$j+1$};
\node at +(-45:2.3) [] {$n$};
\node at +(35:2.3) [] {$\ddots$};
\node at +(-30:2.3) [] {$\iddots$};
\node at (-1.5,0) [] {$f_{i\partial}$};
\node at (1.5,0) [] {$f_{j\partial}$};
}
\]

For (1), after scaling $X_{i\partial}$ by $\alpha$, any path $\gamma$ in $\Gamma$ starting from a source $k$ with $i<k\leq m$ will be affected, and the effect is precisely multiplying $\prod_{f\in \hat{\gamma}}X_f$ by $\alpha$. Therefore if $\chi_\Gamma\left(\left(X_f\right)_f\right)=\left(\overline{w}_0, \left(m_{k,l}\right)_{m\geq k\geq 1}^{m<l\leq n}\right)$, then 
\[
\chi_\Gamma\left(\left(X'_f\right)_f\right)=\begin{pmatrix}0 & \cdots  & 1 & \alpha m_{m,m+1} & \cdots & \alpha m_{m,n} \\
\vdots & \iddots & \vdots & \vdots & \ddots & \vdots\\
0 & \cdots & 0 & \alpha m_{i+1, m+1} & \cdots & \alpha m_{i+1, n} \\
0 & \cdots & 0 & m_{i,m+1} & \cdots & m_{i,n}\\
\vdots & \iddots & \vdots & \vdots & \ddots & \vdots \\
(-1)^{m-1} & \cdots & 0 & m_{1,m+1} & \cdots & m_{1,n} \end{pmatrix}.
\]
This implies that if we rescale the vectors $\xi_k$ by $\alpha^{-1}$ for $i<k\leq m$, the identity \eqref{xi expansion} will still hold, which is precisely the statement of (1).

Similarly for (2), after scaling $X_{j\partial}$ by $\alpha$, any path $\gamma$ in $\Gamma$ ending with a sink $l$ with $m<l\leq j$ will be affected, and the effect is precisely multiplying $\prod_{f\in \hat{\gamma}}X_f$ by $\alpha$. Therefore we have
\[
\chi_\Gamma\left(\left(X'_f\right)_f\right)=\begin{pmatrix} 0 & \cdots & 1 & \alpha m_{m,m+1} & \cdots & \alpha m_{m,j} & m_{m,j+1} & \cdots & m_{m,n}\\
\vdots & \iddots & \vdots & \vdots & \ddots & \vdots & \vdots & \ddots & \vdots\\
(-1)^{m-1} & \cdots & 0 & \alpha m_{1,m+1} & \cdots & \alpha m_{1,j} & m_{1,j+1} & \cdots & m_{1,n} 
\end{pmatrix}.
\]
Thus if we rescale the vectors $\xi_l$ by $\alpha$ for $m<l\leq j$, the identity \eqref{xi expansion} will still hold, which is precisely the statement of (2).

For (3), note that the variable $X_{m\partial}$ does not participate in $\chi_\Gamma$ because it is on the right hand side of no paths. Therefore changing $X_{m\partial}$ does not affect the identity \eqref{xi expansion} at all what so ever.
\end{proof}

Note that the span of $\xi_i$ is precisely $h_{[i-m+1,i-1]}$. Therefore the observation above plus our last proposition implies the following corollary.

\begin{cor}\label{factorization of xi} Let $l_1,\dots, l_n$ be $n$ lines in $\mathbb{C}^m$ such that any $m$ of them are in general position. Let $v_i$ be a non-zero vector in $l_i$ for each $i$, and let $M$ be the matrix with column vectors $v_i$. Let $h_{[i,i+m-2]}$ be the line dual to the hyperplane spanned by the $m-1$ lines $l_i, l_{i+1}, \dots, l_{i+m-2}$ (indices modulo $n$). Then
\[
\left[h_{[2-m,n]}, h_{[3-m,1]}, \dots, h_{[1-m,n-1]}\right]=\left[\chi_\Gamma \circ p_\Gamma \circ \psi_\Gamma[M]\right]
\]
where $[M]$ represents the point it defines in the decorated Grassmannian $\dGr_{m,n}^\times$ via $\SL_m\left\backslash \mat_{m,n}^\times \right.\cong \dGr_{m,n}^\times$, and $\left[\chi_\Gamma \circ p_\Gamma \circ \psi_\Gamma[M]\right]$ is again an $m\times n$ matrix in $\mat_{m,n}^\times$ with non-zero column vectors and hence it defines a configuration in $\conf_n^\times \left(\mathbb{P}^{m-1}\right)$ via $\GL_m\left\backslash \mat_{m,n}^\times\right.\cong \conf_n^\times\left(\mathbb{C}^m\right)\rightarrow \conf_n^\times\left(\mathbb{P}^{m-1}\right)$.
\end{cor}

From this corollary we also see that the image of $\chi_\Gamma$, and hence the image of the map $\chi$ defined at the end of last chapter, not only contains some points in $\conf_n^\times\left(\mathbb{C}^m\right)$ but also a dense open subset, proving the following corollary.

\begin{cor}\label{nontrivial intersection} The map $\chi:\mathcal{X}_{m,n}\dashrightarrow \conf_n^\times\left(\mathbb{C}^m\right)$ defined at the end of last chapter is a dominant rational map.
\end{cor}

Let's draw more connection between the cluster variety $\mathcal{X}_{m,n}^\uf$ and the configuration space $\conf_n^\times\left(\mathbb{P}^{m-1}\right)$.

On the one hand, since changing the cluster $\mathcal{X}$-variables corresponding to the boundary faces (i.e., frozen vertices) does not affect the configuration of lines defined by the image of $\chi$, it follows that the rational map $\chi$ can be passed to a rational map
\[
\chi:\mathcal{X}_{m,n}^\uf\dashrightarrow \conf_n^\times \left(\mathbb{P}^{m-1}\right)
\]
which fits into the commutative diagram
\[
\vcenter{\vbox{\xymatrix{\mathcal{X}_{m,n}\ar@{-->}[r]^(0.3){\chi} \ar[d]_q & \conf_n^\times\left(\mathbb{C}^m\right) \ar[d] \\
\mathcal{X}_{m,n}^\uf \ar@{-->}[r]_(0.4){\chi} & \conf_n^\times\left(\mathbb{P}^{m-1}\right)}}}.
\]

On the other hand, we claim that there is another map $\psi:\conf_n^\times\left(\mathbb{P}^{m-1}\right)\dashrightarrow \mathcal{X}_{m,n}^\uf$ that makes the following diagram commutes:
\begin{equation}\label{psicommdiag}
\vcenter{\vbox{\xymatrix{\dGr_{m,n}^\times \ar@{-->}[r]^\psi \ar@<0.5ex>[d]^\pi & \mathcal{A}_{m,n} \ar[d]^{p} \\
\conf_n^\times\left(\mathbb{P}^{m-1}\right) \ar@<0.5ex>[u]^s \ar@{-->}[r]_\psi & \mathcal{X}_{m,n}^\uf}}},
\end{equation}
where $s$ is an arbitrary section against the projection map $\pi$. To prove this claim, it suffices to show that when one computes the ratio between cluster $\mathcal{A}$-coordinates $A_f(M)$ when taking the $p$ map, the result is not affected by scalings of column vectors $v_i$ of $M$. Note that the cluster $\mathcal{A}$-coordinates are minors of $M$, so it suffices to prove the following claim.

\begin{prop}\label{welldefinedpsi} For any fixed non-boundary face $f$ of a bipartite graph $\Gamma$ and any $i\in \{1,\dots, n\}$, the number of faces $g$ dominated by $\zeta_i$ with $\epsilon_{fg}=1$ is equal to the number of faces $h$ dominated by $\zeta_i$ with $\epsilon_{fh}=-1$.
\end{prop}
\begin{proof} This is obvious if $\zeta_i$ dominates all faces surrounding $f$ since $\sum_g \epsilon_{fg}=0$ for any face $f$. If $\zeta_i$ only dominates some of the faces surrounding $f$, then it can be reduced to one of the following two cases using type II 2-by-2 moves, which also implies the result.
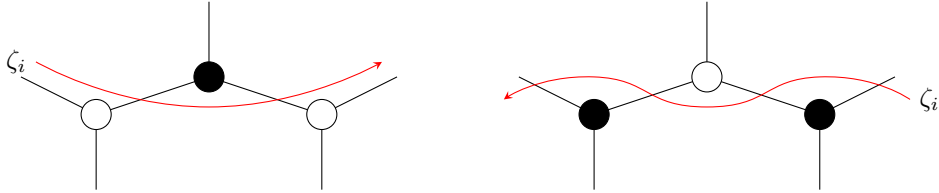
\[
\tikz{
\draw (-1.5,0) -- (-1.5,1) -- (0,1.5) -- (1.5,1) -- (1.5,0);
\draw (0,1.5) -- (0,2.5);
\draw (-1.5,1) -- (-2.5,1.5);
\draw (1.5,1) -- (2.5,1.5);
\draw [fill=white] (1.5,1) circle [radius=0.2];
\draw [fill=white] (-1.5,1) circle [radius=0.2];
\draw [fill] (0,1.5) circle [radius=0.2];
\draw [red, ->]  plot [smooth, tension=1] coordinates {(-2.3,1.7)(0,1.1)(2.3,1.7)};
\node at (-2.3,1.7) [left] {$\zeta_i$};
}\quad \quad \quad \quad
\tikz{
\draw (-1.5,0) -- (-1.5,1) -- (0,1.5) -- (1.5,1) -- (1.5,0);
\draw (0,1.5) -- (0,2.5);
\draw (-1.5,1) -- (-2.5,1.5);
\draw (1.5,1) -- (2.5,1.5);
\draw [fill] (1.5,1) circle [radius=0.2];
\draw [fill] (-1.5,1) circle [radius=0.2];
\draw [fill=white] (0,1.5) circle [radius=0.2];
\draw [red, ->]  plot [smooth, tension=1] coordinates {(2.7,1.2)(1.5,1.5)(0,1.1)(-1.5,1.5)(-2.7,1.2)};
\node at (2.7,1.2) [right] {$\zeta_i$};
}\qedhere
\]
\end{proof}

\begin{cor}\label{dominant of psi} The rational map $\psi$ is dominant.
\end{cor}
\begin{proof} From the proof of Corollary \ref{factorization of xi} we see that the image of $\chi_\Gamma\circ\psi_\Gamma$ contains an dense open subset of $\conf_n^\times\left(\mathbb{P}^{m-1}\right)$; therefore the image of $\psi_\Gamma$ in the seed torus $\mathcal{X}_\Gamma^\uf$ must be of the same dimension as $\conf_n^\times\left(\mathbb{P}^{m-1}\right)$, which is equal to the dimension of $\mathcal{X}_\Gamma^\uf$. Since $\mathcal{X}_\Gamma^\uf$ is an algebraic torus and hence irreducible, the image of $\psi_\Gamma$ must be dense.
\end{proof}

Putting these new rational maps into the picture we get the following commutative diagram:
\[
\vcenter{\vbox{\xymatrix{\dGr_{m,n}^\times \ar@{-->}[r]^{\psi} \ar@<0.5ex>[d]^\pi & \mathcal{A}_{m,n} \ar[d]^{q\circ p} & \\
\conf_n^\times\left(\mathbb{P}^{m-1}\right) \ar@<0.5ex>[u]^s \ar@{-->}[r]^\psi \ar@/_5ex/[rr]_\Xi & \mathcal{X}_{m,n}^\uf \ar@{-->}[r]^(0.4){\chi} & \conf_n^\times\left(\mathbb{P}^{m-1}\right)}}}.
\]

\subsection{Cluster Nature of the Hyperplane Map \texorpdfstring{$\Xi$}{}}\label{cluster nature}

We see in last subsection that the hyperplane map $\Xi$ is birationally equivalent to the composition $\chi\circ \psi$. By reversing the order of composition, we also get a birational map on $\mathcal{X}_{m,n}^\uf$ back to itself. The main goal of this section is to show that $\psi\circ \chi$ is in fact birationally equivalent to a cluster transformation, which by an abuse of notation we also denote as $\Xi$.

To show such a claim, we first need to construct a candidate for the cluster transformation $\Xi$, which requires us to find a pair of isomorphic seeds for $\mathcal{X}_{m,n}^\uf$. 

We begin with the following construction. Given a bipartite graph $\Gamma$ of full rank $m$ on a disk $\mathbb{D}$ whose boundary marked points are evenly distributed, we first draw the diameter bisecting the arcs between the marked points $1$ and $m$; then we reflect $\Gamma$ over this diameter, which results in a new bipartite graph $\Gamma^\circ$ on $\mathbb{D}$. Note that since we assume that the boundary marked points are evenly distributed, each external edges of $\Gamma^\circ$ should connect back to a boundary marked point uniquely. It is not hard to see that the $\Gamma^\circ$ is also a bipartite graph of full rank $m$ on $\mathbb{D}$. We call $\Gamma^\circ$ the \emph{chiral dual} of the bipartite graph $\Gamma$, due to the fact that $\vec{i}_{\Gamma^\circ}=\vec{i}_\Gamma^\circ$. Note that there is a natural one-to-one correspondence between faces of $\Gamma$ and those of $\Gamma^\circ$; we will hence denote the corresponding faces by the same symbol.
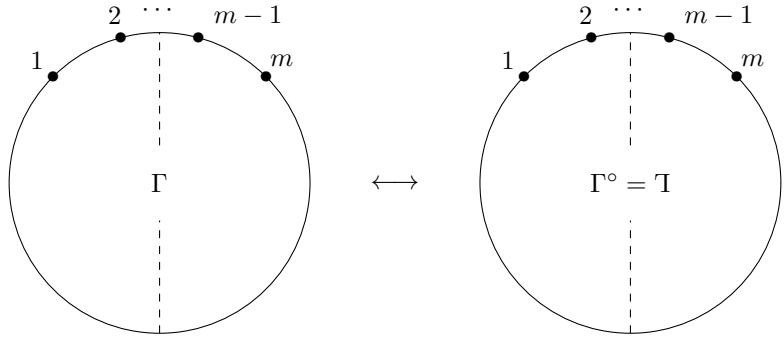
\[
\begin{tikzpicture}[baseline=-0.5ex]
\draw (0,0) circle [radius=2];
\foreach \i in {0,...,3}
    {
    \node at +(135-\i*30:2) [] {$\bullet$};
    }
\draw [dashed] (0,-2) -- (0,-0.5);
\draw [dashed] (0,0.5) -- (0,2);
\node at +(135:2.3) [] {$1$};
\node at +(105:2.3) [] {$2$};
\node at +(75:2.3) [right] {$m-1$};
\node at +(45:2.3) [] {$m$};
\node at (0,2.3) [] {$\cdots$};
\node at (0,0) [] {$\Gamma$};
\end{tikzpicture}\quad \quad \longleftrightarrow \quad \quad 
\begin{tikzpicture}[baseline=-0.5ex]
\draw (0,0) circle [radius=2];
\foreach \i in {0,...,3}
    {
    \node at +(135-\i*30:2) [] {$\bullet$};
    }
\draw [dashed] (0,-2) -- (0,-0.5);
\draw [dashed] (0,0.5) -- (0,2);
\node at +(135:2.3) [] {$1$};
\node at +(105:2.3) [] {$2$};
\node at +(75:2.3) [right] {$m-1$};
\node at +(45:2.3) [] {$m$};
\node at (0,2.3) [] {$\cdots$};
\node at (0,0) [] {$\Gamma^\circ=\text{\reflectbox{$\Gamma$}}$};
\end{tikzpicture}
\]

The seeds we will use to construct the cluster transformation $\Xi$ is the honeycomb bipartite graph $\Gamma:=\Gamma_\hc$ and its chiral dual $\Gamma^\circ=\Gamma_\hc^\circ$. By following the combinatorial procedure we described in Section \ref{section2}, it is not hard to see that the quivers $\vec{i}$ can be drawn as follows.
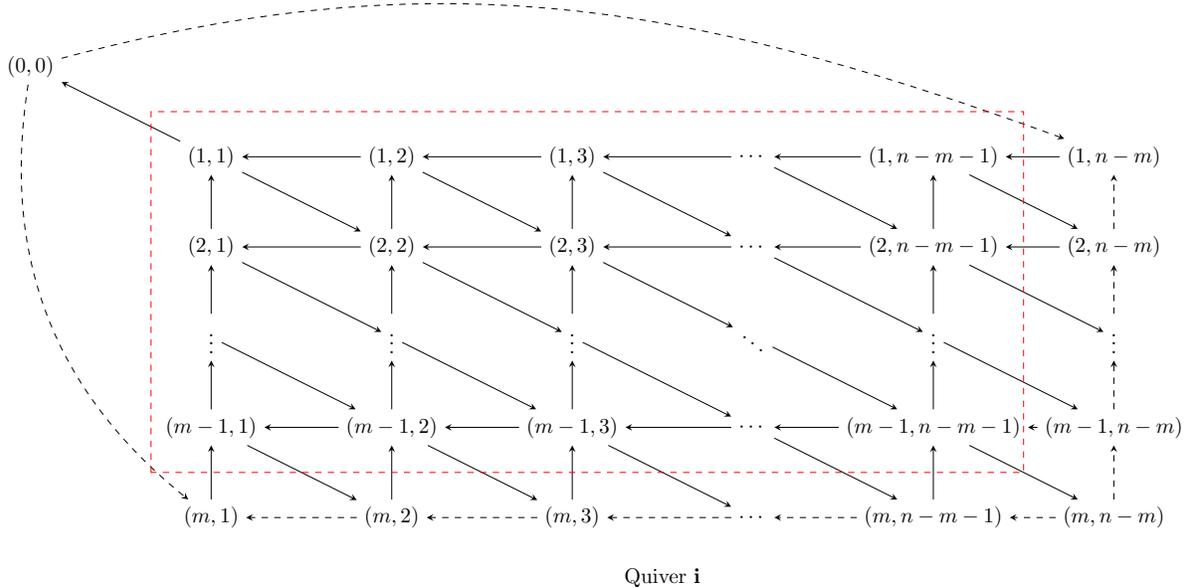
\[
\scalebox{0.8}{\tikz{
\foreach \j in {1,2,3}
    {
    \foreach \i in {1,2}
        {
        \node (\i-\j) at (3*\j,7.5-\i*1.5) [] {$(\i,\j)$};
        }
    \node (3-\j) at (3*\j,3) [] {$\vdots$};
    \node (4-\j) at (3*\j,1.5) [] {$(m-1,\j)$};
    \node (5-\j) at (3*\j,0) [] {$(m,\j)$};
    }
\foreach \i in {1,2,4,5}
    {
    \node (\i-4) at (12,7.5-\i*1.5) [] {$\cdots$};
    }
    \node (3-4) at (12,3) [] {$\ddots$};
\node (1-5) at (15,6) [] {$(1,n-m-1)$};
\node (2-5) at (15,4.5) [] {$(2,n-m-1)$};
\node (3-5) at (15,3) [] {$\vdots$};
\node (4-5) at (15,1.5) [] {$(m-1,n-m-1)$};
\node (5-5) at (15,0) [] {$(m,n-m-1)$};
\node (1-6) at (18,6) [] {$(1,n-m)$};
\node (2-6) at (18,4.5) [] {$(2,n-m)$};
\node (3-6) at (18,3) [] {$\vdots$};
\node (4-6) at (18,1.5) [] {$(m-1,n-m)$};
\node (5-6) at (18,0) [] {$(m,n-m)$};
\node (0-0) at (0,7.5) [] {$(0,0)$};
\foreach \j in {2,...,6}
    {
    \pgfmathtruncatemacro{\l}{\j-1}
    \foreach \i in {1,2,4}
        {
        \draw [->] (\i-\j) -- (\i-\l);
        }
    \draw [dashed, ->] (5-\j) -- (5-\l);
    \foreach \i in {1,...,4}
        {
        \pgfmathtruncatemacro{\k}{\i+1};
        \draw [->] (\i-\l) -- (\k-\j);
        }
    }
\foreach \i in {1,...,4}
    {
    \pgfmathtruncatemacro{\k}{\i+1};
    \foreach \j in {1,2,3,5}
        {
        \draw [->] (\k-\j) -- (\i-\j);
        }
    \draw [dashed, ->] (\k-6) -- (\i-6);
    }
\draw [->] (1-1) -- (0-0);
\draw [->, dashed] (0-0) to [out=15, in=160] (1-6);
\draw [->, dashed] (0-0) to [bend right] (5-1);
\draw [red, dashed] (2,0.75) rectangle (16.5,6.75);
\node at (10.5,-1) [] {Quiver $\vec{i}$};
}}
\]
Note that other than the boundary face at the top left corner, all other faces of $\Gamma$ fit perfectly into an $m\times (n-m)$ grid. Thus we can index the faces of $\Gamma$ using a pair of integers according to the matrix entry convention, leaving the boundary face at the top left corner as $(0,0)$.

On the other hand, we get the following quiver for the chiral dual bipartite graph $\Gamma^\circ$; note that we have labeled the faces of $\Gamma^\circ$ with the same pair of numbers as their counterparts in $\Gamma$.
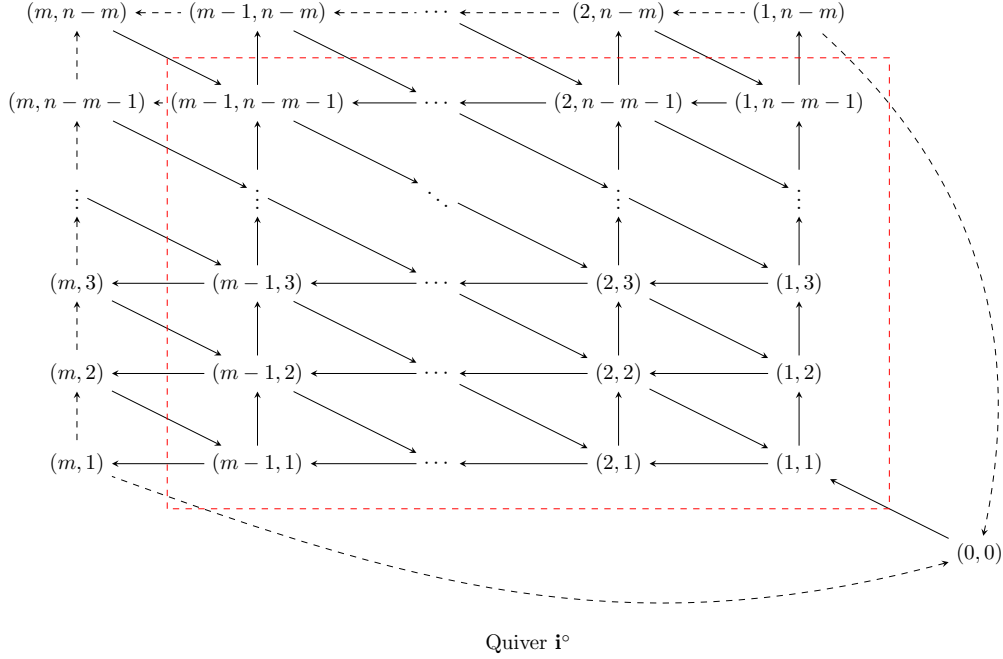
\[
\scalebox{0.8}{\tikz{
\foreach \j in {1,2,3}
    {
    \foreach \i in {1,2}
        {
        \node (\i-\j) at (15-\i*3,1.5*\j) [] {$(\i,\j)$};
        }
    \node (3-\j) at (6,1.5*\j) [] {$\cdots$};
    \node (4-\j) at (3,1.5*\j) [] {$(m-1,\j)$};
    \node (5-\j) at (0,1.5*\j) [] {$(m,\j)$};
    }
\foreach \i in {1,2,4,5}
    {
    \node (\i-4) at (15-\i*3,6) [] {$\vdots$};
    }
    \node (3-4) at (6,6) [] {$\ddots$};
\node (1-5) at (12,7.5) [] {$(1,n-m-1)$};
\node (2-5) at (9,7.5) [] {$(2,n-m-1)$};
\node (3-5) at (6,7.5) [] {$\cdots$};
\node (4-5) at (3,7.5) [] {$(m-1,n-m-1)$};
\node (5-5) at (0,7.5) [] {$(m,n-m-1)$};
\node (1-6) at (12,9) [] {$(1,n-m)$};
\node (2-6) at (9,9) [] {$(2,n-m)$};
\node (3-6) at (6,9) [] {$\cdots$};
\node (4-6) at (3,9) [] {$(m-1,n-m)$};
\node (5-6) at (0,9) [] {$(m,n-m)$};
\node (0-0) at (15,0) [] {$(0,0)$};
\foreach \j in {2,...,6}
    {
    \pgfmathtruncatemacro{\l}{\j-1}
    \foreach \i in {1,2,4}
        {
        \draw [<-] (\i-\j) -- (\i-\l);
        }
    \draw [dashed, <-] (5-\j) -- (5-\l);
    \foreach \i in {1,...,4}
        {
        \pgfmathtruncatemacro{\k}{\i+1};
        \draw [<-] (\i-\l) -- (\k-\j);
        }
    }
\foreach \i in {1,...,4}
    {
    \pgfmathtruncatemacro{\k}{\i+1};
    \foreach \j in {1,2,3,5}
        {
        \draw [<-] (\k-\j) -- (\i-\j);
        }
    \draw [dashed, <-] (\k-6) -- (\i-6);
    }
\draw [<-] (1-1) -- (0-0);
\draw [<-, dashed] (0-0) to [bend right] (1-6);
\draw [<-, dashed] (0-0) to [out=-165,in=-20] (5-1);
\draw [red, dashed] (1.5,0.75) rectangle (13.5,8.25);
\node at (7.5,-1.5) [] {Quiver $\vec{i}^\circ$};
}}
\]

At the first glance the quivers $\vec{i}$ and $\vec{i}$ do not look isomorphic to each other. However, what we need is isomorphic seeds for $\mathcal{X}_{m,n}^\uf$; therefore we only need the unfrozen part of the quivers, which we have identified as the parts inside the dashed rectangles in the above pictures. Then it is obvious that there is a seed isomorphism $\Xi:\left(\vec{i}^\circ\right)^\uf\rightarrow \vec{i}^\uf$ given by
\[
\Xi(m-i,n-m-j)=(i,j).
\]
By Definition \ref{cluster transformation} we get a cluster transformation $\Xi$ on $\mathcal{X}_{m,n}^\uf$ which maps the seed torus $\mathcal{X}_\Gamma^\uf$ to $\mathcal{X}_{\Gamma^\circ}^\uf$.

Now we need to prove the following proposition.

\begin{prop} The cluster transformation $\Xi:\mathcal{X}_{m,n}^\uf\rightarrow \mathcal{X}_{m,n}^\uf$ is birationally equivalent to the composition $\psi\circ \chi$.
\end{prop}
\begin{proof} We make use of the following diagram to prove this proposition.
\[
\xymatrix{
\dGr_{m,n}^\times \ar@{-->}[r]^{\psi_\Gamma} \ar@<0.5ex>[d]^\pi & \mathcal{A}_\Gamma \ar[d]^{p_\Gamma} &  \dGr_{m,n}^\times \ar@{-->}[r]^{\psi_{\Gamma^\circ}} \ar@<0.5ex>[d]^\pi & \mathcal{A}_{\Gamma^\circ} \ar[d]^{ p_{\Gamma^\circ}} \\
\conf_n^\times\left(\mathbb{P}^{m-1}\right) \ar@<0.5ex>[u]^s \ar@{-->}[r]^(0.6){\psi_\Gamma} \ar@/_5ex/[rr]_\Xi & \mathcal{X}_\Gamma^\uf \ar@/_5ex/[rr]_{\Xi?} \ar@{-->}[r]^(0.4){\chi_\Gamma} & \conf_n^\times \left(\mathbb{P}^{m-1}\right) \ar@<0.5ex>[u]^{s'} \ar@{-->}[r]^(0.6){\psi_{\Gamma^\circ}}& \mathcal{X}_{\Gamma^\circ}^\uf
}
\]
Let $\left(X_{i,j}\right)_{0<i<m}^{0<j<n-m}$ be a generic point in $\mathcal{X}_\Gamma^\uf$. Since $\psi_\Gamma$ is dominant (Corollary \ref{dominant of psi}), we can find a configuration $\left[l_1,\dots, l_n\right]$ such that $\psi_\Gamma\left[l_1,\dots, l_n\right]=\left(X_{i,j}\right)_{0<i<m}^{0<j<n-m}$. We will be done if we can show that
\[
\psi_{\Gamma^\circ}\circ \Xi \left[l_1,\dots, l_n\right]=\left(X_{m-i,n-m-j}\right)_{0<i<m}^{0<j<n-m}.
\]

To do this computation, we need to lift to the top level of the commutative diagram. Let $v_i$ be a non-zero vector in each $l_i$ and define dual vectors $\xi_i$ by
\[
\inprod{\xi_i}{v}:=\det \left(v\wedge \overset{\rightarrow}{\bigwedge}_{j\in [i-m+1,i-1]}v_j\right)
\]
Then it follows that
\[
\psi_{\Gamma^\circ}\circ \Xi\left[l_1,\dots, l_n\right]= p_{\Gamma^\circ}\circ \psi_{\Gamma^\circ}\left[\xi_1,\dots, \xi_n\right].
\]
Note that the right hand side of the last equality is asking us to take ratios of minors of the matrix $\left(\xi_1,\dots, \xi_n\right)$ with $\xi_i$ as column vectors. 

But what does it mean to take the determinant of dual vectors? Well, the good news is that we have chosen to be in $\mathbb{C}^m$, so there is a standard inner product $\inprod{\cdot}{\cdot}$ which we can use to identify $\mathbb{C}^m$ and its dual $\left(\mathbb{C}^m\right)^\circ$. After such identification, we can take minors as we used to.

Next we need to figure out the dominating sets for the faces of $\Gamma$ and $\Gamma^\circ$. By some easy combinatorics consideration one can arrive at the following formula (we include a superscript $\circ$ to indicate that these are the dominating sets associated to the bipartite graph $\Gamma^\circ$):
\begin{align*}
&I(i,j)=\{\underbrace{n-m-j+1,\dots, n-m-j+i}_i, \underbrace{n-m+i+1,\dots, n}_{m-i}\},\\
&I^\circ(m-i,n-m-j)=\{\underbrace{1,\dots, i}_i,\underbrace{n-m-j+i+1,\dots, n-j}_{m-i}\}.
\end{align*}

In order to distinguish the minors taken from matrix $\left(v_1,\dots, v_n\right)$ and minors taken from matrix $\left(\xi_1,\dots, \xi_n\right)$, we denote them by $A_{i,j}^v$ and $A_{m-i,n-m-j}^\xi$ respectively. It follows from the dominating set formula above that
\[
A_{m-i,n-m-j}^\xi=\det (\underbrace{\xi_1\wedge\dots \wedge \xi_i}_i\wedge \underbrace{\xi_{n-m-j+i+1}\wedge \dots \wedge \xi_{n-j}}_{m-i}).
\]
This may still seem quite complicated. But then we observe that since $\xi_i$ is a vector that is dual to $v_{i-m+1}\wedge \dots \wedge v_{i-1}$, the minor $A_{m-i,n-m-j}^\xi$ is just a scalar multiple of 
\[
\det (\underbrace{v_{i-m+1}\wedge \dots \wedge v_n}_{m-i}\wedge \underbrace{v_{n-j-m+1}\wedge v_{n-m-j+i}}_i)=A_{i,j}^v,
\]
and the scalar multiple changes linearly with respect to the scaling (different choices) of $v_i$. But then since we know that the ratio (the cluster $\mathcal{X}$-variable) we need does not depend on the choices of $v$, we can effectively replace all $A_{m-i,n-m-j}^\xi$ in our calculation by the minors $A_{i,j}^v$, and the result is only off by at most a sign. Now observe that for a unfrozen vertex $(m-i,n-m-j)$ in $\vec{i}^\circ$ that looks like the following
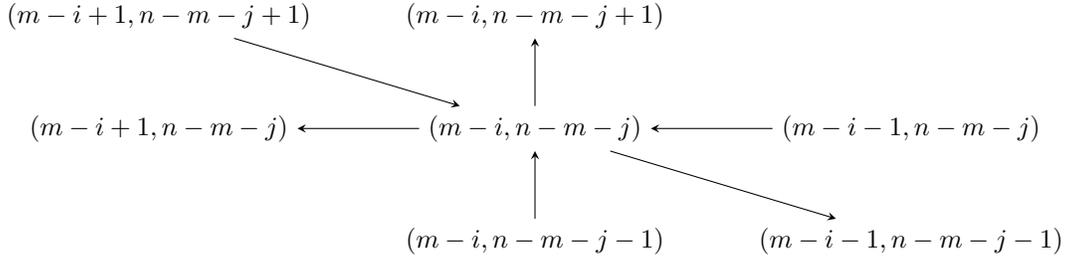
\[
\tikz{
\node (ul) at (0,3) [] {$(m-i+1,n-m-j+1)$};
\node (u) at (5,3) [] {$(m-i,n-m-j+1)$};
\node (l) at (0,1.5) [] {$(m-i+1,n-m-j)$};
\node (c) at (5,1.5) [] {$(m-i,n-m-j)$};
\node (r) at (10,1.5) [] {$(m-i-1,n-m-j)$};
\node (d) at (5,0) [] {$(m-i,n-m-j-1)$};
\node (dr) at (10,0) [] {$(m-i-1,n-m-j-1)$};
\draw [->] (ul) -- (c);
\draw [->] (c) -- (u);
\draw [->] (c) -- (l);
\draw [->] (r) -- (c);
\draw [->] (d) -- (c);
\draw [->] (c) -- (dr);
}
\]
its associated cluster $\mathcal{X}$-variable on the seed torus $\mathcal{X}_{\Gamma^\circ}$ can be computed by
\begin{align*}
X_{m-i,n-m-j}^\xi=&\frac{A_{m-i,n-m-j+1}^\xi A_{m-i+1,n-m-j}^\xi A_{m-i-1,n-m-j-1}^\xi}{A_{m-i+1,n-m-j+1}^\xi A_{m-i-1,n-m-j}^\xi A_{m-i,n-m-j-1}^\xi}\\
=&\pm \frac{A_{i,j-1}^vA_{i-1,j}^vA_{i+1,j+1}^v}{A_{i-1,j-1}^vA_{i+1,j}^vA_{i,j+1}^v};
\end{align*}
but then we notice that the last ratio is just what we get when computing the $p_\Gamma$ map when using the quiver $\vec{i}$, and therefore we can conclude that $\left(\psi_{\Gamma^\circ}\circ \chi_\Gamma\right)^*\left(X^\xi_{m-i,n-m-j}\right)=\pm X_{i,j}^v$. But then Proposition \ref{positive} tells us that the minors of the matrx $\left(\xi_1,\dots, \xi_n\right)$ are always polynomials of the cluster $\mathcal{X}$-variables $X_{i,j}^v$; therefore their ratio must also be a positive rational function in terms of $X_{i,j}^v$. This helps us to fix the sign and conclude that $\left(\psi_{\Gamma^\circ}\circ \chi_\Gamma\right)^*\left(X^\xi_{m-i,n-m-j}\right)=X_{i,j}^v$, which is exactly the same as the pull-back via the cluster transformation $\Xi$.
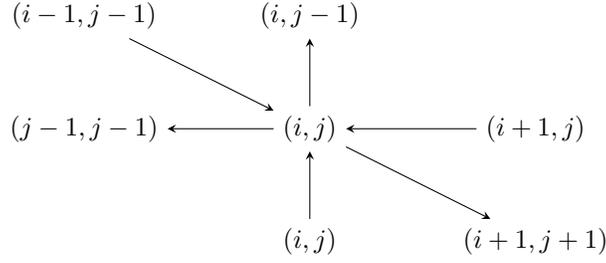
\[
\tikz{
\node (ul) at (0,3) [] {$(i-1,j-1)$};
\node (u) at (3,3) [] {$(i,j-1)$};
\node (l) at (0,1.5) [] {$(j-1,j-1)$};
\node (c) at (3,1.5) [] {$(i,j)$};
\node (r) at (6,1.5) [] {$(i+1,j)$};
\node (d) at (3,0) [] {$(i,j)$};
\node (dr) at (6,0) [] {$(i+1,j+1)$};
\draw [->] (ul) -- (c);
\draw [->] (c) -- (u);
\draw [->] (c) -- (l);
\draw [->] (r) -- (c);
\draw [->] (d) -- (c);
\draw [->] (c) -- (dr);
}
\] 
Other more special cases (e.g. unfrozen vertices $(1,n-m-j)$ and $(m-i,1)$) can be checked in an analogous way as well.
\end{proof}

After knowing that $\psi\circ \chi$ and $\chi\circ \psi$ are both birationally equivalent to isomorphisms, we now give a proof of the following well known result as a byproduct.

\begin{cor} Both $\chi$ and $\psi$ are birational equivalences.
\end{cor}
\begin{proof} Recall that both $\chi$ and $\psi$ are dominant (Corollary \ref{nontrivial intersection} and Corollary \ref{dominant of psi}). Let $U$ be an open dense subset in the intersection between the image of $\psi$ and the domain of $\chi$, let $V:=U\cap \chi^{-1}\circ \Xi\circ \psi^{-1}(U)$, and let $W:=\psi^{-1}(V)$. It is not hard to see that $V$ and $W$ are still open dense subsets of $\conf_n^\times\left(\mathbb{P}^{m-1}\right)$ and $\mathcal{X}_{m,n}^\uf$ respectively. We claim that $\psi$ and $\Xi^{-1}\circ \chi$ are inverses morphisms of each other between $V$ and $W$. It is not hard to see that $\left(\Xi^{-1}\circ \chi\right)\circ \psi=\Xi^{-1}\circ \Xi=\mathrm{Id}_V$. On the other hand, for any $x\in W$, let $x=\psi(y)$ for some $y\in V$; then
\[
\psi\circ \left(\Xi^{-1}\circ\chi\right)(x)=\psi\circ \Xi^{-1}\circ \Xi(y)=\psi(y),
\]
which proves that $\psi$ is a birational equivalence. A similar proof also works for $\chi$.
\end{proof}

\subsection{Proof of \texorpdfstring{$\Xi$ being the Cluster Donaldson-Thomas Transformation}{}}\label{degere property}

In the last section we obtain a cluster analogue of the hyperplane map, which is a cluster transformation $\Xi:\mathcal{X}_{m,n}^\uf\rightarrow \mathcal{X}_{m,n}^\uf$ that is birationally equivalent to the composition $\psi\circ \chi$. In this section we show that the restriction of $\Xi$ on some (and equivalently any) seed torus satisfies the degree property of cluster Donaldson-Thomas transformation $\DT$.

\begin{prop} Let $\Gamma:=\Gamma_\hc$ be the honeycomb bipartite graph and let $\Gamma^\circ$ be its chiral dual. On the seed torus $\mathcal{X}_{\Gamma^\circ}^\uf$, the pull-back via the cluster transformation $\Xi$ satisfies the property
\[
\deg_{X_g}\Xi^*\left(X_f\right)=-\delta_{fg}.
\]
In other words, the cluster transformation $\Xi$ is the cluster Donaldson-Thomas transformation $\DT$ of $\mathcal{X}_{m,n}^\uf$.
\end{prop}
\begin{proof} It is not hard to see that the standard perfect orientation on $\Gamma^\circ$ can be obtained from its counterpart on $\Gamma$ via the same reflection over the diameter bisecting the arcs between $1$ and $m$. Thus from Picture \ref{hc perfect orientation} we get the following schematic picture of the standard perfect orientation on $\Gamma^\circ$.
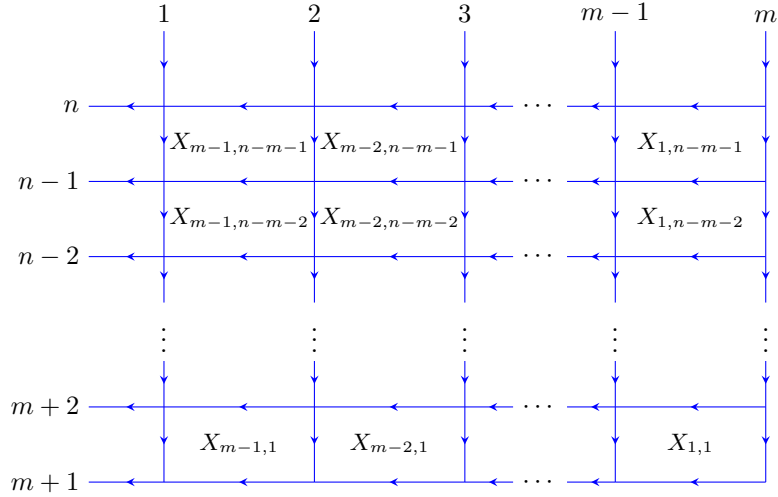
\[
\tikz{
\foreach \i in {1,2,3,4,5}
    {
    \foreach \j in {3,4,5}
        {
        \draw [blue,->-] (2*\i,\j+1) -- (2*\i,\j);
        }
    \node (v\i) at (2*\i,2) [] {$\vdots$};
    \draw [blue,->-] (2*\i,3) -- (v\i);
    \draw [blue,->-] (v\i) -- (2*\i,1);
    \draw [blue,->-] (2*\i,1) -- (2*\i,0);
    }
\foreach \i in {0,1,3,4,5}
    {
    \draw [blue,->-] (4,\i) -- (2,\i);
    \draw [blue,->-] (6,\i) -- (4,\i);
    \node (c\i) at (7,\i) [] {$\cdots$};
    \draw [blue,->-] (c\i) -- (6,\i);
    \draw [blue,->-] (2,\i) -- (1,\i);
    \draw [blue,->-] (8,\i) -- (c\i);
    \draw [blue,->-] (10,\i) -- (8,\i);
    }
\node at (2,6) [above] {$1$};
\node at (4,6) [above] {$2$};
\node at (6,6) [above] {$3$};
\node at (8,6) [above] {$m-1$};
\node at (10,6) [above] {$m$};
\node at (1,5) [left] {$n$};
\node at (1,4) [left] {$n-1$};
\node at (1,3) [left] {$n-2$};
\node at (1,1) [left] {$m+2$};
\node at (1,0) [left] {$m+1$};
\node at (3,4.5) [] {\small{$X_{m-1,n-m-1}$}};
\node at (5,4.5) [] {\small{$X_{m-2,n-m-1}$}};
\node at (3,3.5) [] {\small{$X_{m-1,n-m-2}$}};
\node at (5,3.5) [] {\small{$X_{m-2,n-m-2}$}};
\node at (9,4.5) [] {\small{$X_{1,n-m-1}$}};
\node at (9,3.5) [] {\small{$X_{1,n-m-2}$}};
\node at (3,0.5) [] {\small{$X_{m-1,1}$}};
\node at (5,0.5) [] {\small{$X_{m-2,1}$}};
\node at (9,0.5) [] {\small{$X_{1,1}$}};
}
\]

Note that since $\Xi=\psi_{\Gamma^\circ}\circ \chi_{\Gamma^\circ}= p_{\Gamma^\circ}\circ \psi_{\Gamma^\circ}\circ s\circ \chi_{\Gamma^\circ}$, it follows that in order to compute the degree of $X_g$ in $\Xi^*\left(X_f\right)$, we just need to find out the term in each minor of $\chi_{\Gamma^\circ}$ that gives the highest degree possible, and then take compute the degree according to the $p$ map. 
\[
\xymatrix{
 &  \dGr_{m,n}^\times \ar@{-->}[r]^{\psi_{\Gamma^\circ}} \ar@<0.5ex>[d]^\pi & \mathcal{A}_{\Gamma^\circ} \ar[d]^{ p_{\Gamma^\circ}} \\
 \mathcal{X}_{\Gamma^\circ}^\uf \ar@/_5ex/[rr]_{\Xi} \ar@{-->}[r]^(0.4){\chi_{\Gamma^\circ}} & \conf_n^\times \left(\mathbb{P}^{m-1}\right) \ar@<0.5ex>[u]^{s} \ar@{-->}[r]^(0.6){\psi_{\Gamma^\circ}}& \mathcal{X}_{\Gamma^\circ}^\uf
}
\]
Recall from Proposition \ref{positive} that any minor $\Delta_I$ of $\chi_{\Gamma^\circ}\left(\left(X_g\right)\right)_g$ can be computed by taking 
\[
\Delta_I\left(\chi_{\Gamma^\circ}\left(\left(X_g\right)\right)_g\right)=\sum_{\{\gamma\}:I^s\rightarrow I^t}\prod_{\gamma\in \{\gamma\}}\prod_{f\in \hat{\gamma}}X_f,
\]
where $I^s:=\{1,\dots, m\}\setminus I$ and $I^t:=I\setminus \{1,\dots, m\}$. But then from the picture of the standard perfect orientation above we see that for any collection of $I^s$ and $I^t$ there is a unique collection of paths maximizes the degree of each cluster $\mathcal{X}$-variable, namely the one consisting of $\tikz[baseline=0.5ex]{\draw [blue,->] (0.5,0.5) -- (0.5,0) -- (0,0);}$-shaped paths. Now let's try to compute $\deg_{X_g} \Xi^*\left(X_{m-i,n-m-j}\right)$ for an unfrozen vertex $(m-i,n-m-j)$ that looks like the following.
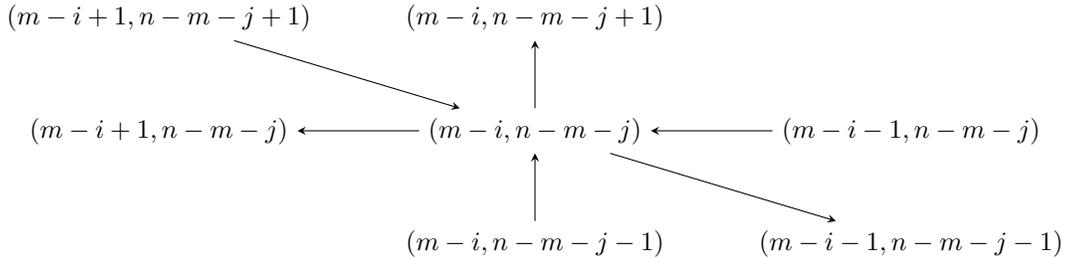
\[
\tikz{
\node (ul) at (0,3) [] {$(m-i+1,n-m-j+1)$};
\node (u) at (5,3) [] {$(m-i,n-m-j+1)$};
\node (l) at (0,1.5) [] {$(m-i+1,n-m-j)$};
\node (c) at (5,1.5) [] {$(m-i,n-m-j)$};
\node (r) at (10,1.5) [] {$(m-i-1,n-m-j)$};
\node (d) at (5,0) [] {$(m-i,n-m-j-1)$};
\node (dr) at (10,0) [] {$(m-i-1,n-m-j-1)$};
\draw [->] (ul) -- (c);
\draw [->] (c) -- (u);
\draw [->] (c) -- (l);
\draw [->] (r) -- (c);
\draw [->] (d) -- (c);
\draw [->] (c) -- (dr);
}
\]

Let's start by breaking the quotient 
\[
\left( p_{\Gamma^\circ}\right)^*\left(X_{m-i,n-m-j}\right)=\frac{A_{m-i,n-m-j+1} A_{m-i+1,n-m-j} A_{m-i-1,n-m-j-1}}{A_{m-i+1,n-m-j+1} A_{m-i-1,n-m-j} A_{m-i,n-m-j-1}}
\]
into a product of three ratios: $\frac{A_{m-i,n-m-j+1}}{A_{m-i-1,n-m-j}}$, $\frac{A_{m-i-1,n-m-j-1}}{A_{m-i+1,n-m-j+1}}$, and $\frac{A_{m-i+1,n-m-j}}{A_{m-i,n-m-j-1}}$. Note that since
\[
I^\circ(m-i,n-m-j)=\{\underbrace{1,\dots, i}_i,\underbrace{n-m-j+i+1,\dots, n-j}_{m-i}\},
\]
it follows that in each of the three ratios, the sets $I^s$ for both the numerator and the denominator only differ by 1 or 2 elements, and so are the sets $I^t$ for the numerator and the denominator. In other words, for highest degree terms in the numerator and the denominator, the collection of paths only differ by 1 or 2 paths.

Let's discuss in more details on the ratio $\frac{A_{m-i,n-m-j+1}}{A_{m-i-1,n-m-j}}$. For the numerator, the set of paths that give the highest degrees is the one consisting of $\tikz[baseline=0.5ex]{\draw [blue,->] (0.5,0.5) -- (0.5,0) -- (0,0);}$-shaped paths that leave $[i+1,m]$ and arrive at $[n-m-j+i+2,n-j+1]$. For the denominator, the set of paths that give the highest degrees is the one consisting of $\tikz[baseline=0.5ex]{\draw [blue,->] (0.5,0.5) -- (0.5,0) -- (0,0);}$-shaped paths that leave $[i+2,m]$ and arrive at $[n-m-j+i+2,n-j]$. Therefore the overall contribution of this ratio to $X_{m-i,n-m-j}$ is a product of $X_g$ with $g$ lying on the right hand side of the $\tikz[baseline=0.5ex]{\draw [blue,->] (0.5,0.5) -- (0.5,0) -- (0,0);}$-shaped path going from $i+1$ to $n-j+1$ in the numerator.

By going through the same analysis, one can find that the overall contribution of the ratio $\frac{A_{m-i-1,n-m-j-1}}{A_{m-i+1,n-m-j+1}}$ is a product of $X_g$ (with multiplicities) with $g$ lying on the right hand side of the $\tikz[baseline=0.5ex]{\draw [blue,->] (0.5,0.5) -- (0.5,0) -- (0,0);}$-shaped path going from $i$ to $n-j+1$ and the one going from $i+1$ to $n-j$ in the denominator, and the overall contribution of the ratio $\frac{A_{m-i+1,n-m-j}}{A_{m-i,n-m-j-1}}$ is a product of $X_g$ with $g$ lying on the right hand side of the $\tikz[baseline=0.5ex]{\draw [blue,->] (0.5,0.5) -- (0.5,0) -- (0,0);}$-shaped path going from $i$ to $n-j$. 
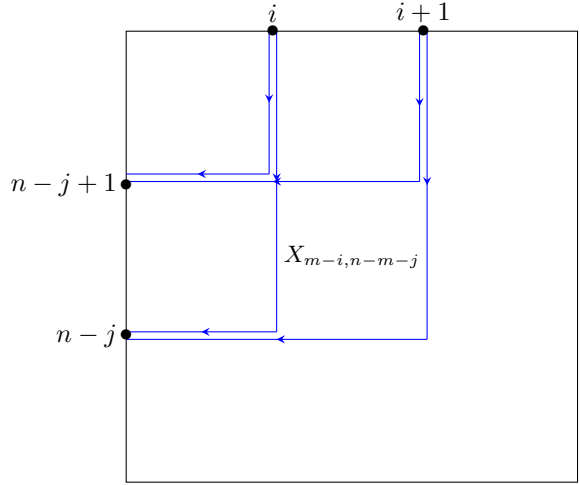
\[
\tikz{
\draw (0,0) rectangle (6,6);
\draw [->-,blue] (1.9,6) -- (1.9,4.1);
\draw [->-,blue] (1.9,4.1) -- (0,4.1);
\draw [->-,blue] (2,6) -- (2,2);
\draw [->-,blue] (2,2) -- (0,2);
\draw [->-,blue] (3.9,6) -- (3.9,4);
\draw [->-,blue] (3.9,4) -- (0,4);
\draw [->-,blue] (4,6) -- (4,1.9);
\draw [->-,blue] (4,1.9) -- (0,1.9);
\node at (1.95,6) [] {$\bullet$};
\node at (1.95,6) [above] {$i$};
\node at (3.95,6) [] {$\bullet$};
\node at (3.95,6) [above] {$i+1$};
\node at (0,1.95) [] {$\bullet$};
\node at (0,1.95) [left] {$n-j$};
\node at (0,3.95) [] {$\bullet$};
\node at (0,3.95) [left] {$n-j+1$};
\node at (3,3) [] {\small{$X_{m-i,n-m-j}$}};
}
\]
Now from the picture above we can easily deduce that $\deg_{X_{m-k,n-m-l}}\Xi^*\left(X_{m-i,n-m-j}\right)=-\delta_{ik}\delta_{jl}$. The statement for other more special cases (e.g. unfrozen vertices $(1,n-m-j)$ and $(m-i,1)$) can be proved in an analogous way as well.
\end{proof}

\subsection{Reflection, Chiral Duality, and the Grassmannian Twist Map}\label{reflection}

We have successfully constructed the $\DT$ version of cluster Donaldson-Thomas transformation on the cluster variety $\mathcal{X}_{m,n}^\uf$, which fits into the commutative diagram below:
\[
\xymatrix{
\conf_n^\times\left(\mathbb{P}^{m-1}\right)\ar@{-->}[r]^(0.6){\psi} \ar@/_5ex/[rr]_\Xi & \mathcal{X}_{m,n}^\uf \ar@{-->}[r]^(0.4){\chi} \ar@/_5ex/[rr]_\DT & \conf_n^\times\left(\mathbb{P}^{m-1}\right) \ar@{-->}[r]^(0.6){\psi} & \mathcal{X}_{m,n}^\uf
}
\]
From this commutative diagram we see that the hyperplane map $\Xi$ can be understood as a geometric realization of the cluster Donaldson-Thomas transformation $\DT$.

Recall that there is another version of cluster Donaldson-Thomas transformation $\dt$, which can be obtained via a conjugation $\dt:=i_{\mathcal{X}}\circ \DT\circ i_\mathcal{X}$. A natural question to ask is whether a geometric realization of $\dt$ also exists; but this is equivalent to asking whether a geometric realization of the chiral duality map $i_\mathcal{X}$ exists, and we will give a positive answer in this section.

Again, let's assume that the boundary marked points on $\mathbb{D}$ are distributed evenly. Recall that we get the chiral dual bipartite graph $\Gamma^\circ$ by flipping $\Gamma$ over the diameter bisecting the arcs between $1$ and $m$ where $m$ is the rank of $\Gamma$. Since $\Gamma^\circ$ also defines a seed torus $\mathcal{X}_{\Gamma^\circ}^\uf$, it follows that the cluster variety $\mathcal{X}_{m,n}^\uf$ is its own chiral dual. In fact, by the same argument we see that $\mathcal{X}_{m,n}$ is its own chiral dual as well without taking the unfrozen part.

Since $\mathcal{X}_{m,n}^\uf$ is chircal dual to itself, so we may expect that there is an auto-isomorphism on $\conf_n^\times\left(\mathbb{P}^{m-1}\right)$ that corresponds to the chiral dual map $i_\mathcal{X}$. It turns out that such an auto-isomorphism does exist, and we have actually seen it already: it is the reflection map
\begin{align*} \iota:\conf_n^\times\left(\mathbb{P}^{m-1}\right)&\rightarrow \conf_n^\times\left(\mathbb{P}^{m-1}\right)\\
\left[l_1,\dots, l_n\right]&\mapsto \left[l_m,\dots, l_1, l_n,\dots, l_{m+1}\right]
\end{align*}
which we have used to show that the hyperplane map $\Xi$ is an isomorphism. 

\begin{prop}\label{11.1} The following diagram commutes.
\[
\xymatrix{ \mathcal{X}_{m,n}^\uf \ar@{-->}[r]^(0.4){\chi} \ar@{<->}[d]_{i_\mathcal{X}} & \conf_n^\times\left(\mathbb{P}^{m-1}\right) \ar@{<->}[d]^\iota \\
\mathcal{X}_{m,n}^\uf \ar@{-->}[r]^(0.4){\chi} &  \conf_n^\times\left(\mathbb{P}^{m-1}\right)
}
\]
\end{prop}
\begin{proof} It suffices to show that the following diagram commutes for any bipartite graph $\Gamma$.
\[
\xymatrix{ \mathcal{X}_\Gamma^\uf \ar@{-->}[r]^(0.4){\chi_\Gamma} \ar@{<->}[d]_{i_\mathcal{X}} & \conf_n^\times\left(\mathbb{P}^{m-1}\right) \ar@{<->}[d]^\iota \\
\mathcal{X}_{\Gamma^\circ}^\uf \ar@{-->}[r]^(0.4){\chi_{\Gamma^\circ}} &  \conf_n^\times\left(\mathbb{P}^{m-1}\right)
}
\]

Recall that in the case of the honeycomb bipartite graph, the standard perfect orientation on $\Gamma_\hc^\circ$ can be obtained from $\Gamma_\hc$ via the same reflection process. Then by using sequences of 2-by-2 moves we can deduce that the standard perfect orientations on any bipartite graph $\Gamma$ and its chiral dual $\Gamma^\circ$ can be related by a reflection through the diameter bisecting 1 and $m$ as well. 

Let $\left(\overline{w}_0,\left(m_{ij}\right)_{m\geq i\geq 1}^{m<j\leq n}\right):=\chi_\Gamma\left(\left(X_f\right)_f\right)$ be the $m\times n$ matrix obtained by following \eqref{chi}, and let $c:=\prod_f X_f$. Note that for a path $\gamma:i\rightarrow j$ in $\Gamma$ compatible with the standard perfect orientation, we also get a corresponding path $\gamma^\circ:m-i+1\rightarrow n-m-j+1$ compatible with the standard perfect orientation in $\Gamma^\circ$.
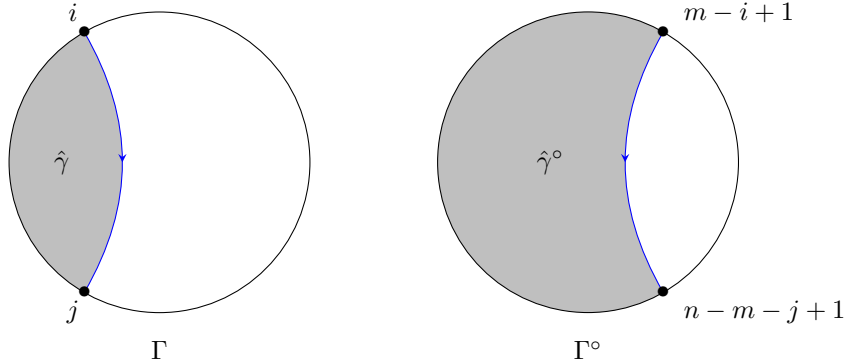
\[
\tikz{
\path [fill=lightgray] (120:2) to [bend left] (-120:2) arc (240:120:2);
\draw [->-, blue] (120:2) to [bend left] (-120:2);
\draw (0,0) circle [radius=2];
\node at (120:2) [] {$\bullet$};
\node at (-120:2) [] {$\bullet$};
\node at (120:2.3) [] {$i$};
\node at (-120:2.3) [] {$j$};
\node at (-1.3,0) [] {$\hat{\gamma}$};
\node at (0,-2.5) [] {$\Gamma$}
}\quad \quad \quad \quad
\tikz{
\path [fill=lightgray] (60:2) to [bend right] (-60:2) arc (300:60:2);
\draw [->-, blue] (60:2) to [bend right] (-60:2);
\draw (0,0) circle [radius=2];
\node at (60:2) [] {$\bullet$};
\node at (-60:2) [] {$\bullet$};
\node at (60:2.3) [right] {$m-i+1$};
\node at (-60:2.3) [right] {$n-m-j+1$};
\node at (-0.5,0) [] {$\hat{\gamma}^\circ$};
\node at (0,-2.5) [] {$\Gamma^\circ$}
}
\]

Let $m_{ij}^\circ$ be the counterpart of $m_{ij}$ for the bipartite graph $\Gamma^\circ$. Then from the above picture we see that 
\[
m_{m-i+1, n-m-j+1}^\circ=\sum_{\substack{\gamma^\circ:m-i+1\\
\rightarrow n-m-j+1}}\prod_{f\in \hat{\gamma}^\circ}X_f^\circ=\sum_{\gamma:i\rightarrow j}\prod_{f\notin \hat{\gamma}} X_f^{-1}=\frac{1}{c}\sum_{\gamma:i\rightarrow j}\prod_{f\in \hat{\gamma}} X_f=\frac{1}{c}m_{ij}.
\]
From this calculation we can deduce that
\begin{align*}
    \chi^\circ_{\Gamma^\circ}\circ i_\mathcal{X}\left(\left(X_f\right)_f\right)=&\begin{tikzpicture}[baseline=7ex]
\draw (0,0) rectangle (5,2);
\draw (2,0) -- (2,2);
\node at (0,2) [below left] {$m$};
\node at (0,0) [above left] {$1$};
\node at (0,2) [above right] {$1$};
\node at (2,2) [above left] {$m$};
\node at (2,2) [above right] {$m+1$};
\node at (5,2) [above left] {$n$};
\node at (0,1) [left] {$\vdots$};
\node at (1,2) [above] {$\cdots$};
\node at (4,2) [above] {$\cdots$};
\node at (1,1) [] {$\overline{w}_0$};
\node at (3.5,1) [] {$m_{ij}^\circ$};
\end{tikzpicture}\\
=&\begin{tikzpicture}[baseline=7ex]
\draw (0,0) rectangle (5,2);
\draw (2,0) -- (2,2);
\node at (0,2) [below left] {$m$};
\node at (0,0) [above left] {$1$};
\node at (0,2) [above right] {$1$};
\node at (2,2) [above left] {$m$};
\node at (2,2) [above right] {$m+1$};
\node at (5,2) [above left] {$n$};
\node at (0,1) [left] {$\vdots$};
\node at (1,2) [above] {$\cdots$};
\node at (4,2) [above] {$\cdots$};
\node at (1,1) [] {$\overline{w}_0$};
\node at (3.5,1) [] {$\frac{1}{c}m_{m-i+1,n-m-j+1}$};
\end{tikzpicture}\\
=&\left(\begin{tikzpicture}[baseline=7ex]
\draw (0,0) rectangle (2,2);
\node at (1,1) [] {$\overline{w}_0$};
\node at (0,0) [above left] {$1$};
\node at (0,2) [below left] {$m$};
\node at (0,2) [above right] {$1$};
\node at (2,2) [above left] {$m$};
\node at (0,1) [left] {$\vdots$};
\node at (1,2) [above] {$\cdots$}; 
\end{tikzpicture}\right)\left(
\begin{tikzpicture}[baseline=7ex]
\draw (0,0) rectangle (5,2);
\draw (2,0) -- (2,2);
\node at (0,2) [below left] {$m$};
\node at (0,0) [above left] {$1$};
\node at (0,2) [above right] {$1$};
\node at (2,2) [above left] {$m$};
\node at (2,2) [above right] {$m+1$};
\node at (5,2) [above left] {$n$};
\node at (0,1) [left] {$\vdots$};
\node at (1,2) [above] {$\cdots$};
\node at (4,2) [above] {$\cdots$};
\node at (1,1) [] {$\mathrm{Id}_m$};
\node at (3.5,1) [] {$\frac{(-1)^{i-1}}{c}m_{i,n-m-j+1}$};
\end{tikzpicture}\right),
\end{align*}
Note that if we ignore the $\GL_m$ matrix $\overline{w}_0$ on the left, each the column vector in remaining rectangular matrix spans the same line as the corresponding one in $\iota\circ \chi_\Gamma\left(\left(X_f\right)_f\right)$. This finishes our proof.
\end{proof}

Since $\dt=i_\mathcal{X}\circ \DT\circ i_\mathcal{X}$ and we have found the geometric realization of each of the maps on the right, we can now define another hyperplane map $H:=\iota\circ \Xi\circ \iota$. By a simple computation one can show that $H$ maps configurations as follows:
\[
H:\left[l_1,\dots, l_n\right]=\left[h_{[4-m,2]}, h_{[5-m,3]},\dots, h_{[3-m,1]}\right].
\]
Similar to the way we view $\Xi$ as a geometric realization of $\DT$, we should also view $H$ as a geometric realization of the $\dt$ version of the cluster Donaldson-Thomas transformation.

Recall that we have also defined a twist map $\tw:=\iota\circ \Xi$. In \cite{GS}, Goncharov and Shen constructed another cluster theoretical involution $D_\mathcal{X}:=i_\mathcal{X}\circ \DT$, and it is not hard to see that the twist map is just a geometric realization of such map $D_\mathcal{X}$. As a summary, we put all the maps we have mentioned so far into the following commutative diagram.
\begin{equation}\label{summary}
\vcenter{\vbox{\xymatrix{&
\dGr_{m,n}^\times \ar@{-->}[r]^\psi \ar@<0.5ex>[d]^\pi & \mathcal{A}_{m,n} \ar[d]^(0.3){q\circ p} & \dGr_{m,n}^\times \ar@{-->}[r]^\psi \ar@<0.5ex>[d]^(0.3){\pi} & \mathcal{A}_{m,n} \ar[d]^{q\circ p} & \\
\cdots \ar@{-->}[r] & \conf_n^\times\left(\mathbb{P}^{m-1}\right) \ar@<0.5ex>[u]^s \ar@{<->}[d]_\iota \ar@{<->}[drr]^(0.6){\tw} \ar@{-->}[r]^(0.6){\psi} \ar@/^5ex/[rr]^(0.3){\Xi} & \mathcal{X}_{m,n}^\uf \ar@{<->}[drr]^(0.6){D_\mathcal{X}} \ar@{-->}[r]^(0.4){\chi} \ar@/^5ex/[rr]^(0.7){\DT} \ar@{<->}[d]_(0.6){i_\mathcal{X}} & \conf_n^\times\left(\mathbb{P}^{m-1}\right)\ar@{-->}[r]^(0.6){\psi} \ar@<0.5ex>[u]^(0.7){s} \ar@{<->}[d]_(0.6){\iota} & \mathcal{X}_{m,n}^\uf \ar@{<->}[d]^{i_\mathcal{X}} \ar@{-->}[r] & \cdots \\
\cdots \ar@{-->}[r] & \conf_n^\times \left(\mathbb{P}^{m-1}\right)  \ar@/_5ex/[rr]_H & \mathcal{X}_{m,n}^\uf  \ar@{-->}[r]_(0.4){\chi} \ar@/_5ex/[rr]_{\dt} & \conf_n^\times\left(\mathbb{P}^{m-1}\right) & \mathcal{X}_{m,n}^\uf \ar@{-->}[r] & \cdots
}}}
\end{equation}


\section{Generalizations} \label{section6}

\subsection{Generalization to Double Bruhat Cells in Semisimple Lie Groups}

Double Bruhat cells were first introduced by Fomin and Zelevinsky in \cite{FZ}; their associated cluster algebra structures were studied by Berenstein, Fomin, and Zelevinsky in \cite{BFZ}; Fock and Goncharov described cluster Poisson structures on double Bruhat cells in \cite{FGamalgamation}. Utilizing these structures, we are able to find the analogue of all the maps in the summarized commutative diagram at the end of last section and construct the cluster Donaldson-Thomas transformations in those cases. We will only briefly describe the relevant maps and structures in this section; the full proof can be found in \cite{Wengdb}.

Let $G$ be a simply connected semisimple Lie group. Fix a pair of opposite Borel subgroups $B_\pm$. Then $H:=B_+\cap B_-$ is a maximal torus in $G$, and $W:=N_GH/H$ is the associated Weyl group. 

There is an open dense subset of elements in a semisimple Lie group $G$ that admits a triangle decomposition of the form $N_-HN_+$ where $N_\pm:=\left[B_\pm, B_\pm\right]$ are the maximal unipotent subgroups corresponding to $B_\pm$. For a group element $x$ that is triangle decomposable, we denote the factors in the factorization as
\[
x=[x]_-[x]_0[x]_+.
\]

The choice of opposite Borel subgroups $B_\pm$ determines a set of simple roots $\{\alpha\}$, which in turn give rise to Coxeter generators $s_\alpha$ of the Weyl group $W$. In particular, the Coxeter generators satisfy a set of relations called the braid relations.

The choice of opposite Borel subgroups $B_\pm$ produce two Bruhat decompositions of the group: 
\[
G=\bigsqcup_{w\in W} B_+wB_+=\bigsqcup_{w\in W} B_-wB_-.
\]
Given a pair of Weyl group elements $(u,v)$, the associated \emph{double Bruhat cell} is defined to be 
\[
G^{u,v}:=B_+uB_+\cap B_-vB_-.
\]
Berenstein, Fomin, and Zelevinsky showed in \cite{BFZ} that the coordinate ring $\mathcal{O}\left(G^{u,v}\right)$ of a double Bruhat cell has the structure of an upper cluster algebra, which can be turned into a map
\[
\psi:G^{u,v}\dashrightarrow \mathcal{A}^{u,v}_G
\]
for certain cluster variety $\mathcal{A}^{u,v}_G$. There is a collection of seed tori in $\mathcal{A}^{u,v}_G$ that are constructed from reduced words of the pair of Weyl group elements $(u,v)$.

Fock and Goncharov introduced an amalgamation procedure in \cite{FGamalgamation} to describe the cluster Poisson structure on the double Bruhat cells, which can be used to define a map 
\[
\chi:\left(\mathcal{X}^{u,v}_G\right)\uf\dashrightarrow H\left\backslash G^{u,v}\right/H.
\]
Note that there is also an obvious projection map
\[
\pi:G\rightarrow H\left\backslash G^{u,v}\right/H.
\]

Let $\left\{E_{\pm \alpha}, H_\alpha\right\}$ be a collection of Chevalley generators of $G$ with respect to the choice of the pair of opposite Borel subgroups $B_\pm$. Exponentiating the generators $E_{\pm \alpha}$ we get group elements $e_{\pm \alpha}$ in $G$. There is an anti-involution $\iota$ on $G$ defined by 
\[
\iota\left(e_{\pm \alpha}\right)=e_{\pm \alpha} \quad \quad \text{and} \quad \quad \iota(h)=h^{-1} \quad \text{for $h\in H$}.
\]
It is not hard to see that this anti-involution descends to an isomorphism
\[
\iota:H\left\backslash G^{u,v}\right/H\rightarrow H\left\backslash G^{u^{-1},v^{-1}}\right/H.
\]

The group elements $e_{\pm \alpha}$ can also be used to produce lifts of the Weyl group elements of $W$ to the group $G$: one first lift the Coxeter generators $s_\alpha$ to $\overline{s}_\alpha=e_\alpha^{-1}e_{-\alpha}e_\alpha^{-1}$, and then use a (and equivalently any) reduced word of a Weyl group element $w$ to define the lift $\overline{w}$ to be the product of the lifts of the letters in the reduced word. In particular, such lift $\overline{w}$ is independent of the choice of reduced word for $w$.

Fomin and Zelevinsky introduced a \emph{twist map} on a double Bruhat cell $G^{u,v}$ in \cite{FZ}, which is an isomorphism between double Bruhat cells defined by
\begin{align*}
    \tw:G^{u,v}&\rightarrow G^{u^{-1},v^{-1}}\\
    x& \mapsto \left(\left[\overline{u}^{-1}x\right]_-^{-1}\overline{u}^{-1}x\overline{v^{-1}}\left[x\overline{v^{-1}}\right]_+^{-1}\right)^{t\circ \iota}.
\end{align*}
Here $t$ denotes the transposition map, which is also an anti-involution and it commutes with the anti-involution $\iota$. At the end, we define $\Xi:=\iota \circ \tw$ and $H:=\tw\circ \iota$.

Now we have all the pieces ready, and we put them into the following commutative diagram, which is analogous to the one we had in the last section.
\[
\xymatrix{&
G^{u,v} \ar@{-->}[r]^\psi \ar@<0.5ex>[d]^\pi & \mathcal{A}^{u,v}_G \ar[d]^(0.3){q\circ p} & G^{u,v} \ar@{-->}[r]^\psi \ar@<0.5ex>[d]^(0.3){\pi} & \mathcal{A}^{u,v}_G \ar[d]^{q\circ p} & \\
\cdots \ar@{-->}[r] & H\left\backslash G^{u,v}\right/H \ar@<0.5ex>[u]^s \ar@{<->}[d]_\iota \ar@{<->}[drr]^(0.6){\tw} \ar@{-->}[r]^(0.6){\psi} \ar@/^5ex/[rr]^(0.3){\Xi} & \left(\mathcal{X}^{u,v}_G\right)^\uf \ar@{<->}[drr]^(0.6){D_\mathcal{X}} \ar@{-->}[r]^(0.4){\chi} \ar@/^5ex/[rr]^(0.7){\DT} \ar@{<->}[d]_(0.6){i_\mathcal{X}} & H\left\backslash G^{u,v}\right/H \ar@{-->}[r]^(0.6){\psi} \ar@<0.5ex>[u]^(0.7){s} \ar@{<->}[d]_(0.6){\iota} & \left(\mathcal{X}^{u,v}_G\right)^\uf \ar@{<->}[d]^{i_\mathcal{X}} \ar@{-->}[r] & \cdots \\
\cdots \ar@{-->}[r] & H\left\backslash G^{u^{-1},v^{-1}}\right/H  \ar@/_5ex/[rr]_H & \left(\mathcal{X}^{u,v}_G\right)^\uf  \ar@{-->}[r]_(0.4){\chi} \ar@/_5ex/[rr]_{\dt} & H\left\backslash G^{u^{-1},v^{-1}}\right/H & \left(\mathcal{X}^{u,v}_G\right)^\uf \ar@{-->}[r] & \cdots}
\]
Note that the cluster varieties $H\left\backslash G^{u,v}\right/H$ and $\left(\mathcal{X}^{u,v}_G\right)\uf$ are actually not chiral dual to themselves, but to their counterparts with subscripts $\left(u^{-1},v^{-1}\right)$.

\subsection{Bipartite Graphs on a Disk and Double Bruhat Cells in \texorpdfstring{$\GL_n$}{}}

The parallel between Grassmannian and the double Bruhat cells in semisimple Lie groups in general may stop at the end of last subsection, but in the special case of $\GL_n$ (which is reductive but very close to being semisimple), our bipartite graph method still applies and the maps $\psi$ and $\chi$ in the last commutative diagram arise in the exact same way as we have introduced them in Section \ref{section4}. Such connection between Grassmannian and double Bruhat cells in $\GL_n$ has been observed by Gekhtman, Shapiro, and Vainshtein in \cite{GSVdirectednetwork}. For the rest of this section we will describe such connection.

Fix $B_+$ to be the group of upper triangular matrices and $B_-$ to be the group of lower triangular matrices. Then the corresponding maximal torus $H$ is the group of diagonal matrices. Let $\left\{\alpha_i\right\}_{i=1}^{n-1}$ be the set of simple roots of $\GL_n$ defined by 
\[
\left(\diag\left(a_1,\dots, a_n\right)\right)^{\alpha_i}=\frac{a_i}{a_{i+1}}.
\]

The Weyl group $W$ associated to $\GL_n$ is the permutation group $S_n$, which acts on the maximal torus $H$ by permuting its entries. The Coxeter generator $s_i$ corresponding to the simple root $\alpha_i$ can be seen as an adjacent transposition between $i$ and $i+1$. A \emph{reduced word} of a pair of Weyl group elements $(u,v)$ is then a sequence formed by elements from the set $\{\pm 1,\dots, \pm (n-1)\}$ such that when all positive entries are deleted and change the sign of the remaining letters we get a reduced word of $u$ and when all negative entries are deleted we get a reduced word of $v$.

Given a reduced word $\vec{i}=\left(i_1,\dots, i_l\right)$ of a pair of Weyl group elements $(u,v)$, we draw a bipartite graph (not full rank) on a disk with $2n$ boundary marked points as follows. First we deform the disk into a rectangle with $n$ boundary marked points on each of the left edge and the right edge. Label the boundary marked points on the left edge by $1, \dots, n$ from top to bottom and label the boundary marked points on the right edge by $1',\dots, n'$ from top to bottom (similar to the labellings we have used in Section \ref{quantum}). Draw $n$ parallel lines connecting $i$ and $i'$ for each $i$; we call the one that connects $i$ and $i'$ the \emph{$i$th horizontal lines}. These horizontal lines decompose the disk into $n+1$ connected components, and we call the component between the $i$th and the $(i+1)$th horizontal lines the \emph{$i$th spacing}; in particular the component above the first horizontal line is the $0$th spacing and the component below the last horizontal line is the $n$th spacing. Now as we go from the first letter to the last letter through the reduced word $\vec{i}$ and from left to right through the rectangular disk, we draw one of two vertical patterns across a spacing for each letter: if $i_k$ is positive, then we draw \begin{tikzpicture}[baseline=2ex] \draw (0,0) -- (0,1); \draw [fill=white] (0,1) circle [radius=0.2]; \draw [fill] (0,0) circle [radius=0.2];\end{tikzpicture} across the $i_k$th spacing, and if $i_k$ is negative, then we draw \begin{tikzpicture}[baseline=2ex] \draw (0,0) -- (0,1); \draw [fill] (0,1) circle [radius=0.2]; \draw [fill=white] (0,0) circle [radius=0.2];\end{tikzpicture} across the $i_k$th spacing. Note that the resulting graph on a disk may not be bipartite, but we can always make it into bipartite by adding a bivalent vertex of the opposite color to an edge connecting two vertices of the same color (recall that bivalent vertices do not affect the associated quiver or dominating sets or perfect orientations). We denote the resulting bipartite graph by $\Gamma_\vec{i}$.

\begin{exmp}\label{12.1} Consider the double Bruhat cell $\GL_n^{w_0,w_0}$ where $w_0$ is the longest Weyl group element and consider the reduced word $\vec{i}=(1,-2,2,1,-1,-2)$. The associated bipartite graph $\Gamma_\vec{i}$ looks like the following.
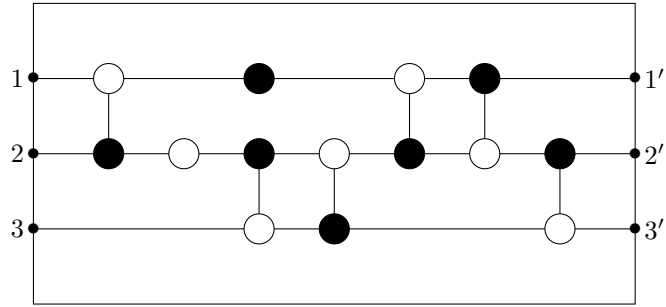
\[
\tikz{
\draw (-1,-1) rectangle (7,3);
    \draw (-1,0) -- (7,0);
    \draw (-1,1) -- (7,1);
    \draw (-1,2) -- (7,2);
    \draw (0,1) -- (0,2);
    \draw (2,1) -- (2,0);
    \draw (3,1) -- (3,0);
    \draw (4,1) -- (4,2);
    \draw (5,1) -- (5,2);
    \draw (6,1) -- (6,0);
    \draw[fill=white] (0,2) circle [radius=0.2];
    \draw[fill=black] (0,1) circle [radius=0.2];
    \draw[fill=white] (2,0) circle [radius=0.2];
    \draw[fill=black] (2,1) circle [radius=0.2];
    \draw[fill=white] (3,1) circle [radius=0.2];
    \draw[fill=black] (3,0) circle [radius=0.2];
    \draw[fill=white] (4,2) circle [radius=0.2];
    \draw[fill=black] (4,1) circle [radius=0.2];
    \draw[fill=white] (5,1) circle [radius=0.2];
    \draw[fill=black] (5,2) circle [radius=0.2];
    \draw[fill=white] (6,0) circle [radius=0.2];
    \draw[fill=black] (6,1) circle [radius=0.2];
    \draw[fill=white] (1,1) circle [radius=0.2];
    \draw[fill=black] (2,2) circle [radius=0.2];
    \foreach \i in {1,2,3}
        {
        \node at (-1,3-\i) [] {$\bullet$};
        \node at (7,3-\i) [] {$\bullet$};
        \node at (-1,3-\i) [left] {$\i$};
        \node at (7,3-\i) [right] {$\i'$};
        }
}
\]

If we draw zig-zag strands on a bipartite graph $\Gamma_\vec{i}$, we see that there are two types: the ones that goes from left to right and the ones that goes from right to left. Moreover we observe the following.

\begin{prop} Let $\vec{i}$ be a reduced word of $(u,v)$. Then the zig-zag strands in $\Gamma_\vec{i}$ going from left to right goes from $i$ to $u(i)$ and the zig-zag strands going from right to left goes from $v(i)$ to $i$ for each $1\leq i<n$.
\end{prop}
\begin{proof} It suffices to note that the zig-zag strands going from left to right tangle at vertical edges of the form \begin{tikzpicture}[baseline=2ex] \draw (0,0) -- (0,1); \draw [fill] (0,1) circle [radius=0.2]; \draw [fill=white] (0,0) circle [radius=0.2];\end{tikzpicture} whereas the zig-zag strands going from right to left tangle at vertical edges of the form \begin{tikzpicture}[baseline=2ex] \draw (0,0) -- (0,1); \draw [fill=white] (0,1) circle [radius=0.2]; \draw [fill] (0,0) circle [radius=0.2];\end{tikzpicture}.
\end{proof}
\end{exmp}

Recall that one important property of bipartite graph that we use is Thurston's theorem \ref{thurston}, which says that two bipartite graphs of full rank $m$ on a disk with $n$ boundary marked points can be transformed into one another via a sequence of 2-by-2 moves. In the case of $\Gamma_\vec{i}$, 2-by-2 moves can occur in three occasions, which we list as follows.
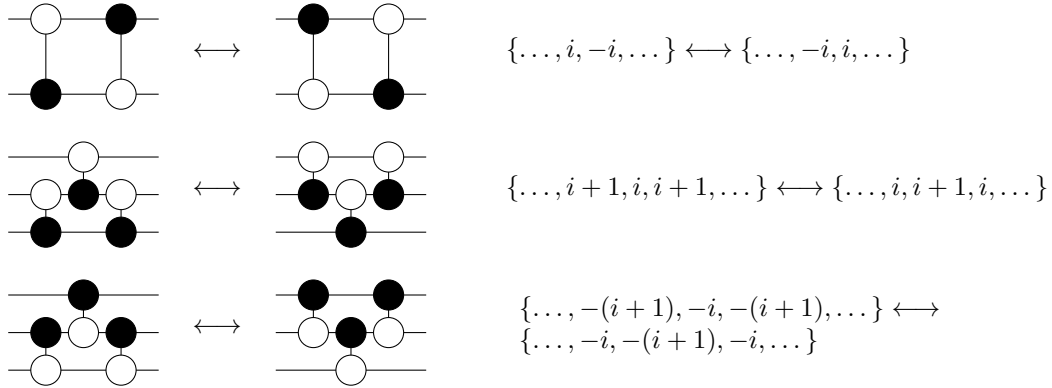
\[
\begin{array}{ll}
        \begin{tikzpicture}[baseline=0ex]
        \draw (-1,-0.5) -- (1,-0.5);
        \draw (-1,0.5) -- (1,0.5);
        \draw (-0.5,0.5) -- (-0.5,-0.5);
        \draw (0.5,0.5) -- (0.5,-0.5);
        \draw [fill] (-0.5,-0.5) circle [radius=0.2];
        \draw [fill=white] (-0.5,0.5) circle [radius=0.2];
        \draw [fill=white] (0.5,-0.5) circle [radius=0.2];
        \draw [fill] (0.5,0.5) circle [radius=0.2];
        \end{tikzpicture} \quad \longleftrightarrow \quad \begin{tikzpicture}[baseline=0ex]
        \draw (-1,-0.5) -- (1,-0.5);
        \draw (-1,0.5) -- (1,0.5);
        \draw (-0.5,0.5) -- (-0.5,-0.5);
        \draw (0.5,0.5) -- (0.5,-0.5);
        \draw [fill=white] (-0.5,-0.5) circle [radius=0.2];
        \draw [fill] (-0.5,0.5) circle [radius=0.2];
        \draw [fill] (0.5,-0.5) circle [radius=0.2];
        \draw [fill=white] (0.5,0.5) circle [radius=0.2];
        \end{tikzpicture} & \quad \quad \left\{\dots, i,-i,\dots\right\}\longleftrightarrow \left\{\dots, -i,i,\dots\right\} \\
        & \\
        \begin{tikzpicture}[baseline=0ex]
        \draw (-1,-0.5) -- (1,-0.5);
        \draw (-1,0) -- (1,0);
        \draw (-1,0.5) -- (1,0.5);
        \draw (-0.5,-0.5) -- (-0.5,0);
        \draw (0,0) -- (0,0.5);
        \draw (0.5,-0.5) -- (0.5,0);
        \draw [fill=white] (0,0.5) circle [radius=0.2];
        \draw [fill=white] (0.5,0) circle [radius=0.2];
        \draw [fill=white] (-0.5,0) circle [radius=0.2];
        \draw [fill] (0,0) circle [radius=0.2];
        \draw [fill] (-0.5,-0.5) circle [radius=0.2];
        \draw [fill] (0.5,-0.5) circle [radius=0.2];
        \end{tikzpicture} \quad \longleftrightarrow \quad \begin{tikzpicture}[baseline=0ex]
        \draw (-1,-0.5) -- (1,-0.5);
        \draw (-1,0) -- (1,0);
        \draw (-1,0.5) -- (1,0.5);
        \draw (-0.5,0.5) -- (-0.5,0);
        \draw (0,0) -- (0,-0.5);
        \draw (0.5,0.5) -- (0.5,0);
        \draw [fill=white] (0,0) circle [radius=0.2];
        \draw [fill=white] (-0.5,0.5) circle [radius=0.2];
        \draw [fill=white] (0.5,0.5) circle [radius=0.2];
        \draw [fill] (0,-0.5) circle [radius=0.2];
        \draw [fill] (-0.5,0) circle [radius=0.2];
        \draw [fill] (0.5,0) circle [radius=0.2];
        \end{tikzpicture}
        & \quad \quad \left\{\dots, i+1,i,i+1,\dots\right\}\longleftrightarrow \left\{\dots, i,i+1, i, \dots\right\}\\
        & \\
         \begin{tikzpicture}[baseline=0ex]
        \draw (-1,-0.5) -- (1,-0.5);
        \draw (-1,0) -- (1,0);
        \draw (-1,0.5) -- (1,0.5);
        \draw (-0.5,-0.5) -- (-0.5,0);
        \draw (0,0) -- (0,0.5);
        \draw (0.5,-0.5) -- (0.5,0);
        \draw [fill] (0,0.5) circle [radius=0.2];
        \draw [fill] (0.5,0) circle [radius=0.2];
        \draw [fill] (-0.5,0) circle [radius=0.2];
        \draw [fill=white] (0,0) circle [radius=0.2];
        \draw [fill=white] (-0.5,-0.5) circle [radius=0.2];
        \draw [fill=white] (0.5,-0.5) circle [radius=0.2];
        \end{tikzpicture} \quad \longleftrightarrow \quad \begin{tikzpicture}[baseline=0ex]
        \draw (-1,-0.5) -- (1,-0.5);
        \draw (-1,0) -- (1,0);
        \draw (-1,0.5) -- (1,0.5);
        \draw (-0.5,0.5) -- (-0.5,0);
        \draw (0,0) -- (0,-0.5);
        \draw (0.5,0.5) -- (0.5,0);
        \draw [fill] (0,0) circle [radius=0.2];
        \draw [fill] (-0.5,0.5) circle [radius=0.2];
        \draw [fill] (0.5,0.5) circle [radius=0.2];
        \draw [fill=white] (0,-0.5) circle [radius=0.2];
        \draw [fill=white] (-0.5,0) circle [radius=0.2];
        \draw [fill=white] (0.5,0) circle [radius=0.2];
        \end{tikzpicture}
        & \quad \quad \begin{array}{l} \left\{\dots, -(i+1),-i,-(i+1),\dots\right\}\longleftrightarrow \\
        \left\{\dots, -i,-(i+1), -i, \dots\right\}
        \end{array}
    \end{array}
    \]
Note that for any two reduced words $\vec{i}$ and $\vec{i}'$ of the same pair of Weyl group elements $(u,v)$, we can transform $\Gamma_\vec{i}$ and $\Gamma_{\vec{i}'}$ into one another via a sequence of moves like the ones above and moves that exchange letters $i$ with $-j$ for $i\neq j$, which can be seen as a merging-splitting move \eqref{splitting}. Therefore Thurston's theorem still holds for bipartite graphs $\Gamma_\vec{i}$.

By following the recipe of quiver drawing (Definition \ref{quiver}), we can get a quiver $\vec{i}$ (this is an abuse of notation on purpose) associated to the bipartite graph $\Gamma_\vec{i}$, which is the same as the quiver one will get from the reduced word $\vec{i}$ by following the amalgamation procedure described in \cite{FGamalgamation}. In particular, the 2-by-2 moves listed above correspond to quiver mutations at the center face in the same way as the case of full-ranked bipartite graphs. After we obtain quivers related by mutations, we can use them to construct cluster varieties; we denote the resulting cluster varieties by $\mathcal{A}^{u,v}:=\mathcal{A}_{\GL_n}^{u,v}$ and $\mathcal{X}^{u,v}:=\mathcal{X}_{\GL_n}^{u,v}$; by removing the frozen vertices we also get an unfrozen version $\left(\mathcal{X}^{u,v}\right)^\uf$.

Next let's move to the level of geometric spaces and cluster varieties. The key to relate double Bruhat cells of $\GL_n$ and decorated Grassmannian $\dGr_n\left(\mathbb{C}^{2n}\right)$ is to think of an element $x$ of the former, which is an $n\times n$ matrix, as an element of the latter represented by the matrix $\left(\overline{w}_0, x\right)$ (note that there is no superscript $\times$ on the decorated Grassmannian because some consecutive minors may vanish). Let's now follow our recipe for $\psi$ and $\chi$ in the Grassmannian story and see how we recover the maps $\psi:\GL_n^{u,v}\dashrightarrow \mathcal{A}^{u,v}_G$ and $\chi:\left(\mathcal{X}^{u,v}_G\right)^\uf\dashrightarrow \GL_n^{u,v}$. 

First for $\psi$, we recall that the map $\psi_\Gamma$ defined by a bipartite graph $\Gamma$ basically maps a decorated Grassmannian, which is represented by a matrix, into the minors corresponding to the dominating sets associated to faces of $\Gamma$. In the case of the bipartite graph $\Gamma$, suppose we define dominating sets $I(f)$ the same was as before (Definition \ref{dominating set}); then for each face $f$ we should assign the minor of the matrix $\left(\overline{w}_0,x\right)$ picked out by elements of the dominating set $I(f)$. But we have seen in the proof of Proposition \ref{positive} that, if we let $I^s(f):=\{1,\dots, m\}\setminus I(f)$ and $I^t(f):=I(f)\setminus \{1,\dots, m\}$, then the minor $\Delta_{I(f)}\left(\overline{w}_0,x\right)=\Delta_{I^s,I^t}(x)$, which is the minor of $x$ whose rows are picked by elements of $I^s$ and whose columns are picked by elements of $I^t$. But then we also observe that $I^s$ precisely consists of indices of the zig-zag strands going from left to right that is above $f$, and $I^t$ precisely consists of indices of the zig-zag strands going from right to left that is above $f$! Therefore we see that the faces in the $i$th spacing corresponds to $i\times i$ minors of a $\GL_n^{u,v}$ matrix $x$, and to read off the rows (resp. columns) of such minor one only need to look at which zig-zag strands going from left to right (resp. right to left) are above the face. In particular, the top face corresponds to the minor $\Delta_{\emptyset, \emptyset}=1$ and the bottom face corresponds to the minor $\Delta_{[1,n],[1,n]}=\det$. One can verify that the minors obtained this way are exactly the same as those decribed by Berenstein, Fomin, and Zelevinsky in \cite{BFZ}. From this we see that the map $\psi:\GL_n^{u,v}\dashrightarrow \mathcal{A}^{u,v}$ arise naturally in the same way as the map $\psi:\dGr_{m,n}^\times\dashrightarrow \mathcal{A}_{m,n}$ we introduced in Section \ref{section4}.

Next for $\chi$, we recall that the map $\chi_\Gamma$ defined by a bipartite graph $\Gamma$ basically maps a generic point $\left(\left(X_f\right)_f\right)$ to a matrix $\left(\overline{w}_0,\left(m_{ij}\right)\right)$ where $m_{ij}:=\sum_{\gamma:i\rightarrow j}\prod_{f\in \hat{\gamma}}X_f$. Since we would like to think of a $\GL_n^{u,v}$ matrix $x$ as $\left(\overline{w}_0,x\right)$, we just need to show that the matrix $\left(m_{ij}\right)=\chi_\vec{i}\left(\left(X_f\right)_f\right)$ is precisely the one that people get when follow the amalgamation map described in \cite{FGamalgamation}. Let's start by observing that the standard perfect orientation on $\Gamma_\vec{i}$ always travel from left to right on the horizontal lines and going from the white vertex to the black vertex on each vertical edge; therefore all sources of the standard perfect orientation on $\Gamma_\vec{i}$ are on the left edge and all sources of the standard perfect orientation on $\Gamma_\vec{i}$ are on the right edge. Now we can assign the standard basis vector $b_i$ of $\mathbb{C}^n$ to boundary marked point $i$ on the left and assign the $j$th column vector $\xi_j$ of $\chi_\vec{i}\left(\left(X_f\right)_f\right)$ to the boundary marked point $j'$ on the right, and view $\chi_\vec{i}\left(\left(X_f\right)_f\right)$ as a matrix that tries to express $\left\{\xi_j\right\}$ in terms of the basis $\left\{b_i\right\}$. Now let's proceed with an induction on the length of the reduced word $\vec{i}$. When $\vec{i}$ is the empty word, it is obvious that $\chi_\vec{i}\left(\left(X_f\right)_f\right)$ agrees with the amalgamation map, which says that
\[
\begin{tikzpicture}[baseline=9ex]
\draw (0,0) node [left] {$n$} -- (1,0);
\draw (0,1) node [left] {$i+1$} -- (1,1);
\draw (0,2) node [left] {$i$} -- (1,2);
\draw (0,3) node [left] {$1$} -- (1,3);
\node at (0.5,1.5) [] {$X_i$};
\node at (0.5,0.5) [] {$\vdots$};
\node at (0.5,2.5) [] {$\vdots$};
\end{tikzpicture}  \quad= \quad  \begin{pmatrix} X_1 & 0 & 0 & \cdots & 0 \\
0 & X_1X_2 & 0  & \cdots & 0 \\
0 & 0 & X_1X_2X_3 & \cdots & 0 \\
\vdots & \vdots & \vdots & \ddots & \vdots\\
0 & 0 & 0 & \cdots & \prod_i X_i\end{pmatrix}.
\]
Suppose then we add a new letter $i$ to the reduced word $\vec{i}$ to get a longer reduced word $\vec{i}'$; then from the construction of $\chi_{\vec{i}'}$ using the standard perfect orientation we can express the collection $\left\{\xi'_i\right\}$ in terms of $\left\{\xi_i\right\}$, which should be
\[
\begin{tikzpicture}[baseline=9ex]
\draw[blue, ->] (0,3) -- (2,3);
\draw[blue, ->] (0,2) -- (2,2);
\draw[blue, ->] (0,1) -- (2,1);
\draw[blue, ->] (0,0) -- (2,0);
\draw[->-, blue] (1,2) -- (1,1);
\draw[fill=white] (1,2) circle [radius=0.2];
\draw[fill=black] (1,1) circle [radius=0.2];
\node at (1,2.5) [] {$\vdots$};
\node at (1,0.5) [] {$\vdots$};
\node at (-0.2,0) [left] {$\xi_n$};
\node at (-0.2,1) [left] {$\xi_{i+1}$};
\node at (-0.2,2) [left] {$\xi_{i}$};
\node at (-0.2,3) [left] {$\xi_1$};
\node at (2.2,0) [right] {$\xi'_n$};
\node at (2.2,1) [right] {$\xi'_{i+1}$};
\node at (2.2, 2) [right] {$\xi'_{i}$};
\node at (2.2,3) [right] {$\xi'_1$};
\node at (1.5,1.5) [] {$X$};
\end{tikzpicture} \quad \quad \quad \quad \quad \quad \left(\xi_1,\dots, \xi_n\right)\begin{pmatrix}X &  \cdots & 0 & 0 & \cdots & 0 \\
\vdots & \ddots & \vdots & \vdots & \ddots & \vdots \\
0 & \cdots & X & 1 & \cdots & 0 \\
0 & \cdots & 0 & 1 & \cdots & 0 \\
\vdots & \ddots & \vdots & \vdots & \ddots & \vdots \\
0 & \cdots & 0 & 0 & \cdots & 1 \end{pmatrix}=\left(\xi'_1,\dots, \xi'_n\right). 
\]
Note that the big square matrix above is equal to $e_{i}\diag(\underbrace{X,\dots, X}_i, 
\underbrace{1, \dots, 1}_{n-i})$. Similarly, if we add a new letter $-i$ instead, we get
\[
\begin{tikzpicture}[baseline=9ex]
\draw[blue, ->] (0,3) -- (2,3);
\draw[blue, ->] (0,2) -- (2,2);
\draw[blue, ->] (0,1) -- (2,1);
\draw[blue, ->] (0,0) -- (2,0);
\draw[->-, blue] (1,1) -- (1,2);
\draw[fill] (1,2) circle [radius=0.2];
\draw[fill=white] (1,1) circle [radius=0.2];
\node at (1,2.5) [] {$\vdots$};
\node at (1,0.5) [] {$\vdots$};
\node at (-0.2,0) [left] {$\xi_n$};
\node at (-0.2,1) [left] {$\xi_{i+1}$};
\node at (-0.2,2) [left] {$\xi_{i}$};
\node at (-0.2,3) [left] {$\xi_1$};
\node at (2.2,0) [right] {$\xi'_n$};
\node at (2.2,1) [right] {$\xi'_{i+1}$};
\node at (2.2, 2) [right] {$\xi'_{i}$};
\node at (2.2,3) [right] {$\xi'_1$};
\node at (1.5,1.5) [] {$X$};
\end{tikzpicture} \quad \quad \quad \quad \quad \quad \left(\xi_1,\dots, \xi_n\right)\begin{pmatrix}X &  \cdots & 0 & 0 & \cdots & 0 \\
\vdots & \ddots & \vdots & \vdots & \ddots & \vdots \\
0 & \cdots & X & 0 & \cdots & 0 \\
0 & \cdots & X & 1 & \cdots & 0 \\
\vdots & \ddots & \vdots & \vdots & \ddots & \vdots \\
0 & \cdots & 0 & 0 & \cdots & 1 \end{pmatrix}=\left(\xi'_1,\dots, \xi'_n\right). 
\]
where the big square matrix is $e_{-i}\diag(\underbrace{X,\dots, X}_i, 
\underbrace{1, \dots, 1}_{n-i})$. Now we see that the matrix $\chi_\vec{i}\left(\left(X_f\right)_f\right)$ constructed by using the standard perfect orientation can be obtained by first decomposing $\Gamma_\vec{i}$ into parts that are of one of the following three forms, and then multiplying the corresponding group elements from left to right: 
\[
\begin{tikzpicture}[baseline=9ex]
\draw (0,0) node [left] {$n$} -- (1,0);
\draw (0,1) node [left] {$i+1$} -- (1,1);
\draw (0,2) node [left] {$i$} -- (1,2);
\draw (0,3) node [left] {$1$} -- (1,3);
\node at (0.5,0.5) [] {$\vdots$};
\node at (0.5,2.5) [] {$\vdots$};
\draw (0.5,1) -- (0.5,2);
\draw [fill=white] (0.5,2) circle [radius=0.2];
\draw [fill] (0.5,1) circle [radius=0.2];
\end{tikzpicture}\quad= \quad e_i, \quad \quad 
\begin{tikzpicture}[baseline=8ex]
\draw (0,0) node [left] {$n$} -- (1,0);
\draw (0,1) node [left] {$i+1$} -- (1,1);
\draw (0,2) node [left] {$i$} -- (1,2);
\draw (0,3) node [left] {$1$} -- (1,3);
\node at (0.5,0.5) [] {$\vdots$};
\node at (0.5,2.5) [] {$\vdots$};
\draw (0.5,1) -- (0.5,2);
\draw [fill] (0.5,2) circle [radius=0.2];
\draw [fill=white] (0.5,1) circle [radius=0.2];
\end{tikzpicture} \quad = \quad  e_{-i}, \quad  \quad 
\begin{tikzpicture}[baseline=8ex]
\draw (0,0) node [left] {$n$} -- (1,0);
\draw (0,1) node [left] {$i+1$} -- (1,1);
\draw (0,2) node [left] {$i$} -- (1,2);
\draw (0,3) node [left] {$1$} -- (1,3);
\node at (0.5,1.5) [] {$X$};
\node at (0.5,0.5) [] {$\vdots$};
\node at (0.5,2.5) [] {$\vdots$};
\end{tikzpicture}  \quad= \quad  \diag(\underbrace{X,\dots, X}_i, 
\underbrace{1, \dots, 1}_{n-i}).
\]
But this is exactly the amalgamation map for $\GL_n^{u,v}$ described in \cite{FGamalgamation}! Now to get a map from $\left(\mathcal{X}^{u,v}\right)^\uf$ to $H\left\backslash \GL_n^{u,v}\right/H$ we just need to notice that the boundary faces on the left correspond to a left multiplication by $H$ and the boundary faces on the right correspond to a right multiplication by $H$. Therefore we can conclude that the amalgamation map $\chi:\left(\mathcal{X}^{u,v}\right)^\uf\dashrightarrow H\left\backslash \GL_n^{u,v}\right/H$ arise naturally in the same way as the map $\chi:\mathcal{X}_{m,n}^\uf\dashrightarrow \conf_n^\times\left(\mathbb{P}^{m-1}\right)$ introduced at the end of Section \ref{quantum}.

\begin{exmp} Let's end by computing $\psi_\vec{i}$ and $\chi_\vec{i}$ for the reduced word $\vec{i}=(1,-2,2,1,-1,-2)$ used in Example \ref{12.1}. By drawing the zig-zag strands one can find that the minors defined by the dominating sets of the faces are as follows, and the map $\psi_\vec{i}$ simply sends a $\GL_3^{w_0,w_0}$ element to the 10-tuple of minors shown in the picture (note that $\Delta_{\emptyset,\emptyset}=1$ by convention).
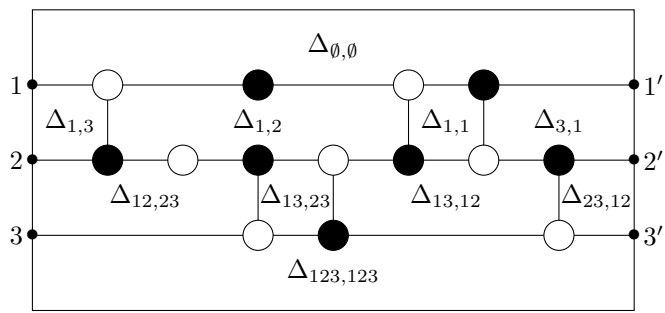
\[
\tikz{
\draw (-1,-1) rectangle (7,3);
    \draw (-1,0) -- (7,0);
    \draw (-1,1) -- (7,1);
    \draw (-1,2) -- (7,2);
    \draw (0,1) -- (0,2);
    \draw (2,1) -- (2,0);
    \draw (3,1) -- (3,0);
    \draw (4,1) -- (4,2);
    \draw (5,1) -- (5,2);
    \draw (6,1) -- (6,0);
    \draw[fill=white] (0,2) circle [radius=0.2];
    \draw[fill=black] (0,1) circle [radius=0.2];
    \draw[fill=white] (2,0) circle [radius=0.2];
    \draw[fill=black] (2,1) circle [radius=0.2];
    \draw[fill=white] (3,1) circle [radius=0.2];
    \draw[fill=black] (3,0) circle [radius=0.2];
    \draw[fill=white] (4,2) circle [radius=0.2];
    \draw[fill=black] (4,1) circle [radius=0.2];
    \draw[fill=white] (5,1) circle [radius=0.2];
    \draw[fill=black] (5,2) circle [radius=0.2];
    \draw[fill=white] (6,0) circle [radius=0.2];
    \draw[fill=black] (6,1) circle [radius=0.2];
    \draw[fill=white] (1,1) circle [radius=0.2];
    \draw[fill=black] (2,2) circle [radius=0.2];
    \foreach \i in {1,2,3}
        {
        \node at (-1,3-\i) [] {$\bullet$};
        \node at (7,3-\i) [] {$\bullet$};
        \node at (-1,3-\i) [left] {$\i$};
        \node at (7,3-\i) [right] {$\i'$};
        }
    \node at (3,2.5) [] {$\Delta_{\emptyset,\emptyset}$};
    \node at (-0.5,1.5) [] {$\Delta_{1,3}$};
    \node at (2,1.5) [] {$\Delta_{1,2}$};
    \node at (4.5,1.5) [] {$\Delta_{1,1}$};
    \node at (6,1.5) [] {$\Delta_{3,1}$};
    \node at (0.5,0.5) [] {$\Delta_{12,23}$};
    \node at (2.5,0.5) [] {$\Delta_{13,23}$};
    \node at (4.5,0.5) [] {$\Delta_{13,12}$};
    \node at (6.5,0.5) [] {$\Delta_{23,12}$};
    \node at (3,-0.5) [] {$\Delta_{123,123}$};
}
\]
On the other hand, by assigning cluster $\mathcal{X}$-variables to the non-boundary faces of $\Gamma_\vec{i}$ as follows
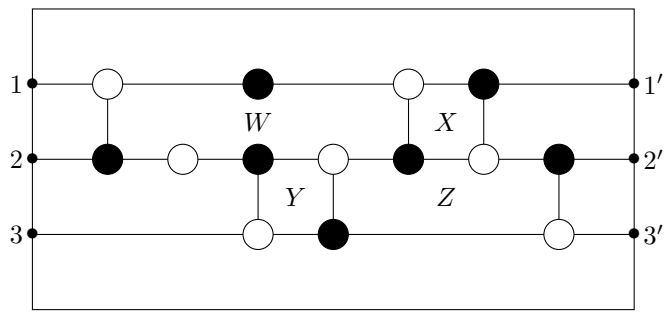
\[
\tikz{
\draw (-1,-1) rectangle (7,3);
    \draw (-1,0) -- (7,0);
    \draw (-1,1) -- (7,1);
    \draw (-1,2) -- (7,2);
    \draw (0,1) -- (0,2);
    \draw (2,1) -- (2,0);
    \draw (3,1) -- (3,0);
    \draw (4,1) -- (4,2);
    \draw (5,1) -- (5,2);
    \draw (6,1) -- (6,0);
    \draw[fill=white] (0,2) circle [radius=0.2];
    \draw[fill=black] (0,1) circle [radius=0.2];
    \draw[fill=white] (2,0) circle [radius=0.2];
    \draw[fill=black] (2,1) circle [radius=0.2];
    \draw[fill=white] (3,1) circle [radius=0.2];
    \draw[fill=black] (3,0) circle [radius=0.2];
    \draw[fill=white] (4,2) circle [radius=0.2];
    \draw[fill=black] (4,1) circle [radius=0.2];
    \draw[fill=white] (5,1) circle [radius=0.2];
    \draw[fill=black] (5,2) circle [radius=0.2];
    \draw[fill=white] (6,0) circle [radius=0.2];
    \draw[fill=black] (6,1) circle [radius=0.2];
    \draw[fill=white] (1,1) circle [radius=0.2];
    \draw[fill=black] (2,2) circle [radius=0.2];
    \foreach \i in {1,2,3}
        {
        \node at (-1,3-\i) [] {$\bullet$};
        \node at (7,3-\i) [] {$\bullet$};
        \node at (-1,3-\i) [left] {$\i$};
        \node at (7,3-\i) [right] {$\i'$};
        }
    \node at (2,1.5) [] {$W$};
    \node at (4.5,1.5) [] {$X$};
    \node at (2.5,0.5) [] {$Y$};
    \node at (4.5,0.5) [] {$Z$};
}
\]
we can compute the image $\chi_\vec{i}$ in $H\left\backslash \GL_3^{w_0,w_0}\right/H$, which is represented by the matrix
\[
e_1e_{-2}\diag(WY,Y,1)e_2e_1\diag(XZ,Z,1)e_{-1}e_{-2}=\begin{pmatrix} (1+W+WX)YZ & Y(1+Z+WZ) & Y \\
YZ & Y(1+Z) & Y \\
YZ & 1+Y+YZ & 1+Y
\end{pmatrix}.
\]
\end{exmp}

\begin{rmk} We have seen that bipartite graphs on a disk are closely tied to Grassmannian $\Gr_{m,n}^\times$ and double Bruhat cells $\GL_n^{u,v}$, and these two cases have $\GL_n^{w_0,w_0}$ as an intersection. In fact, there is a bigger family of geometric spaces called \emph{positroid cells} of Grassmannian that can be studied by using bipartite graphs on a disk. Postnikov has written down the combinatorics as well as the map $\chi$ in \cite{Pos}. We believe that a commutative diagram similar to \eqref{summary} should exist for these positroid cells, and cluster Donaldson-Thomas transformations should also exist in these cases, but more studies need to be done.
\[
\tikz{
\draw (0,0.25) ellipse (6 and 2);
\draw (-2,0.5) ellipse (3 and 1);
\draw (2,0.5) ellipse (3 and 1);
\node at (0,0.5) [] {$\GL_n^{w_0,w_0}$};
\node at (-3,0.5) [] {$\Gr_{m,n}^\times$};
\node at (3,0.5) [] {$\GL_n^{u,v}$};
\node at (0,-1) [] {Positroid Cells};
}
\]
\end{rmk}

\bibliographystyle{amsalpha-a}

\bibliography{biblio}

\end{document}